\documentclass[11pt]{amsart}
\usepackage{amssymb}
\usepackage{verbatim}
\usepackage{amscd}
\usepackage{amsmath}
\usepackage{graphics,graphicx,psfrag,caption}
\usepackage[toc,page]{appendix}

\author{}
\thanks{}
 \keywords{}

\newtheorem{theorem}{Theorem}[section]

\newtheorem{lemma}[theorem]{Lemma}
\newtheorem{proposition}[theorem]{Proposition}

\theoremstyle{definition}
\newtheorem{remark}[theorem]{Remark}

\numberwithin{equation}{section}

\def\boxit#1{\vbox{\hrule\hbox{\vrule\kern3pt
     \vbox{\kern3pt#1\kern3pt}\kern3pt\vrule}\hrule}}

\newfont{\msam}{msam10}            
\newfont{\msym}{msbm10 scaled\magstep1}            

\newfont{\gotic}{eufm10 scaled\magstep1}

\newcommand{\ra}{\rightarrow}
\newcommand{\lra}{\mbox{\Huge $\longrightarrow$}}

\newcommand{\kra}{\kern-7pt\rightarrow\kern-7pt}
\newcommand{\na}{NA}
\newcommand{\sono}{SO_0(n,1)}
\newcommand{\frakk}{\mathfrak{k}}
\newcommand{\fsono}{\mathfrak{so}(n,1)}

\def\bbh{{\mathbb H}}

\def\ra{\rightarrow}

\def\inv{^{-1}}

\def\lra{\longrightarrow}

\numberwithin{equation}{section}

\newcommand{\angles}[1]{{\langle #1 \rangle}}
\def\begtab{\begin{tabbing} WW\=23/02: \= point 1\kill}

\def\NIL1{{\mathcal H^{3}}}

\def\so{\mathrm{SO}}

\def\n2{\mathfrak{N}_2}

\def\frakg{{\mathfrak g}}

\begin{document}

\title{Horizontal displacement of curves in bundle $\text{SO}(n)
\rightarrow \text{SO}_0 (1,n) \rightarrow \mathbb{H}^n$}

\author{Taechang Byun}
\address{Dept.\ of Mathematics \\
         University of Oklahoma \\
         Norman, OK, 73019 }
\email{tcbyun@math.ou.edu}

\begin{abstract}

  The Riemannian submersion
  $ \pi : \text{SO}_0(1,n) \rightarrow \mathbb{H}^n $ is
  a principal bundle and its fiber at $ \pi (e) $ is the imbedding
  of $\text{SO}(n)$ into $ \text{SO}_0(1,n) $ , where $e$ is the
  identity of both $\text{SO}_0(1,n)$ and $\text{SO}(n)$. In this study, we
  associate a curve, starting from the identity, in $\text{SO}(n)$ to
  a given surface with boundary, diffeomorphic to the closed disk
  $D^2$, in $ \mathbb{H}^n $ such that the starting point and the
  ending point of the curve agree with those of the horizontal lifting
  of the boundary curve of the given surface with boundary,
  respectively, and that the length of the curve is as same as the
  area of the given surface with boundary.

\end{abstract}

\maketitle

\section*{Introduction}

  \medskip
  Let
  $
    \text{O}(1,n) =
    \{A \in \text{GL}(n+1; \mathbb{R}) \; | \;A^{t}SA = S \},
  $
  where
  $ S =
    \left(
            \begin{matrix}
                 -1 & 0             \\
                  0 & \textbf{I}_n
            \end{matrix}
    \right) .
  $

  Let $\text{SO}_0(1,n)$ be the identity component of
  $\text{O}(1,n)$, which is also the identity component of
  $\text{SO}(1,n)$, and consider a subgroup of $\text{SO}_0(1,n)$
  consisting of all matrices of the form
  $
    \left(
            \begin{matrix}
                1 & 0       \\
                0 & B
            \end{matrix}
    \right),
  $
  where $B \in \text{SO}(n).$
  Call the embedded subgroup $\text{SO}(n)$ again.

  Note the Lie algebra $\mathfrak{o}(1,n)$ is given by

  $$
  \mathfrak{o}(1,n) =
  \{ X \in \mathfrak{gl}(n+1;\mathbb{R}) | X^t S + SX = 0 \} .
  $$

  Now, think of a left-invariant metric on $\text{SO}_0 (1,n)$,
  induced from an inner product $\langle \cdot \; , \; \cdot \rangle$
  on the Lie algebra, $\mathfrak{so}(1,n)$, defined as follows :

  $$
   \langle A,B \rangle = \frac{1}{2} \text{ trace}(A^t \, B)
   \qquad \text{for } A,B \in \mathfrak{so}(1,n).
  $$

  If $\phi$ is a Killing-Cartan form, then
   $$
   \langle X,Y \rangle  = \: \frac{n-1}{2} \: \phi(X,Y)
   \hspace{1cm} \text{ for } X,Y \in \mathfrak{so}(n)^{\perp} \subset
   \mathfrak{so}(1,n) \: .
   $$

  And the right actions of $\text{SO}(n)$ become isometries and
  $ \text{SO}_0(1,n) / \text{SO}(n)$ becomes isometric to
  $\mathbb{H}^n.$

  Under this metric, we have a principal bundle structure
  $$
  \text{SO}(n) \rightarrow \text{SO}_0 (1,n) \rightarrow \mathbb{H}^n
  \: ,
  $$
  where $\pi : \text{SO}_0 (1,n) \rightarrow \mathbb{H}^n $ is a
  Riemannian submersion.

  If n=2, then it can be easily shown that for a given geodesic
  triangle in $\mathbb{H}^2$, the distance by the horizontal
  displacement of the boundary curve of the given geodesic triangle
  in the fiber is as same as the area of the triangle. Furthermore,
  the direction of the boundary curve of the given geodesic triangle
  in $\mathbb{H}^2$ will determine the direction of its holonomy
  displacement. All of these are dealt with in Section ~\ref{sec:n=2}.

  Can a similar result be obtained in a topological disk in
  $\mathbb{H}^n$?

  If it is a geodesic triangle, something similar can be easily said
  from the result for the case n=2 and `Fact 2', mentioned in
  section \ref{sec:n-general}. But, what can be done for a general disk in
  $\mathbb{H}^n$?

  To answer this question, we intend to approximate the given disk with geodesic
  triangles, since there exists a unique totally geodesic triangle
  for any 3 different points in $\mathbb{H}^n.$ And then we intend to
  construct a curve in the fiber by using the property for the case
  n=2. But how can we approximate it? Though each geodesic triangle
  and its boundary curve determine the direction of each horizontal
  displacement, some linear ordering of geodesic triangles and the induced
  ordering of their boundary curves may not represent the boundary curve of
  their union. If the given disk is contained in an isometrically
  embedded plane $\mathbb{H}^2$ in $\mathbb{H}^n \; ,$ something
  similar can be said from a curve in the fiber SO(n), made from the
  result for the case n=2 and `Fact 2', mentioned in section
  \ref{sec:n-general},
  since horizontal displacements are happening in the one-dimensional
  vertical subgroup. Though the different orderings of triangles give
  different curves in the vertical space, they will meet at the same
  ending point. So, with respect to any ordering,  the horizontal
  displacement of the boundary curve of the given disk can be
  approximated. But in other cases, what can be obtained? Something
  similar could be done if the fiber SO(n) were abelian, which would
  make the ending points of any other different two curves in the
  fiber, induced from different linear orderings, be the same. But the
  fiber SO(n) is not abelian for $n \geq 3$. The difficult part is that not only the
  approximation of the area but also the linear ordering of the
  triangles on each step for the approximation of the
  boundary curve of the disk should be considered at the same time.
  This is one of the hardest parts in this paper, which is dealt with
  in Section ~\ref{sec:strategy} and Appendices ~\ref{sec:triangles}
  and ~\ref{sec:curves}. Furthermore, can holonomy displacements by the
  lifts of piecewise geodesics approaching to the boundary of the given
  topological disk in the base space converge to the holonomy displacement
  by the lift of the boundary? It will be discussed in
  Subsection \ref{sec-converge-lift}.

  After the case n=2 is explained in section ~\ref{sec:n=2}, our
  following main result for the general case will be explained in
  section ~\ref{sec:n-general}.

 \begin{theorem} \label{thm}
  Let $ \pi : \text{SO}_0(1,n) \rightarrow \mathbb{H}^n $ be the
  Riemannian submersion given as before.
       \label{lenth-area-thm} Then, given a smooth disk $S$, with
         $\bar{e} = \pi(e)$ on its boundary, in $\mathbb{H}^n$,
         there is a $C^1$- curve $ f : [0,1] \rightarrow
         \text{SO}(n) \subset \text{SO}_0(1,n)$ with $f(0) = e$
         such that
       \begin{itemize}
         \item $f(1) = f(0)^{-1} f(1) = $ the difference by
             the holonomy induced from the boundary of $S$  in
             view of right multiplication
         \item the length of the curve  $f$ = the area of
             $S$.
       \end{itemize}
 \end{theorem}

 \bigskip

\section{Strategy for approximation}\label{sec:strategy}

  \medskip

  `Factorization Lemma', given by Lichnerowicz,
  \emph{
        Theorie Globale des Connexions et des Groupes
        d'Holonomie
        },
  [3, vol 1, p.284],
  will be helpful to understand this section. For the difference,
  focus on properties of triangles mentioned in number 6. The
  reason for introducing another approximation
  will be given in subsection \ref{factorization}.

  \medskip

  1. For any 3 points in $\mathbb{H}^n ,$ there exists a unique
  totally geodesic triangle with these vertices.

  \medskip

  2. Let $ \bigtriangleup ABC $ be a totally geodesic triangle in
  $\mathbb{H}^n $ and consider a piecewise geodesic from
  $ \bar{e} = \pi (e)$ to $A,$ where $e$ is the identity of SO(1,n).

    \begin{figure}[h]
     \centering{\includegraphics[width=2in]{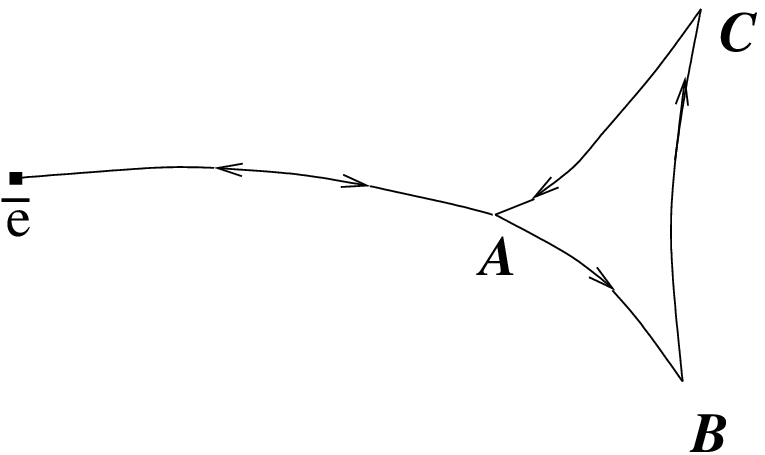}}
    \end{figure}

    Then, horizontal displacement of
    $
     \gamma =
     \overline{\bar{e}A} \cdot \overline{AB} \cdot \overline{BC}
     \cdot \overline{CA} \cdot \overline{A \bar{e}}
    $
    is $g \in $ SO(n), where
    $$
      \text{the length of } \: \overline{eg} \:
      =  \: \text{ the area of } \: \bigtriangleup ABC \: .
    $$

  \medskip

  3. Let $ \bigtriangleup ABC $ and $ \bigtriangleup ACD $ be two
  given geodesic triangles in $\mathbb{H}^n$ and consider a piecewise
  geodesic curve from $\bar{e}$ to $A.$

    \begin{figure}[h]
       \centering{\includegraphics[width=2in]{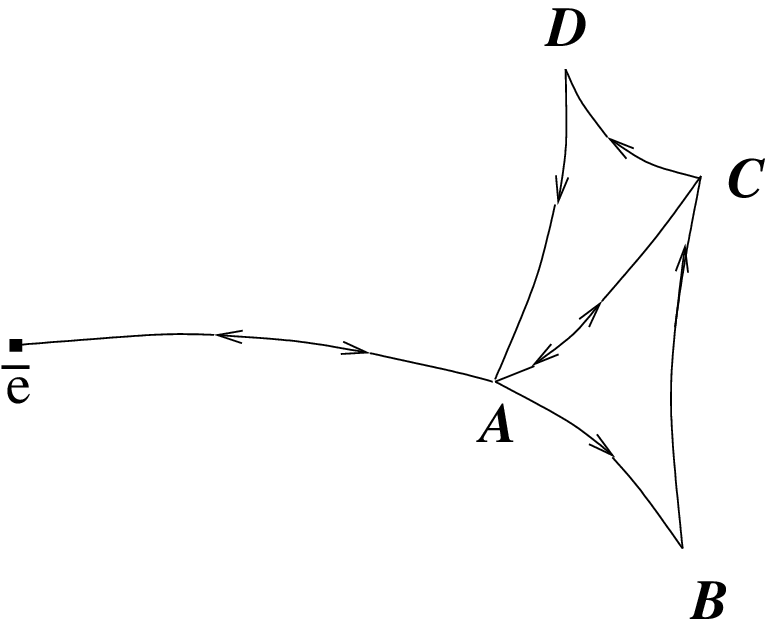}}
    \end{figure}

    Consider two curves
    $
      \gamma_1 =
      \overline{\bar{e}A} \cdot \overline{AB} \cdot \overline{BC}
      \cdot \overline{CA} \cdot \overline{A \bar{e}}
    $
    and
    $
      \gamma_2 =
      \overline{\bar{e}A} \cdot \overline{AC} \cdot \overline{CD}
      \cdot \overline{DA} \cdot \overline{A \bar{e}} \; .
    $
    Then the horizontal displacement of
    ${\gamma_1} * { \gamma_2}$ equals to that of
    $
      \gamma_3 =
      \overline{\bar{e}A} \cdot \overline{AB} \cdot \overline{BC}
      \cdot \overline{CD} \cdot \overline{DA} \cdot
      \overline{A \bar{e}} \; .
    $

    \begin{figure}[h]
      \centering{\includegraphics[width=2.2in]{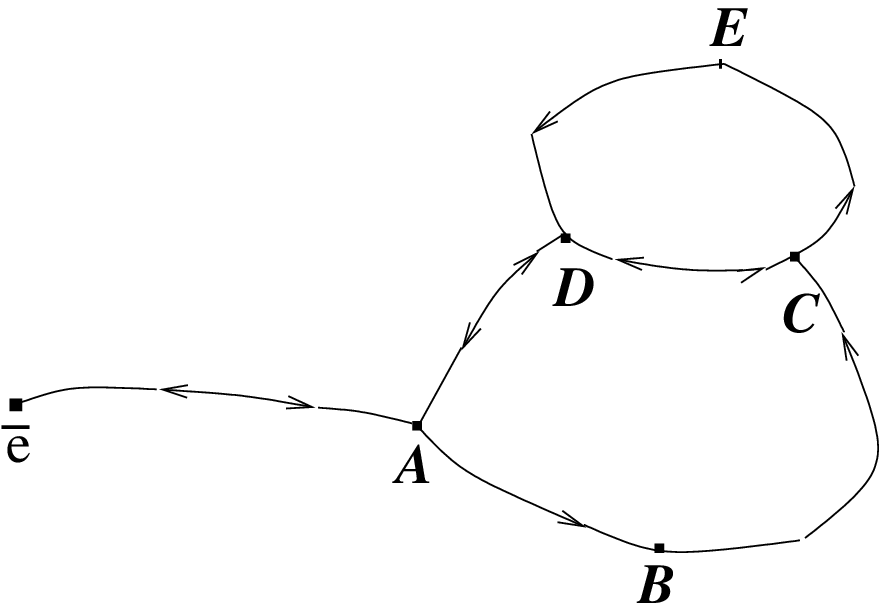}}
    \end{figure}

    In general, if
     $
      \gamma_1 =
      \overline{\bar{e}A} \cdot \overline{AB} \cdot \overline{BC}
      \cdot \overline{CD} \cdot \overline{DA} \cdot \overline{A
      \bar{e}}
     $
     and
     $
      \gamma_2 =
      \overline{\bar{e}A} \cdot \overline{AD} \cdot \overline{DC}
      \cdot \overline{CE} \cdot \overline{ED} \cdot \overline{DA}
      \cdot \overline{A \bar{e}} \;
     $
    are two curves in $\mathbb{H}^n \; ,$
    then the horizontal displacement of ${\gamma_1} * { \gamma_2}$ equals
    to that of
     $
      \gamma_3 = \overline{\bar{e}A} \cdot \overline{AB} \cdot
      \overline{BC} \cdot \overline{CE} \cdot \overline{ED} \cdot
      \overline{DA} \cdot \overline{A \bar{e}} \; .
     $

  \medskip

  4. What's the difficulty of the approximation?

    Consider three given geodesic triangles
    $
    \bigtriangleup ABC \; , \bigtriangleup ACE \; ,
    \bigtriangleup CDE \;
    $
    in $\mathbb{H}^n \; . $

   \begin{figure}[h]
     \centering{\includegraphics[width=2.2in]{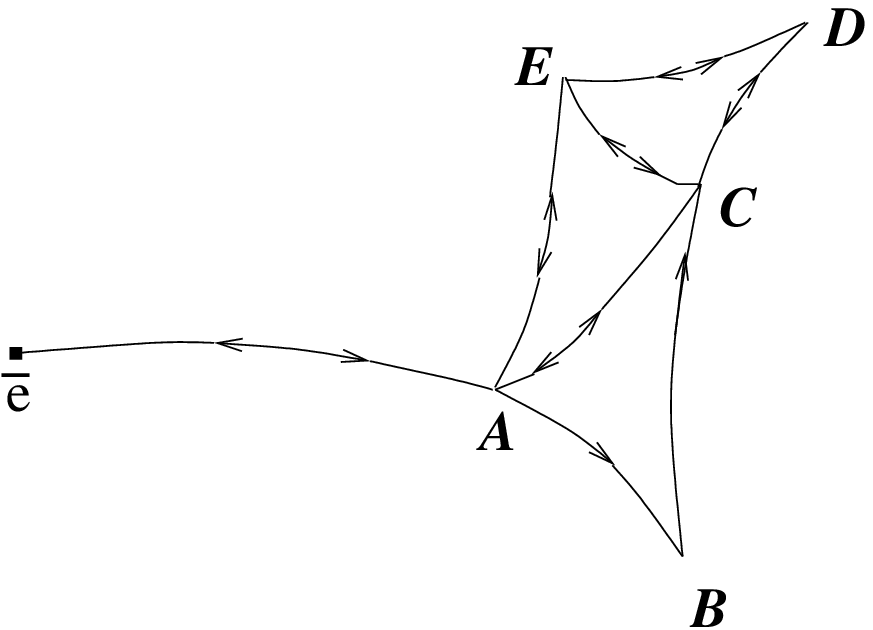}}
   \end{figure}

   Then, for three curves
   $
   \gamma_1 =
   \overline{\bar{e}A} \cdot \overline{AB} \cdot \overline{BC} \cdot
   \overline{CA} \cdot \overline{A \bar{e}} \; , \;
   $
   $ \gamma_2 =
   \overline{\bar{e}A} \cdot \overline{AC} \cdot \overline{CE} \cdot
   \overline{EA} \cdot \overline{A \bar{e}}
   $
   and
   $
   \gamma_3 = \overline{\bar{e}A} \cdot \overline{AE} \cdot
   \overline{EC} \cdot \overline{CD} \cdot \overline{DE} \cdot
   \overline{EA} \cdot \overline{A \bar{e}} \; ,
   $
   the horizontal lift of $\gamma_1 * \gamma_2 * \gamma_3$ equals to that of
   $
   \gamma_4 = \overline{\bar{e}A} \cdot \overline{AB} \cdot
   \overline{BC} \cdot \overline{CD} \cdot \overline{DE} \cdot
   \overline{EA} \cdot \overline{A \bar{e}} \; ,
   $
   which relates to the boundary of the polygon $ABCDE \: .$ But the
   horizontal lift of $\gamma_1 * \gamma_3 * \gamma_2$ equals to
   $
   \gamma_5 =
   \overline{\bar{e}A} \cdot \overline{AB} \cdot \overline{BC} \cdot
   \overline{CA} \cdot \overline{AE} \cdot \overline{EC} \cdot
   \overline{CD} \cdot \overline{DE} \cdot \overline{EA} \cdot
   \overline{AC} \cdot \overline{CE} \cdot \overline{EA} \cdot
   \overline{A \bar{e}}\; ,
   $
   which does not relate to the boundary of the polygon $ABCDE \; .$
   Thus, for our object, the order of curves is important, which
   relates to the order of triangles.

  \medskip

  5. Refer to the number 4.

   Consider a curve
   $
   \tilde{\gamma}_3 =
   \overline{\bar{e}A} \cdot \overline{AC} \cdot \overline{CD} \cdot
   \overline{DE} \cdot \overline{EC} \cdot \overline{CA} \cdot
   \overline{A \bar{e}} \; .
   $
   Though the order
   $ ( \gamma_1 , \; \gamma_2 , \; \tilde{\gamma}_3)$ of curves
   relates to the order of the triangles, induced by the order
   $ ( \gamma_1 , \; \gamma_2 , \; \gamma_3 )$ of curves,
   the horizontal lift of
   $ \gamma_1 * \gamma_2 * \tilde{\gamma}_3$ equals to
   $
    {\tilde{\gamma}_4} =
    \overline{\bar{e}A} \cdot \overline{AB} \cdot \overline{BC} \cdot
    \overline{CE} \cdot \overline{EA} \cdot \overline{AC} \cdot
    \overline{CD} \cdot \overline{DE} \cdot \overline{EC} \cdot
    \overline{CA} \cdot \overline{A \bar{e}} \; ,
   $
   which does not relate to the boundary of the polygon $ABCDE \: .$
   Thus, it is also important how to make a curve that represents a
   given triangle. This problem in the construction of a curve for
   each triangle will be solved by introducing the starting point and
   the ending point of each triangle.

  \medskip

  6. Instead of approximating a given topological disk in
  $\mathbb{H}^n$ directly, we will approximate $D^2$ by triangles,
  and approximate the given disk in $\mathbb{H}^n$ by the diffeomorphism
  from $D^2$ to it. In fact, in Appendix,
  for each $n=0,1,2, \cdots$, we will
  construct a subdivision $D_n$ of the interval $[0,1]$ and
  an ordered set $A_n$ consisting of triangles having
  the following properties:

  \bigskip

  Property 1.) Given a non-first element $L$ in $A_n$, the boundary
  of $ \bigcup \{ M \in A_n | M < L \}  $ contains a side of L, which
  will be divided into two line segments in its barycentric
  subdivision, where one of two line segments will become a side of
  the first triangle and the other one will become a side of the
  second triangle in the barycentric subdivision of L.

  \bigskip

  Property 2.) Given $L \in A_n$ ,
  $ \bigcup \{ M \in A_n | M \leq L \}  $ is diffeomorphic to the
  disk $D^2$.
  \bigskip

  Property 3.)
  Assume $L \in A_n$ and six triangles
  $M_1, M_2, \cdots, M_6 \in A_{n+1}$, obtained from the
  barycentric subdivision of $L$, follows the order
  of $i=1,2, \cdots 6$ in $A_{n+1}$. Then the starting points of
  $M_1$ and $L$ are same. Also are the ending points of $M_6$ and
  $L$.

  \bigskip

  Property 4.)
  Assume $ L, M \in A_n $ and that $M$ is the next element of $L$ in
  $A_n$ for $n \geq 1$.

  Then,  The ending point of $L$ and the starting point of $M$ are
  same.
  \bigskip

\section{Definitions, Triangles and Curves}

  \medskip

  All materials in this section will be dealt with in Appendix
  concretely.

  \bigskip

  \subsection{Notations}

  $f*g : [0,1] \rightarrow \mathbb{H}^n$ is an ordinary
  juxtaposition of curves $f,g : [0,1] \rightarrow \mathbb{H}^n$.
  And, for a given curve $c:[0,1] \rightarrow \mathbb{H}^n$, $\bar{c}$
  represents a curve whose direction is opposite to that of c, that is,
  $\bar{c}:[0,1] \rightarrow \mathbb{H}^n$ is given by
  $\bar{c}(t) = c(1-t).$

  \subsection{Simplification $ \mathbf{\gamma} $ of a curve
  $ \mathbf{g : [a,b] \rightarrow \mathbb{H}^n}$}

  Given a curve $ g : [a,b] \rightarrow S $ , we can think of a
  curve $ \gamma : [a,b] \rightarrow S $ whose direction is
  one-sided as follows :

  If we can find $ c, d, e \in (a,b) $ such that $ a<c<d<e<b $ and
  $ Im (g| _{[c,d]}) = Im (g| _{[d,e]}) $ and that the directions of
  $ g| _{[c,d]} $ and $ g| _{[d,e]} $ are one-sided but opposite
  from each other, then we can think of the new curve
  $ \tilde{g} : [a,b] \rightarrow D^2 $ from the remaining part
  $ g| _{[a,c]} $ and $ g| _{[e,b]} $ by translating in the domain
  and rescaling as follows :

  Note $g(c) = g(e).$

  Consider two curves $g_1 : [a, d] \rightarrow \mathbb{H}^n  $ and
  $g_2 : [d, b] \rightarrow \mathbb{H}^n $ given by
  $$
    \displaystyle
    g \left( \frac{c-a}{d-a} (t-a) + a \right)
    = g_1 (t)  \text{ for } t \in [a, d]
  $$
  and
  $$
    \displaystyle
    g \left( \frac{b-e}{b-d} (t-b) + b \right)
    = g_2 (t)  \text{ for } t \in [d, b],
  $$
  and then let $\tilde{g} = g_1 * g_2.$

  From a curve obtained by doing this work again and again and
  by reparametrizing it, we can
  think of a constant speed curve
  $ \gamma : [a,b] \rightarrow S $ which we want.

  \subsection{The definition of $\mathbf{D_n, j_n, t^n_1, t^n_2}$}

  \begin{align*}
        D_n =
        \displaystyle
        \left\{
           \dfrac{1}{2} \cdot \frac{j}{6^n} \mid j
           = 0, 1, 2, \cdots , 6^n \right
        \}
        \bigcup
        \\
        \left(
          \cup ^n_{k=1}
          \left\{
            \sum ^k_{i=1} \frac{1}{2^i}
            + \frac{1}{2^{k+1}} \cdot
            \frac{j}{2^{k-1} \cdot 6^{n-k+1}}
            \mid
            j=0,1,2, \cdots , 2^{k-1} \cdot 6^{n-k+1}
          \right\}
        \right)
  \end{align*}

  Think of the usual order $D_n$ and regard
  $$
    0, \: \frac{1}{2} \cdot \frac{1}{6^n}, \:
    \frac{1}{2} \cdot \frac{2}{6^n}, \:
    \cdots , \:
    \frac{1}{2} = \frac{1}{2} \cdot \frac{6^n}{6^n}, \:
    \frac{1}{2} + \frac{1}{2^2} \cdot \frac{1}{2^0 \cdot 6^n}, \:
    \cdots \quad
    \in D_n
  $$

  as 0th, 1st, 2nd, $ \cdots , \: 6^n $th, $6^{n+1}$th, $\cdots$
  element, respectively.

  Now, define functions
  $$
    j_n : D_n \rightarrow \{ 0, 1,2,3, \cdots \}
  $$
  $$
    t^n_1 : \big(D_n - \{ 0 \}\big) \cup \{ 1 \}  \rightarrow D_n
  $$
  $$
    t^n_2 :
    D_n - \{ \text{ the last element of } D_n  \} \rightarrow D_n
  $$

  as follows : \\

  $j_n(s) = j$ \quad for the $j$-th element $s \in D_n.$ \\

  $ t^n_1 (s) $ is the $(j-1)$-th element in $ D_n  $ for
  a given $j$-th element $ s \in D_n -\{ 0 \}$ and $t^n_1 (1)$ is the
  last element in $D_n$.\\

  $ t^n_2 (s) $ is the $(j+1)$-th element in $ D_n $ for a given $j$-th
  element $ s \in D_n -  \{ $ the last element of $ D_n  \} .$ \\

  \subsection{Definition of $ \mathbf{\gamma ^n_{t_0} , c^n_{t_0} ,
  \bar{c}^n_{t_0}, {_1}c^n_{t_0}, {_1}\bar{c}^n_{t_0},
  \varphi ^n_{t_0}} $ and $\mathbf{\psi ^n_{t_0}}$ on the disk $ \mathbf{D^2} $}

  Recall from Proposition \ref{prop} that the union $U_i$ of triangles from 1st one
  to $i$-th one is diffeomorphic to a disk.

  Let $ n \in \{1,2,3, \cdots \} $ and $ t_0 \in D_n $ be given.

  With respect to the ordering of $D_n$, we will define
  $ \gamma ^n_{t_0} , c^n_{t_0} , \bar{c}^n_{t_0} $ and
  $ \varphi ^n_{t_0} $ inductively for each fixed n:

  Case 1) $ t_0 $ is the first element in $ D_n $ , in fact,
          $ t_0 = \frac{1}{2} \cdot \frac{1}{{6^n}} $

    \bigskip

    The orientation at the barycenter of $ T_0 \in A_0 $ will give the
    direction of the boundary curve of the first triangle in $A_n.$

    Then
    $$ c^n_{t_0} : [0, 1] \rightarrow \{ basepoint \} \subset D^2 $$

    $$
      \bar{c} ^n_{t_0} : [ 0, 1] \rightarrow
      \{ basepoint \} \subset D^2
    $$

    $$ \varphi ^n_{t_0} : [ 0, 1] \rightarrow D^2 $$

    and

    $$ \gamma ^n_{t_0} : [0, 1 ] \rightarrow D^2 $$

    \noindent
    can be thought, where $ \varphi ^n_{t_0} $ and
    $ \gamma ^n_{t_0} $ are the piecewise smooth boundary curve of
    the first triangle in $A_n$ with constant speed
    and the direction of the boundary curve is induced from the given
    orientation.

    Note $ \gamma ^n_{t_0} $ can be regarded as the simplification of  $ c^n_{t_0} * \varphi ^n_{t_0} *
    \bar{c} ^n_{t_0} \: .$

    We will call $ \gamma ^n_{t_0} $ the \emph{holonomy curve at time
    $ t = t_0 $}.

    Now,  consider the path from the basepoint to the ending point of
    the first triangle in $n$-step along the opposite
    direction of the holonomy curve $ \gamma ^n_{t_0} $ at $t={t_0}$
    , which is a piecewise smooth curve with constant speed. Then from
    the path, we can define a piecewise smooth curve

    $$ _1 c ^n_{t_0} : [0,1] \rightarrow D^2 $$

    \noindent
    with constant speed.
    And its opposite direction can make us define

    $$ _1 \bar{c} ^n_{t_0} : [0,1] \rightarrow D^2 .$$

    Define a piecewise smooth curve

    $$ \psi ^n_{t_0} : [0,1] \rightarrow D^2 $$

    \noindent
    with constant speed as the boundary curve of the 1st triangle in
    the $n$-th step, where the curve is a loop at the ending point of
    the first triangle and the direction of the boundary curve is
    induced from the given orientation.

  \bigskip

  Case 2) $ t_0 $ is the $j$-th element in $ D_n $,
          \emph{i.e.}, $j_n (t_0) = j$,
          where $ j \geq 2 $

    \bigskip

    Let $ t_1 $ be the $(j-1)$-th element in $ D_n $,
    \emph{i.e.}, $t^n_1 (t_0) = t_1$ and $j_n (t_1) = j-1 $,
    where $ j-1 \geq 1 .$

    Consider the path from the basepoint  to the starting point of
    the $j$-th triangle in the $n$-th step along the opposite
    direction of the holonomy curve $ \gamma ^n_{t_1} $ at
    $ t = t_1 $ , which is a piecewise smooth one with constant speed
    .  Then from the path, we can define a piecewise smooth curve

    $$ c^n_{t_0} : [0, 1] \rightarrow  \partial U_{j-1} \subset D^2 $$

    \noindent  with constant speed, where $U_{j-1}$ is the union
    of triangle in $A_n$ from the 1st one to the $(j-1)$-th one.

    And its opposite direction can make us define

    $$ \bar{c} ^n_{t_0} : [0 , 1] \rightarrow  \partial U_{j-1} \subset D^2 .$$

    Define a piecewise smooth curve

    $$ \varphi ^n_{t_0} : [0,1] \rightarrow D^2 $$

    \noindent
    with constant speed as the boundary curve of the $j$-th
    triangle in the $n$-th step, where the curve is a loop at the
    starting point of the triangle and  the direction of the boundary
    curve is induced from the given orientation.

    Now define a piecewise smooth curve

    $$ \gamma ^n_{t_0} : [0, 1 ] \rightarrow  \partial U_{j} \subset D^2 $$

    \noindent
    with constant speed from the simplification of
    $
      \gamma ^n_{t_1} * c^n_{t_0} * \varphi ^n_{t_0} * \bar{c}
      ^n_{t_0},
    $
    where $U_j$ is the union
    of triangle in $A_n$ from the 1st one to the $j$-th one.
    The new curve will be also
    called the \emph{holonomy curve at time $ t=t_0 $ }.


    Now,  consider the path from the basepoint to the ending point of
    the $j$-th triangle in the $n$-th step along the opposite
    direction of the holonomy curve $ \gamma ^n_{t_0} $ at $t={t_0}$
    , which is a piecewise smooth one with constant speed. Then from
    the path, we can define a piecewise smooth curve

    $$ _1 c ^n_{t_0} : [0,1] \rightarrow \partial U_j \subset D^2 $$

    \noindent
    with constant speed.
    And its opposite direction can make us define

    $$ _1 \bar{c} ^n_{t_0} : [0,1] \rightarrow \partial U_j \subset D^2 .$$

    Define a piecewise smooth curve

    $$ \psi ^n_{t_0} : [0,1] \rightarrow D^2 $$

    \noindent
    with constant speed as the boundary curve of the $j$-th triangle in
    the $n$-th step, where the curve is a loop at the ending point of the
    $j$-th triangle and the direction of the boundary curve is induced
    from the given orientation.
    \bigskip

  \subsection{
               the simplification of
               $\mathbf{\bar{c}^n_{t_0} * {_1 c^n_{t_0}}}$
             }

  \bigskip


  For each $n \geq 1$ and $0 \neq t_0 \in D_n$, where  $t_0$ is the
  $j_n (t_0)$-th element in $D_n$, the simplification of
  $\bar{c}^n_{t_0} * {_1 c^n_{t_0}}$ is a curve along the boundary
  of $j_n (t_0)$-th triangle in $A_n$ with opposite
  direction to the given orientation such that it
  starts from the starting point of the triangle and that its image
  consists of the following sets :

  one point, one side, two sides or the boundary of the triangle.

  \bigskip

  \subsection{The induced curves on the surface
  $ S  \subset \mathbb{H}^n $ and totally geodesic planes in
  $ \mathbb{H}^n $}

    Let $ \Phi : D^2 \rightarrow S $ be a given diffeomorphism. Then
    we can think of triangles in $ S $ induced from the barycentric
    subdivision on $ D^2 $ on each $ n$-th  step. We will use
    ' $ \sim $ ' notation for the induced triangles and curves in
    $ S $ , that is ,

    $$
      \tilde{T} = \Phi (T)  \qquad \text{for}  \quad  T \in A_n
    $$

    $$ \tilde{A}_n = \{ \Phi (T) \mid T \in A_n \}  $$

    \noindent
    and

    $$
      \tilde{\gamma}^n_{t_0} ,
      \tilde{c}^n_{t_0} ,
      \tilde{\varphi}^n_{t_0} ,
      \tilde{\bar{c}}^n_{t_0} ,
      {_1 \tilde{c}}^n_{t_0} ,
      \tilde{\psi}^n_{t_0},
      {_1 \tilde{\bar{c}}}^n_{t_0},
    $$

    \noindent
    which are piecewise smooth curves with constant speed such that

    $$
      Im (\tilde{\gamma}^n_{t_0})  = Im (\Phi \circ \gamma^n_{t_0})
    $$

    $$
      Im ( \tilde{c}^n_{t_0})      = Im (\Phi \circ c^n_{t_0})
    $$

    $$
      Im (\tilde{\varphi}^n_{t_0}) = Im (\Phi \circ \varphi ^n_{t_0})
    $$

    $$
      Im(\tilde{\bar{c}}^n_{t_0})  = Im (\Phi \circ \bar{c}^n_{t_0})
    $$

    $$
      Im ( {_1 \tilde{c}}^n_{t_0}) = Im (\Phi \circ {_1 c}^n_{t_0})
    $$

    $$
      Im (\tilde{\psi}^n_{t_0})   = Im (\Phi \circ \psi ^n_{t_0})
    $$

    $$
      Im({_1 \tilde{\bar{c}}}^n_{t_0})
      = Im (\Phi \circ {_1 \bar{c}}^n_{t_0})
    $$

    \noindent
    and whose direction relates to that of
    $
      \gamma ^n_{t_0} , c^n_{t_0} , \varphi ^n_{t_0} ,
      \bar{c}^n_{t_0}, {_1 c}^n_{t_0} , \psi ^n_{t_0} ,
      {_1 \bar{c}}^n_{t_0} \: ,
    $
    respectively.

    Now with respect to each triangle in $ S $, we can think of a
    totally geodesic triangle with same vertices in $ \mathbb{H}^n $
    . So, each step will induce the similar concept , $ i.e. $
    triangles and curves,  on the induced pleated surface consisting
    of totally geodesic triangles and we'll use ' $ \wedge $ '
    notation for them. In other words, we can think of

    $$
      \hat{T} \in \hat{A}_n , \:
      \hat{\gamma}^{n}_{t_0} , \: \hat{c}^{n}_{t_0} , \:
      \hat{\varphi}^{n}_{t_0} , \: \hat{\bar{c}}^{n}_{t_0} , \:
      {_1 \hat{c}}^n_{t_0} , \: \hat{\psi}^n_{t_0} , \:
      {_1 \hat{\bar{c}}}^n_{t_0} \: ,
    $$

    \noindent
    where the curves
    $
      \hat{\gamma}^{n}_{t_0} , \hat{c}^{n}_{t_0} ,
      \hat{\varphi}^{n}_{t_0} , \hat{\bar{c}}^{n}_{t_0} ,
      {_1 \hat{c}}^n_{t_0} , \hat{\psi}^n_{t_0} ,
      {_1 \hat{\bar{c}}}^n_{t_0}
    $
    are piecewise geodesics in $ \mathbb{H}^n, $ induced from the
    boundaries of totally geodesic triangles $ \hat{T}$,
    and are relating  to the previous curves
    $
      \tilde{\gamma}^n_{t_0} , \tilde{c}^n_{t_0} ,
      \tilde{\varphi}^n_{t_0} , \tilde{\bar{c}}^n_{t_0},
      {_1 \tilde{c}}^n_{t_0} , \tilde{\psi}^n_{t_0} ,
      {_1 \tilde{\bar{c}}}^n_{t_0}
    $
    in $ S $ and
    $
      \gamma ^n_{t_0} , c ^n_{t_0} , \varphi ^n_{t_0} ,
      \bar{c} ^n_{t_0}, {_1 c}^n_{t_0} , \psi ^n_{t_0} ,
      {_1 \bar{c}}^n_{t_0}
    $
    in $ D^2 .$

  \pagebreak

\section{In
         $
           \text{SO}(2) \rightarrow \text{SO}_0 (1,2) \rightarrow
           \mathbb{H}^2
         $
         }\label{sec:n=2}

  For

    $$
        E_1  =
        \left(
          \begin{matrix}
            0 & 0 & 1   \\ 0 & 0 & 0    \\ 1 & 0 & 0
          \end{matrix}
        \right)
        ,
        E_2  =
        \left(
          \begin{matrix}
            0 & 1 & 0    \\ 1 & 0 & 0   \\ 0 & 0 & 0
          \end{matrix}
        \right)
        ,
        \text{and }
        E_3  =  [E_1, E_2]  =
        \left(
          \begin{matrix}
            0 & 0 & 0    \\ 0 & 0 & -1  \\ 0 & 1 & 0
          \end{matrix}
        \right)
           ,
   $$
   \noindent
   let $\{ E_1, E_2, E_3 \}$ be an ordered orthonormal basis of
   $\mathfrak{so}(1,n),$ which induces the canonical orientation on
   $\text{SO}_0 (1,2)$, and so induces
   the canonical orientation on $\mathbb{H}^2, $ that is , the
   counterclockwise one.

  For $t \in \mathbb{R}$, put

  $$
    \Psi(\alpha) \; = \; \text{exp} (t E_3) \; = \;
    \left(
      \begin{matrix}
         1 &   0     &    0       \\
         0 & \cos{t} & - \sin{t}  \\
         0 & \sin{t} & \cos{t}
      \end{matrix}
    \right) .
  $$

  Let $ c : [t_0, t_3] \rightarrow \mathbb{H}^2 $ be a simple-closed
  arc-length parameterized piecewise-smooth curve representing a
  geodesic triangle in $  \mathbb{H} ^2 : $

  \begin{figure}[h]
    \centering
      {
        \includegraphics[width=1.5in]{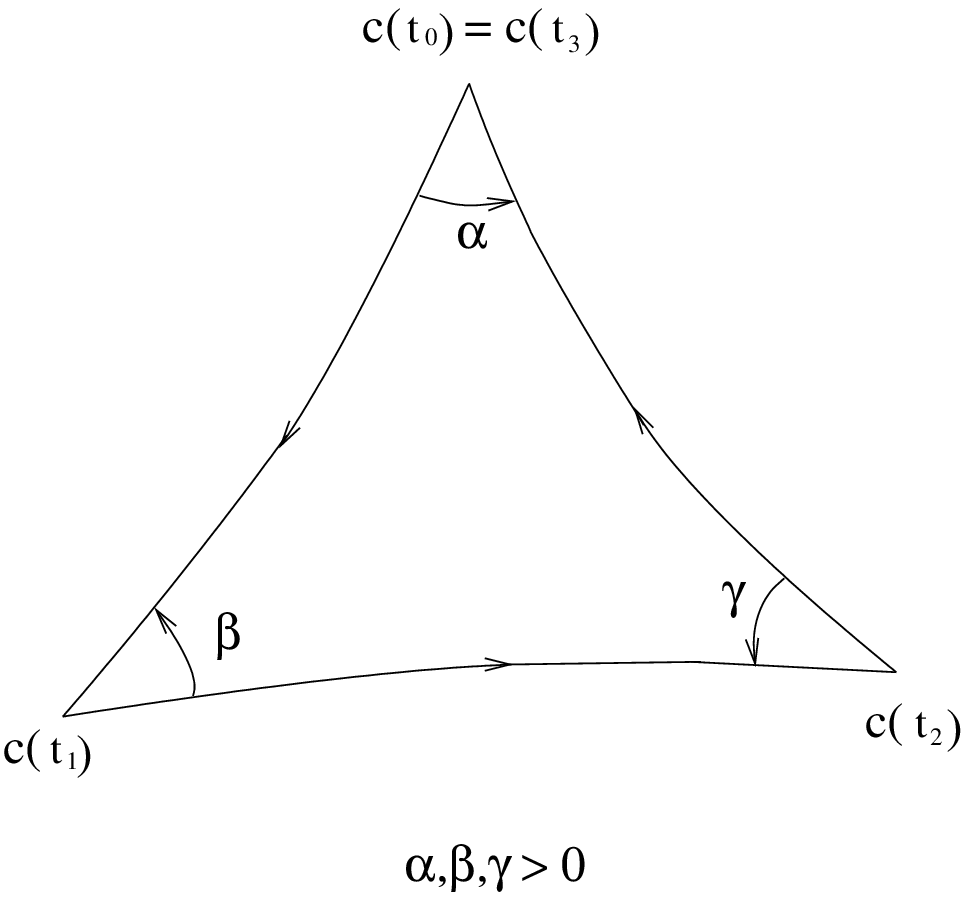}
        \hspace{1cm} or \hspace{1cm}
        \includegraphics[width=1.5in]{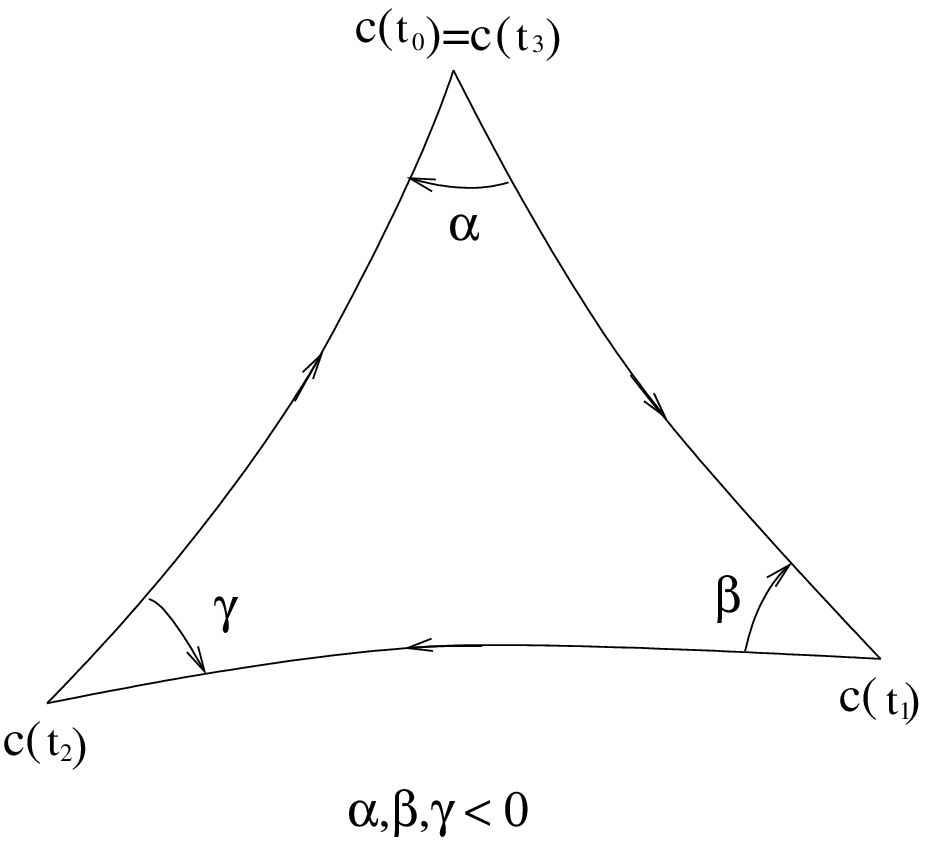}
      }
  \end{figure}

  More precisely, \\
  c is continuous on $[t_0, t_3]$ and smooth on
  $(t_0, t_1) \cup (t_1, t_2) \cup (t_2, t_3)$ , where
  $ c(t_0) = c(t_3), \text{ } c(t_1) \text{ and } c(t_2)$ are
  vertices of the given geodesic triangle.

  Let
  \begin{align*}
   &\alpha \text{ be the angle from } \dot{c}({t_0}^+)
   \text{ to } - \dot{c} ({t_3}^-) ,
   \\
   &\beta  \text{ be the angle from } \dot{c}({t_1}^+)
   \text{ to } - \dot{c} ({t_1}^-) , \text{ and }
   \\
   &\gamma \text{ be the angle from } \dot{c}({t_2}^+)
   \text{ to } - \dot{c} ({t_2}^-) \: .
  \end{align*}

  Then, either
  $
    \alpha , \beta , \gamma > 0
    \text{ or }
    \alpha , \beta , \gamma < 0
  $
  holds.

  \begin{lemma}
  Under the above condition, let
  $ \tilde{c} : [t_0 , t_3] \rightarrow \text{ SO} _0 (1, 2) $
  be a horizontal lift of $c$.
  Then, the relation between the holonomy and the area of a geodesic
  triangle is given by
  $$
    {{\tilde{c} (t_0)}^{-1}}  \cdot  \tilde{c} (t_3)  =
    \Big( \Psi \big(\pi- |\, \alpha + \beta + \gamma \,| \big) \Big)
    ^{\delta}
    \hspace{5mm} ,
  $$
  where
  \begin{equation*}
    \delta =
    \left\{
      \begin{array}{rl}
         1 & \text{ if } \alpha, \beta, \gamma > 0 \\
        -1 & \text{ if } \alpha, \beta, \gamma < 0
      \end{array}
    \right.
  \end{equation*}

  Furthermore, $\pi- |\, \alpha+\beta+\gamma\, |$ is the area of
  the geodesic triangle.
  \end{lemma}

  \bigskip

  \begin{figure}[h]
    \centering
      {
        \includegraphics[width=2.1cm]{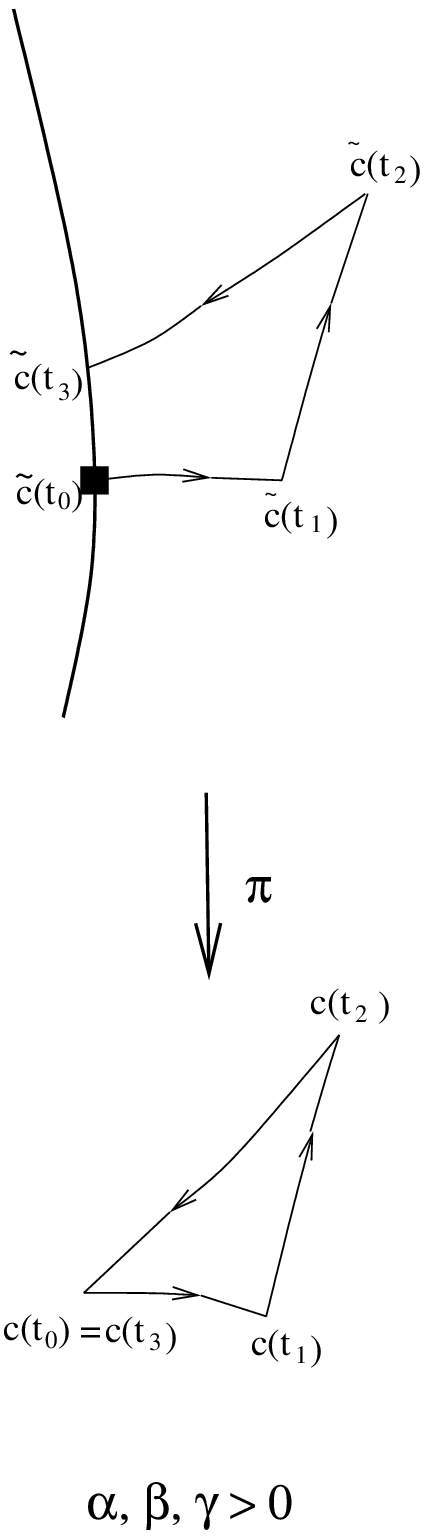}
        \hspace{1cm} or \hspace{1cm}
        \includegraphics[width=1.9cm]{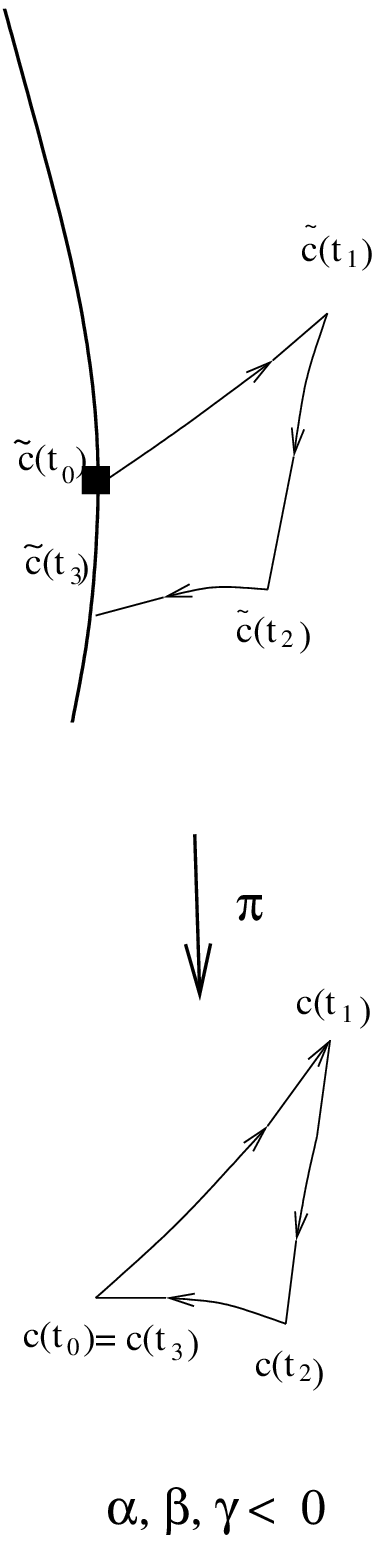}}
  \end{figure}

 \begin{proof}
  Let $ \pi : \text{ SO} _0 (1,2) \rightarrow \mathbb{H} ^2 $ be the
  given Riemannian submersion.

  Recall that for any $k \in \text{SO(2)},$ the restriction
  $Ad_k |_{\mathfrak{so}(2)^{\perp}}$ of
  $Ad_k (\cdot) : \mathfrak{so}(1,2) \rightarrow \mathfrak{so}(1,2)$
  to $\mathfrak{so}(2)^{\perp}$ is an automorphism of
  $\mathfrak{so}(2)^{\perp}.$
  For $-\dot{c} ({t_3}^{-})$ and its horizontal lift $x$ at
  $\tilde{c} (t_0)$, find $A \in \mathfrak{so}(2) ^{ \perp}$ satisfying
  $$
    {L_{\tilde{c}(t_0)^{-1}}}_{*} x = A_e.
  $$

  Then
  $
    \dot{\tilde{c}} ({t_0}^+) ,
    \dot{\tilde{c}} ({t_1}^+)
    \text{ and }
     \dot{\tilde{c}} ({t_2}^+)
  $
  will be related to
  $$
    Ad_{\Psi(\alpha)} A , \hspace{5mm}
    Ad_{ \Psi(\pi) \cdot \Psi(\alpha) \cdot \Psi(\beta) } A
    \hspace{5mm} \text{and} \hspace{5mm}
    Ad_{ \Psi(\alpha) \cdot \Psi(\beta) \cdot \Psi(\gamma) } A
    \hspace{5mm},
  $$
  respectively, by

  \begin{align*}
    {{L_{\tilde{c}(t_0)^{-1}}}_{*}} \, \dot{\tilde{c}} ({t_0} ^{+})
    &= (Ad_{\Psi(\alpha)} A )_e
    \\
    {{L_{\tilde{c}(t_1)^{-1}}}_{*}} \, \dot{\tilde{c}} ({t_1} ^{+})
    &=
    \big(Ad_{\Psi(\beta)} (- Ad_{\Psi(\alpha)} A )\big)_e
    = (Ad_{\Psi(\pi) \cdot \Psi(\alpha) \cdot {\Psi(\beta)} } A)_e
    \\
    {{L_{c(t_2)^{-1}}}_{*}} \, \dot{\tilde{c}} ({t_2} ^{+})
    &=
    \big(
      Ad_{\Psi(\gamma)}
      (- Ad_{\Psi(\pi) \cdot \Psi(\alpha) \cdot \Psi(\beta)}  A )
    \big)_e
    = (Ad_{\Psi(\alpha) \cdot \Psi(\beta) \cdot \Psi(\gamma)} A)_e,
  \end{align*}

  \noindent
  in other words,

  \begin{align*}
    \dot{\tilde{c}} ({t_0} ^{+})
    &= {{L_{\tilde{c}(t_0)}}_{*}}_e (Ad_{\Psi(\alpha)} A )_e
    \\
    \dot{\tilde{c}} ({t_1} ^{+})
    &= {{L_{\tilde{c}(t_1)}}_{*}}_e
      (Ad_{\Psi(\pi) \cdot \Psi(\alpha) \cdot {\Psi(\beta)} } A)_e
    \\
    \dot{\tilde{c}} ({t_2} ^{+})
    &= {{L_{\tilde{c}(t_2)}}_{*}}_e
    (Ad_{\Psi(\alpha) \cdot \Psi(\beta) \cdot \Psi(\gamma)} A)_e
  \end{align*}

  And so

  $ - \dot{\tilde{c}} ({t_3}^-)$ will be related to
  $$
    - Ad_{\Psi(\alpha) \cdot \Psi(\beta) \cdot \Psi(\gamma)} A
    =
    Ad
    _{
      \Psi(\pi) \cdot \Psi(\alpha) \cdot \Psi(\beta) \cdot
      \Psi(\gamma)
     }
    A
    \hspace{5mm},
  $$
  that is ,
  $$
    - \dot{\tilde{c}} ({t_3}^-) =
    {{L_{\tilde{c} (t_3)}}_{*}}_{e}
    (
      Ad
      _{
        \Psi(\pi) \cdot \Psi(\alpha) \cdot \Psi(\beta) \cdot
        \Psi(\gamma)
       }
      A
    )_e
    \hspace{5mm} .
  $$

  Therefore,

  \begin{align*}
    \pi \big(\tilde{c} (t_0) \cdot \text{e} ^{tA}\big)
    &= \pi
    \big(
      \tilde{c} (t_3) \cdot
      \text{e}
      ^{
         t
         (
           Ad
           _{
             \Psi(\pi) \cdot \Psi(\alpha) \cdot \Psi(\beta) \cdot
             \Psi(\gamma)
            }
           A
         )
      }
    \big)
    \\
    &= \pi
    \big(
      \tilde{c} (t_0) \cdot
      \text{e}
      ^{
         t
         Ad
         _{
           {\tilde{c}(t_0)}^{-1} \cdot \tilde{c} (t_3) \cdot
           \Psi(\pi) \cdot \Psi(\alpha) \cdot \Psi(\beta) \cdot
           \Psi(\gamma)
          }
         A
      }
      \cdot {\tilde{c} (t_0)}^{-1} \cdot \tilde{c} (t_3)
    \big)
    \\
    &= \pi
    \big(
      \tilde{c} (t_0) \cdot
      \text{e}
      ^{
         t
         Ad
         _{
           {\tilde{c}(t_0)}^{-1} \cdot \tilde{c} (t_3) \cdot
           \Psi(\pi) \cdot \Psi(\alpha) \cdot \Psi(\beta) \cdot
           \Psi(\gamma)
         }
         A
      }
    \big)
  \end{align*}
  \noindent
  because
  $$ {\tilde{c}(t_0)}^{-1} \cdot \tilde{c}(t_3) \in \text{ SO(2) } $$
  and
  $$
   B  \in  {\text{ so}(2)}^{\perp} \text{ and } k \in \text{ SO} (2)
   \Rightarrow k \cdot e^{tB} \cdot k^{-1} = e^{t Ad_{k} B}
   \hspace{5mm} .
  $$

  Thus, we get
  $$
    A =
    Ad
    _{
      {\tilde{c}(t_0)}^{-1} \cdot \tilde{c} (t_3) \cdot \Psi(\pi)
      \cdot \Psi(\alpha) \cdot \Psi(\beta) \cdot \Psi(\gamma)
     }
    A
  $$
  and so
  $$
    {\tilde{c} (t_0)}^{-1} \cdot \tilde{c} (t_3) \cdot
    \Psi \big( \pi + ( \alpha + \beta + \gamma ) \big) \;
     = \;
     \Psi(2n \pi) \hspace{5mm} \text{ for some } n \in \mathbb{Z}
     \: .
  $$

  Therefore,
  \begin{align*}
   {\tilde{c}(t_0)}^{-1} \cdot \tilde{c}(t_3)
   &= \Psi(2n \pi) \cdot
   \Big( \Psi \big(\pi + (\alpha + \beta + \gamma) \big) \Big)^{-1} \\
   &= \Big(
        \Psi \big(\pi + (\alpha + \beta + \gamma) \big)
      \Big)^{-1}
   \\
   &= \left\{
        \begin{array}{rl}
          \Psi(- \pi - (\alpha + \beta + \gamma))
          \hspace{5mm }
          & \text{ if } \alpha, \beta, \gamma > 0
          \\
          \Big(
            \Psi \big(\pi + (\alpha + \beta + \gamma) \big)
          \Big)^{-1}
          & \text{ if } \alpha, \beta, \gamma < 0
        \end{array}
      \right.
   \\
   &= \left\{
        \begin{array}{rl}
          \Psi(\pi - (\alpha + \beta + \gamma))
          \hspace{29mm }
          & \text{ if } \alpha, \beta, \gamma > 0
          \\
          \bigg(
            \Psi
            \Big(\pi - \big((-\alpha) + (-\beta) + (-\gamma)\big)
            \Big)
          \bigg)^{-1}
          & \text{ if } \alpha, \beta, \gamma < 0
        \end{array}
      \right.
   \\
   &= \Big( \Psi
      \big( \pi - |\alpha + \beta + \gamma| \big)\Big)^{\delta} ,
      \text{ where }
      \delta =
      \left\{
        \begin{array}{rl}
           1 & \text{ if } \alpha, \beta, \gamma > 0 \\
          -1 & \text{ if } \alpha, \beta, \gamma < 0
        \end{array}
                                                                                                                                     \right.
  \end{align*}

 \end{proof}

\section{Liftings in
          $
            \text{SO}(n) \rightarrow \text{SO}_0 (1,n) \rightarrow
            \mathbb{H}^n
          $
        }\label{sec:n-general}

  This section is the proof of Theorem ~\ref{thm}.

  \subsection{Preliminaries on the Riemannian submersion
  $ \pi : \text{SO}_0 (1,n) \rightarrow \mathbb{H}^n $}

  \bigskip

  Let $ \pi : \text{SO}_0 (1,n) \rightarrow \mathbb{H}^n $ be the
  given Riemannian submersion.
  In fact, this is the quotient of the isometric right translation
  by $\text{SO}(n)$ and $ \mathbb{H}^n $ is isometric to
  $ \text{SO}_0 (1,n) / \text{SO}(n) .$

  Let $ G = \text{SO}_0 (1,n) $ , $ K = \text{SO}(n) $  , and
  $ \mathfrak{g} $ and $ \mathfrak{k} $ be their Lie algebras,
  respectively.

  \bigskip

  Fact 1) $ X \in \mathfrak{k}^{\bot} $ $\Rightarrow $
  $t \mapsto g\cdot$  exp $(tX) : (- \infty , \infty) \rightarrow G $
  is a horizontal geodesic for any $ g \in G $ .

  \bigskip

  Fact 2) $ X,Y \in \mathfrak{k}^{\bot} $  $ \Rightarrow $
  $ [X,Y] \in \mathfrak{k} $ and
  Span$ \{ X,Y, [X,Y] \} \subset \mathfrak{g} $ is a Lie subalgebra
  $\mathfrak{h}$ of $\mathfrak{g}$ and its related subgroup $H$
  of $\text{SO}_0(1,n)$ is isometric
  to $\text{SO}_0(1,2)$. Furthermore, the riemannian submersion
  $\text{ SO}_0(1,n) \rightarrow \, \text{SO}(n)$  can be restricted to
  $\big( H= \text{ SO}_0(1,2) \big) \rightarrow \, \text{SO}(2).$

  \bigskip

  Fact 3) $ U \in \mathfrak{k} $   $ \Rightarrow $
  $ t \mapsto k \cdot $
  exp $ (tU) : ( - \infty , \infty) \rightarrow K \subset G $ is a
  vertical geodesic for any $ k \in K .$

  \bigskip

  Fact 4) For any $ k \in K $ , the right translation
  $R_k : G \rightarrow G $ by $ k , R_k (g) = gk, $ is an isometry .
  Or, equivalently, $ Ad_k : \mathfrak{g} \rightarrow \mathfrak{g}  $
  is a linear isometry for any $ k \in K .$

 \bigskip

 \subsection{
              Definition of
              $
               \mathbf{
                        \bar{f}:
                        \displaystyle{\bigcup ^{\infty}_{m=1}} D_m
                        \rightarrow K = \text{SO}(n),
                        \bar{f}_m: D_m \rightarrow  \text{SO}(n)
                      }
              $
              and
              $\mathbf{\hat{f}_m : [0,1] \rightarrow \text{SO}(n)}$
              and the property of $\mathbf{\hat{f}_m}$
            }

 \subsubsection{\textbf{Definition of} $\, \mathbf{\bar{f}, \bar{f}_m}$}

  Let $\bar{f}(0)=e.$
  Fix
  $ t_0 \in \displaystyle{\bigcup ^{\infty}_{m=1}} D_m - \{ 0 \} .$
  Then we can find a positive integer
  $ m_0 $ =
  min $\{ m_1   \mid m+1 \geq m_1 \Rightarrow  t_0 \in D_m \} \: .$

  Note that on the given surface $S$,
  $$
    \tilde{\gamma} ^n_{t_0} = \tilde{\gamma} ^{n_0}_{t_0}
    \text{ and } \;
    \tilde{c}^n_{t^n_2 (t_0)}  = {_1 \tilde{c}}^n_{t_0} =
    {_1 \tilde{c}}^{n_0}_{t_0} = \tilde{c}^{n_0}_{t^{n_0}_2 (t_0)}
  $$
  for all $ n \geq n_0 $ . So let
  $$
    \tilde{\gamma} _{t_0} := \tilde{\gamma} ^{n_0}_{t_0}
    \text { and } \;
    {_1 \tilde{c}}_{t_0} :={_1 \tilde{c}}^{n_0}_{t_0}
  $$

  Define
  \begin{align*}
    \bar{f}(t_0) :
    & = \text{ the value , at } t=1 ,
        \text{ of the horizontal lifting of } \tilde{\gamma} _{t_0}
        \text{ at } e \: .
  \end{align*}

  \medskip
  Put $\bar{f}_n$ as the restriction
  $\hat{f}_n \mid_{{\bigcup ^{\infty}_{n=1}} D_n }$
  of $\hat{f}_n$, defined below, to $ D_n $.

 \pagebreak
 \subsubsection{
                 \textbf{
                         Definition of $\,\mathbf{\hat{f}_n}$
                         and its property
                        }
               }

  Let's define a curve $ \hat{f}_n : [0,1] \rightarrow K = \text{SO}(n) $
  with
  $ \hat{f}_n (0) = e $ inductively as follows:

  \bigskip

  Step 1) Assume $ t_0 \in D_n $ is the 1st element in $ D_n $, in
  fact, $ t_0 = \frac{1}{2} \cdot \frac{1}{6^n} $. Then,
  $t^n_1 (t_0) = 0 .$

  \begin{figure}[h]
       \centering{\includegraphics[width=3.5cm]{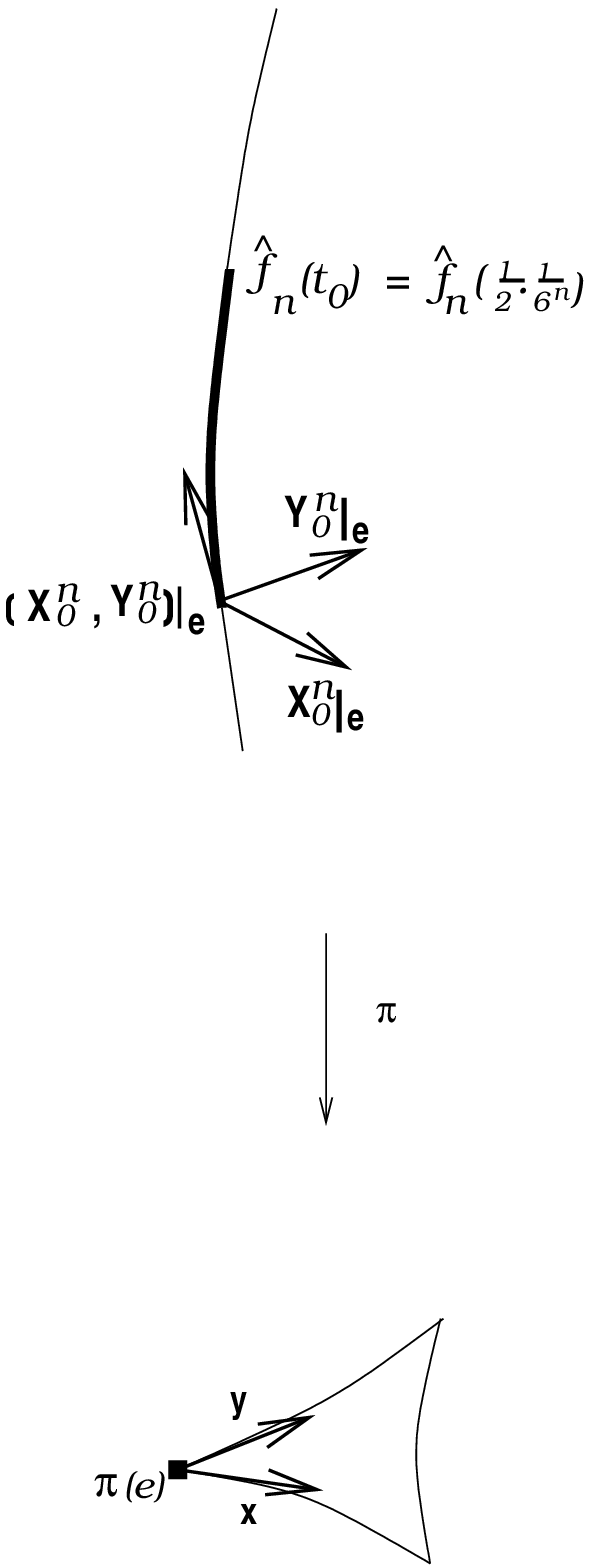}}
  \end{figure}

  Consider the first triangle in $ \hat{A}_n $, its starting point
  and  the horizontal lifting of
   $$
     x :=
     \lim _{t \rightarrow 0^{+}}
     \frac{1}{\mid \dot{ \hat{\varphi} } ^n_{t_0} (t) \mid } \cdot
     \dot{ \hat{\varphi} } ^n_{t_0} (t)
   $$
   and
   $$
     y :=
     - \lim _{t \rightarrow 1^{-}}
     \frac{1}{\mid \dot{ \hat{\varphi} } ^n_{t_0} (t) \mid } \cdot
     \dot{ \hat{\varphi} } ^n_{t_0} (t)
   $$
  at $ e $ , respectively and find
  $$
    X^n_0 = X^n_{t^n_1(t_0)}  , \:
    Y^n_0 = Y^n_{t^n_1(t_0)}
    \in \mathfrak{k}^{ \bot}
  $$
  with
  $$
    \pi _{*} \; X^n_0 \mid _{e} \:
    = \: \pi _{*} \; X^n_{t^n_1(t_0)} \mid _{e} \:
    = \: x
  $$
  and
  $$
    \pi _{*} \; Y^n_0 \mid _{e} \:
    = \: \pi _{*} \; Y^n_{t^n_1(t_0)} \mid _{e} \:
    = \: y  \: .
  $$
  Then, define
  $$
     \hat{f}_n (t) :=
     \text{ exp }
     \Big(
       t \cdot
       \frac { ( \text{Area of the 1st triangle in } \hat{A}_n )}
             { t_0 \cdot \mid  [X^n_0 , Y^n_0 ] \mid }
       \cdot [X^n_0 , Y^n_0 ]
     \Big)
     \qquad
     \text{for }  t \in [0, t_0 ] \; ,
  $$
  which is a geodesic in $ K= \text{SO}(n) $ from Fact 3.

  \bigskip

  Step 2) Assume  $ t_0 \in D_n $ is the $j$-th element in $ D_n $,
  where $ j \geq 2 \; .$

  Note $ t^n_1 (t_0) $ is the $(j-1)$-th element in $ D_n $, where $
  j-1 \geq 1 \; .$

  \begin{figure}[h]
       \centering{\includegraphics[width=2.5in]{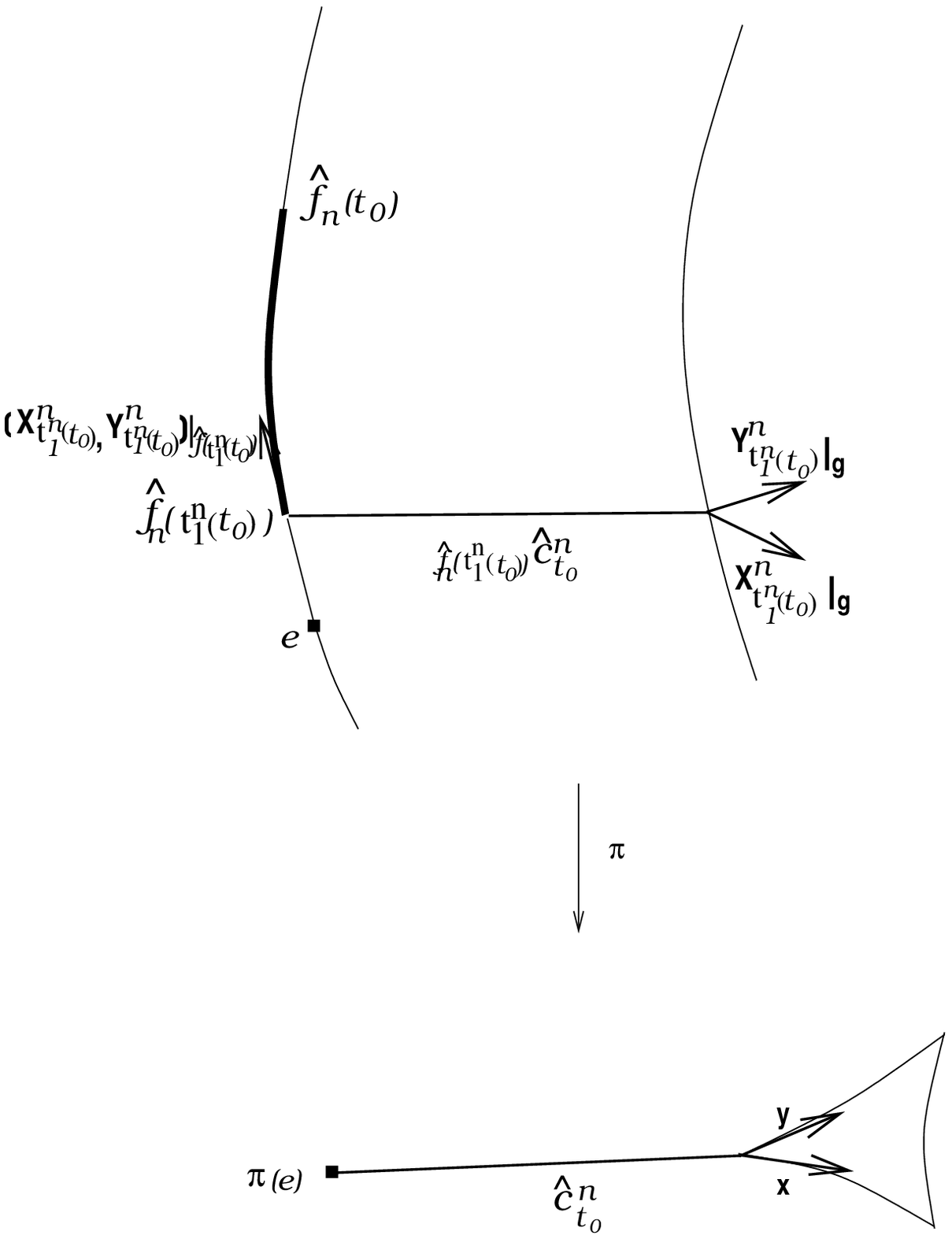}}
  \end{figure}

  Let
  $
    _{\hat{f}_n (t^n_1 (t_0) ) } \hat{c} ^n_{t_0} :
    [0,1] \rightarrow \text{SO}_0 (1,n)
  $
  be the horizontal lifting of $ \hat{c} ^n_{t_0} $ at
  $ \hat{f}_n (t^n_1 (t_0)) $ and then consider the $j$-th triangle in
  $\hat{A}_n$, its starting point and the horizontal lifting of
  $$
    x :=
    \lim _{t \rightarrow 0^{+}}
    \frac{1}{\mid \dot{ \hat{\varphi} } ^n_{t_0} (t) \mid } \cdot
    \dot{ \hat{\varphi} } ^n_{t_0} (t)
  $$
  and
  $$
    y :=
    - \lim _{t \rightarrow 1^{-}}
    \frac{1}{\mid \dot{ \hat{\varphi} } ^n_{t_0} (t) \mid } \cdot
    \dot{ \hat{\varphi} } ^n_{t_0} (t)
  $$
  at $ g := {_{\hat{f}_n (t^n_1 (t_0) ) } \hat{c} ^n_{t_0} (1)} $ ,
  respectively and find

  $$
    X^n_{t^n_1(t_0)} , \; Y^n_{t^n_1(t_0)}  \in \mathfrak{k}^{ \bot}
  $$
  with
  $$ \pi _{*} \: {X^n_{t^n_1(t_0)} \mid _{ g }} \: = \: x  $$
  and
  $$ \pi _{*} \: {Y^n_{t^n_1(t_0)} \mid _{ g }} \: = \: y \: .$$
  Then define

  \vspace{0.5cm}
  \noindent
  $
    \hat{f}_n (t) \: := \\
    \hat{f}_n(t^n_1 (t_0)) \cdot
    \text { exp }
    \Big(
      \big( t- t^n_1 (t_0) \big) \cdot
      \frac {(\text{Area of }  j \text{-th triangle in } \hat{A} _n)}
            {
             \big( t_0 - t^n_1 (t_0) \big) \cdot \:
               \mid \, [X^n_{t^n_1 (t_0)} , Y^n_{t^n_0 (t_0)} ] \, \mid
            }
      \cdot [X^n_{t^n_1 (t_0)} , Y^n_{t^n_1 (t_0)} ]
    \Big)
  $
  \\
  for $ t \in [t^n_1 (t_0), t_0 ] .$

  \bigskip

  Step3 ) $ t_0 = 1 $

  Note $ t^n_1 (1) $ be the last element in $ D_n $ , in other words,
  $t^n_1 (1) = \displaystyle {\sum ^{n+1}_{i=1}} \:  \frac{1}{2^i} .$

  Then define
  $$ \hat{f}_n (t) := \hat{f}_n (t^n_1 (1)) $$ for $ t \in [t^n_1 (1), 1] .$
  \\

  Now check the property of $ \hat{f}_n .$

  Assume $ 0 \neq t_0 \in D_n $ is a $j$-th element in $ D_n $, where
  $ j \geq 1 $ . Then $ t^n_1 (t_0) $ is the $ (j-1) $-th elements in
  $ D_n $ , where $ j-1 \geq 0 $ , and from Facts, mentioned early in
  this section, and from the property in Section ~\ref{sec:n=2}, we get

  \begin{align*}
    \hat{f}_n (t_0)
    & = \text { the value , at } t=1 ,
        \text{ of the horizontal lifting of }
        \hat{c} ^n_{t_0} * \hat{\varphi} ^n_{t_0} *
        \bar{\hat{c}} ^n_{t_0}
        \text { at }  \hat{f}_n (t^n_1 (t_0)) \\
    & = \text{ the value , at } t=1,
        \text{ of the horizontal lifting of }
        \hat{\gamma} ^n_{t_0} \text{ at  e } .
  \end{align*}

   \begin{figure}[h]
       \centering{\includegraphics[width=7cm]{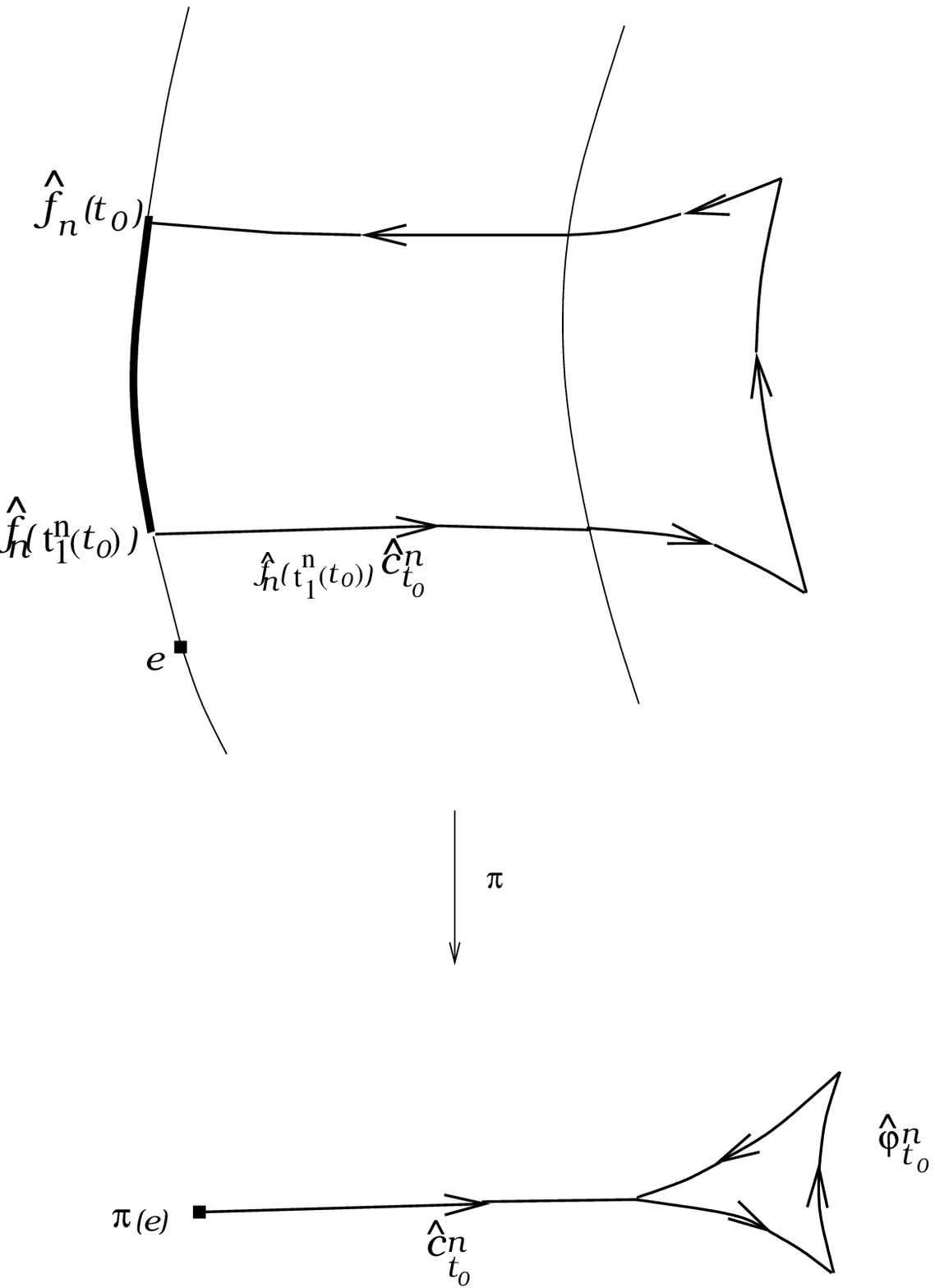}}
   \end{figure}

  Define, for any $ g \in G $, $ l_g : K \rightarrow G $  by
  $ l_g (k) = gk $ , which  is an isometric imbedding of $K$
  onto the fiber $ gK .$

  And let $ \omega $ and $ \Omega $ be the connection form and the
  curvature form of the connection of the principal bundle
  $ \pi : \text{SO}_0 (1,n) \rightarrow \mathbb{H} ^n ,$
  respectively.

  Then, under the identification of $T_e G$ and $\mathfrak{g},$  for
  $ t \in (t^n_1 (t_0), t_0) $ and
  $ g = {_{\hat{f}_n (t^n_1 (t_0))} \hat{c} ^n_{t_0} (1)} ,$  which is the
  value, at $ t =1$, of the horizontal lifting of
  $ \hat{c} ^n_{t_0} $ at $ \hat{f}_n (t^n_1 (t_0)) \; ,$

  \begin{align*}
    \omega \;
    (
      \; \tfrac{1}{\mid \dot{\hat{f}} _n (t) \mid } \cdot \dot{\hat{f}} _n (t) \;
    )
    \:
    & = ({{l_{\hat{f}_n (t)}}_{*}}_e)^{-1} \;
        (
          \; \tfrac{1}{\mid \dot{\hat{f}} _n (t) \mid } \cdot
          (\dot{\hat{f}} _n (t))^v \;
        )
    \\
    & = ({{l_{\hat{f}_n (t)}}_{*}}_e)^{-1} \;
        (
          \; \tfrac{1}{\mid \dot{\hat{f}} _n (t) \mid } \cdot
          \dot{\hat{f}} _n (t) \;
        )
    \\
    & = ({{l_{\hat{f}_n (t)}}_{*}}_e)^{-1} \;
        (
          \; \tfrac{1}
                   {
                     \mid
                       [ X^n_{t^n_1 (t_0)}, \; Y^n_{t^n_1 (t_0)} ]
                       \mid _{\hat{f}_n (t)}
                     \mid
                   }
        \cdot
        [ X^n_{t^n_1 (t_0)}, \; Y^n_{t^n_1 (t_0)} ] \mid _{\hat{f}_n (t)})
    \\
    & = \tfrac{1}
              {
                \mid
                  ({{l_{\hat{f}_n (t)}}_{*}}_e)^{-1}
                  (
                    [ X^n_{t^n_1 (t_0)}, \; Y^n_{t^n_1 (t_0)} ]
                    \mid _{\hat{f}_n (t)}
                  )
                \mid
              }
        \cdot
        ({{l_{\hat{f}_n (t)}}_{*}}_e)^{-1} \;
        (
          \;
          [ X^n_{t^n_1 (t_0)}, \; Y^n_{t^n_1 (t_0)} ]
          \mid _{\hat{f}_n (t)}
        )
    \\
    & = \tfrac{1}
              {\mid [X^n_{t^n_1 (t_0)}, \; Y^n_{t^n_1 (t_0)}] \mid}
        \cdot  \; [ X^n_{t^n_1 (t_0)}, \; Y^n_{t^n_1 (t_0)} ]
    \\
    & = \; \tfrac {-1}
                  {
                    \mid
                      \Omega
                        (
                           X^n_{t^n_1 (t_0)} \mid _g ,\;
                           Y^n_{t^n_1 (t_0)} \mid _g
                        )
                    \mid
                  }
           \cdot
           \Omega \;
             (
               X^n_{t^n_1 (t_0)} \mid _g ,\;
               Y^n_{t^n_1 (t_0)} \mid _g
             )
  \end{align*}

  \noindent
  and
  $$
    \omega \;
    (
      \;
      \tfrac{1}{\mid \dot{\hat{f}} _n (t) \mid } \cdot \dot{\hat{f}} _n (t)
      \;
    )
    \mid _e \;
    = \;
    {L_{(\hat{f}_n (t))^{-1}}}_{*}
    (\tfrac{1}{\mid \dot{\hat{f}} _n (t) \mid } \cdot \dot{\hat{f}} _n (t))
  $$

  Roughly speaking, the unit tangent vector
  $
    \frac{1}{\mid \dot{f}_n (t) \mid} \cdot \dot{f}_n  (t),
  $
  $ t \in (t^n_1 (t_0), t_0), $ is the negative of the unit curvature
  of the 2-dimensional horizontal plane
      \quad
  $$
      \hat{H}^n_g
      = \text{Span}
          \{
            X^n_{t^n_1 (t_0)} \mid _g ,\; Y^n_{t^n_1 (t_0)} \mid _g
          \},
      \qquad          \text{where }
      g = _{\hat{f}_n (t^n_1 (t_0))} \hat{c} ^n_{t_0} (1),
  $$
  which projects to the tangent plane of the
         $j_n(t_0)$-th triangle in $ \hat{A}_n$
  at $ \pi (g) = \hat{c}^n_{t_0}(1) =$
        the starting point of the $j_n(t_0)$-th triangle in
      $ \hat{A}_n$.
  And, the length of  $\hat{f}_n \mid _{[t^n_1 (t_0), t_0]}$
     is the area of the
         $j_n(t_0)$-th triangle in $ \hat{A}_n$.

  \bigskip

  \subsection{
              The convergence of $ \mathbf{\bar{f}_n}(t_0) $
              to $\mathbf{\bar{f}}(t_0)$
             } \label{sec-converge-lift}

  Recall
  $$
    \bar{f}(t_0)
     = \text{ the value , at } t=1 ,
        \text{ of the horizontal lifting of } \tilde{\gamma} _{t_0}
        \text{ at } e
  $$
  and
  $$
    \bar{f}_n (t_0) = \hat{f}_n (t_0)
     = \text{ the value , at } t=1,
        \text{ of the horizontal lifting of }
        \hat{\gamma} ^n_{t_0} \text{ at  e } .
  $$

  Consider our Riemannian submersion
  $$
  \CD
  \so(n) @>>> &\text{SO}_0(1,n) @>>> &\text{ SO}_0(1,n)/\so(n).
  \endCD
  $$
  This bundle has a global cross section $s: \bbh\ra \na\subset G$,
  which comes from the Iwasawa decomposition $\na K$,
  where $K = \text{ SO}(n)$. That is, every
  element of $G$ is uniquely written as $n a k$, and the projection
  maps this to $n a K\in \bbh$.

  The cross section $s$ provides us with a one-to-one correspondence
  between the space of all piecewise $C^k$-curves in $\bbh^n$ and in
  $\text{ SO}_0 (1,n)$, with initial points $\bar e$ and $e$, by
  $$
    h \longleftrightarrow s\circ h.
  $$
  By abusing of notation, express $s \circ h$ by $h$.
  For a curve $h: [0,1]\ra \bbh^n$, the unique horizontal lift
  $\tilde h: [0,1]\ra \sono$ is given by
  $$
    h(t)\cdot a(t)=\tilde h(t)
  $$
  for a unique curve $a(t)$ in $\so(n)$. Such an $a(t)$ is obtained
  by solving the differential equation
  \begin{equation}
  \label{DE1}
  \angles{h\inv h' + a' a\inv,\ V}=0
  \end{equation}
  for every $V\in\frakk$, where $'$ means the derivative with
  respect to $t$. Note that the first entry $h\inv h' + a'\cdot
  a\inv$ is an element of the Lie algebra $\fsono$. The equation
  (\ref{DE1}) comes about as follows. The curve $\tilde h(t)$ being
  horizontal implies the following equalities should hold.
  \begin{align*}
    0
    &=\angles{(h(t) a(t))',(h(t) a(t)) V}\\
    &=\angles{h'(t) a(t) + h(t) a'(t) ,(h(t) a(t)) V}\\
    &=\angles{(h(t) a(t))\left(a(t)\inv h(t)\inv  h'(t) a(t) +  a\inv(t)
   a'(t)\right),(h(t) a(t)) V},
  \end{align*}
  for every $V\in \frakk$, on the tangent space at $h(t) a(t)$.
  Since the metric on $G$ is left-invariant, this implies
  \begin{align*}
   0&=\angles{a(t)\inv h(t)\inv  h'(t) a(t) +  a\inv(t) a'(t), V},
  \end{align*}
  for every $V\in \frakk$, on the tangent space at $e$,
  $G_e=\frakg$. Since this holds for all $V\in\frakk$ and $a(t)\in
  K$, by taking conjugation by $a(t)$, the above is equivalent to
  the equality \ref{DE1} above.
  \bigskip

  We examine the equalities (\ref{DE1}) more closely. The equality
  holds for every $V\in\frak k$ implies that $h(t)\inv  h'(t)  +
  a'(t) a\inv(t)$ does not have any vertical component. That is,
  $-a'(t) a\inv(t)$ is the vertical component of $h(t)\inv  h'(t)$
  so that
  $$
    h(t)\inv  h'(t)  = - a'(t) a\inv(t) + X_1
    \in\frakk \oplus \frakk^\perp.
  $$
  is a vertical and horizontal splitting.

  Let $g(t)$ be another path with a unique horizontal lift $\tilde
  g(t)=g(t) b(t)$, satisfying
  \begin{align}
    \label{DE3}
   0&=\angles{g\inv  g' +   b' b\inv, V},
  \end{align}
  for every $V\in \frakk$. Again, we have a splitting
  $$
    g(t)\inv  g'(t)  = - b'(t) b\inv(t) + X_2
    \in\frakk \oplus \frakk^\perp.
  $$
  From
  $$
    ||h(t)\inv  h'(t)  - g(t)\inv  g'(t) ||
    = || a'(t) a\inv(t) - b'(t) b\inv(t)||+ ||X_1 - X_2||,
  $$
  we get
  \begin{align}
    \label{DE4}
    ||a'(t) a\inv(t)-b'(t) b\inv(t)||\leq ||h(t)\inv  h'(t)-g(t)\inv  g'(t)||.
  \end{align}
  These are norms on the Lie algebra $\mathfrak{so}(1,n)$.
  \bigskip

  On the space of piecewise $C^k$-curves ($k\geq 1$)
  in $\text{ SO}_0 (1,n)$ with
  initial point $e$, we define a distance function by
  $$
    \rho(h,g)= \int_0^1 ||h'(t)\cdot h(t)\inv - g'(t)\cdot g(t)\inv||\ dt.
  $$
  Note that $h'(t)\cdot h(t)\inv \in\fsono$ and $||.||$ is the norm
  there. We argue that this is a metric. Suppose $\rho(h,g)=0$.
  Then, by continuity, $h'(t)\cdot h(t)\inv = g'(t)\cdot g(t)\inv$
  for every $t$. Now we apply the following Lemma to the
  $C^1$-curves piece by piece to conclude $h(t)=g(t)$ for all $t\in
  [0,1]$, see [3], vol 1, p69.

  \begin{lemma} 
    Let $G$ be a Lie group and $\frak g$ its Lie algebra identified
    with $T_e(G)$. Let $Y_t$, $0\leq t\leq 1$, be a continuous curve
    in $T_e(G)$. Then there exists in $G$ a unique curve $a_t$ of
    class $C^1$ such that $a_0=e$ and $\dot a_t a_t\inv= Y_t$ for
    $0\leq t\leq 1$.
  \end{lemma}

  Let $h$ be a curve in $\bbh^n$ (or in $\na$, by abuse of
  notation). The unique curve $a: [0,1]\ra \so(n)$ such that
  $h(t)\cdot a(t)$ is the horizontal lift of $h(t)$ will be called
  $w_h$.

  For two curves $h$ and $g$, the inequality (\ref{DE4}) shows that
  $\rho(w_h,w_g)\leq \rho(h,g)$. Let $\frak P$ be the space of all
  piecewise $C^k$-curves on $\na$ with the initial point $e$.

  \begin{proposition} \label{converge-lift}
    The map $\frak P\lra G$ sending $f$ to $w_h(1)$ is continuous.
    More precisely, let $h: [0,1]\ra \na$ be a piecewise $C^k$-curve.
    For every $\epsilon>0$, there exists $\delta>0$ such that, if
    $g\in\frak P$ and $\rho(h,g)<\delta$, then
    $d(e, w_h(1)^{-1} \cdot w_g(1))=d(w_h(1),w_g(1))<\epsilon$.
  \end{proposition}
  \bigskip

  \begin{proof}
    For simplicity, we write $w_h(t)$, $w_g(t)$ by $a(t)$, $b(t)$,
    respectively.
    \begin{align}
      \label{DAA}
      0=(a a^T)'=a' a^T + a (a^T)'
    \end{align}
    Clearly $a(a^t b)' b^T$ is a well-defined element of the Lie
    algebra $\frak k$ and
    \begin{align*}
      a(a^T b)' b^T
      &=a\left((a^T)' b + a^T  b'\right )b^T\\
      &=a (a^T)' + b' b^T\\
      &=-a' a^T + b' b^T \quad\text{(from the equality (\ref{DAA}))}\\
      &=-a' a\inv + b' b\inv
    \end{align*}
  Thus,
  $$
    ||a' a\inv - b' b\inv|| = ||a(a^T b)' b^T||.
  $$
  Observe that $(a^T b)'\in T_{a^T b}(K)$. The left translation
  $\ell_a$ and the right translation $r_{b^T}$ maps this vector to a
  tangent vector at $0\in T_e(K)$. However, both these translations
  are isometries so that they preserve the norms. We have,
  $$
    ||a' a\inv - b' b\inv|| = ||a(a^T b)' b^T|| = ||(a^T b)'|| = ||(a\inv b)'||.
  $$
  Consequently, if
  $$
    \int ||(a\inv b)'||\ dt = \int ||a' a\inv - b' b\inv||\ dt
  $$
  is small, the arc-length of the path $a(t)\inv b(t)$ is small.
  Therefore, if $a(0)$ and $b(0)$ were close (or if $a(0)=b(0)$),
  then $a(1)$ and $b(1)$ are close. This finishes the proof.
  \end{proof}

  Since
  $$
    \bar{f}_n (t_0) =
    \text{ the value , at } t=1 ,
    \text{ of the horizontal lifting of } \hat{\gamma} ^n_{t_0}
    \text{ at } e
  $$
  and
  $$
    \hat{\gamma} ^n_{t_0} \text{ converges to }
    \tilde{\gamma} ^{n_0}_{t_0}
    = \tilde{\gamma} _{t_0} \text{ as n goes to } \infty \: ,
  $$
  we get, by reparameterizing $\tilde{\gamma} ^{n_0}_{t_0}$ and
  $\tilde{\gamma} _{t_0}$
  respectively if needed,

 \begin{align*}
   \bar{f}(t_0)
   & = \text{ the value , at } t=1 ,
       \text{ of the horizontal lifting of } \tilde{\gamma} _{t_0}
       \text{ at } e  \\
   & = \lim _{n \rightarrow \infty} \bar{f}_n (t_0) .
 \end{align*}

\bigskip
  \subsection{Preliminaries for the main proof}

  \bigskip

  Fix
  $ t_0 \in \displaystyle{\bigcup ^{\infty}_{n=1}} D_n - \{ 0 \} $
  and find a positive integer
  $ n_0 $ =
  min $\{ n_1   \mid n+1 \geq n_1 \Rightarrow  t_0 \in D_n \} \: .$

  Assume $ n \geq n_0 $ .

  Note $ t_0 $ is not the last element in $ D_n $ for $ n \geq n_0 .$
  Notice that with respect to totally geodesic triangles,
  $
    \hat{c} ^n_{t^n_2 (t_0)} (1)
    = {_1 \hat{c} ^n_{t_0} (1)}
    = {_1 \hat{c} ^{n_0}_{t_0}(1)}
    = \hat{c}^{n_0}_{t^{n_0}_2(t_0)}(1)
  $
  for all $ n \geq n_0 ,$  which is the ending point of
  the $ j_n (t_0) $-th triangle in $ \hat{A} _n $ and also the
  starting point of the $ (j_n (t^n_2 (t_0)) = j_n (t_0) + 1) $-th
  triangle in $ \hat{A} _n $ for all $ n \geq n_0 .$

  \subsubsection{ \textbf{
                    A new curve $ \mathbf{\hat{c} ^{short}_{t_0}}$
                    for the comparison of triangles
                   }
                 }

  Define
  $ \hat{c} ^{short}_{t_0} : [0,1] \rightarrow \mathbb{H} ^n $ as
  the shortest geodesic from $ \pi (e) \in \mathbb{H} ^n $ to
  $
    \hat{c} ^{n}_{t^{n}_2 (t_0)} (1)
    = \hat{c} ^{n_0}_{t^{n_0}_2 (t_0)} (1)
    = {_1 \hat{c} ^{n_0}_{t_0} (1)}
    = {_1 \hat{c} ^{n}_{t_0} (1)} ,
  $
  in other words, to the starting point of
  $ (j_n (t^n_2 (t_0)) = j_n (t_0) +1) $
  -th triangle in $ \hat{A}_n$, which is also the ending point of
  $j_n (t_0)$-th triangle in $\hat{A}_n$.
  Consider its horizontal lift
  $$
    _e \hat{c} ^{short}_{t_0} :
    [0,1] \rightarrow \text{SO}_0(1,n)
  $$
  at
  $e .$

  \subsubsection{ \textbf{
                  the comparison of
                  $\mathbf{(j_n (t^n_2 (t_0)) = j_n (t_0) + 1)}$-th
                  totally geodesic triangles
                 }
                }
  For each $ n \geq n_0 $, consider
  $$
    _e \hat{c} ^n_{t^n_2 (t_0)} :=
    R_{(\hat{f}_n (t_0))^{-1}} \circ
    {_{\hat{f}_n (t_0)} \hat{c} ^n_{t^n_2 (t_0)}} :
    [0,1] \rightarrow \text{SO}_0(1,n) \: ,
  $$
  which is the horizontal lifting of $ \hat{c} ^n_{t^n_2 (t_0)} $ at
  $e$, that is ,
  $$
    \pi \circ {_e \hat{c} ^n_{t^n_2 (t_0)}}
    = \; \hat{c} ^n_{t^n_2 (t_0)}
    = \; \pi \circ {_{\hat{f}_n (t_0)} \hat{c} ^n_{t^n_2 (t_0) } } \: .
  $$

  Note
  $
    _e \hat{c} ^n_{t^n_2 (t_0)}
    \text{ and } \;
    {_{\hat{f}_n (t_0)} \hat{c} ^n_{t^n_2 (t_0)}}
  $
  are piecewise geodesics, since   the right translation
  $ R_k : G \rightarrow G $ by $k$ is an isometry for any
  $ k \in K = \text{SO}(n) $ and $ \hat{c} ^n_{t_0} $ are piecewise
  geodesics.

  \begin{figure}[h]
       \centering{\includegraphics[width=2.2in]{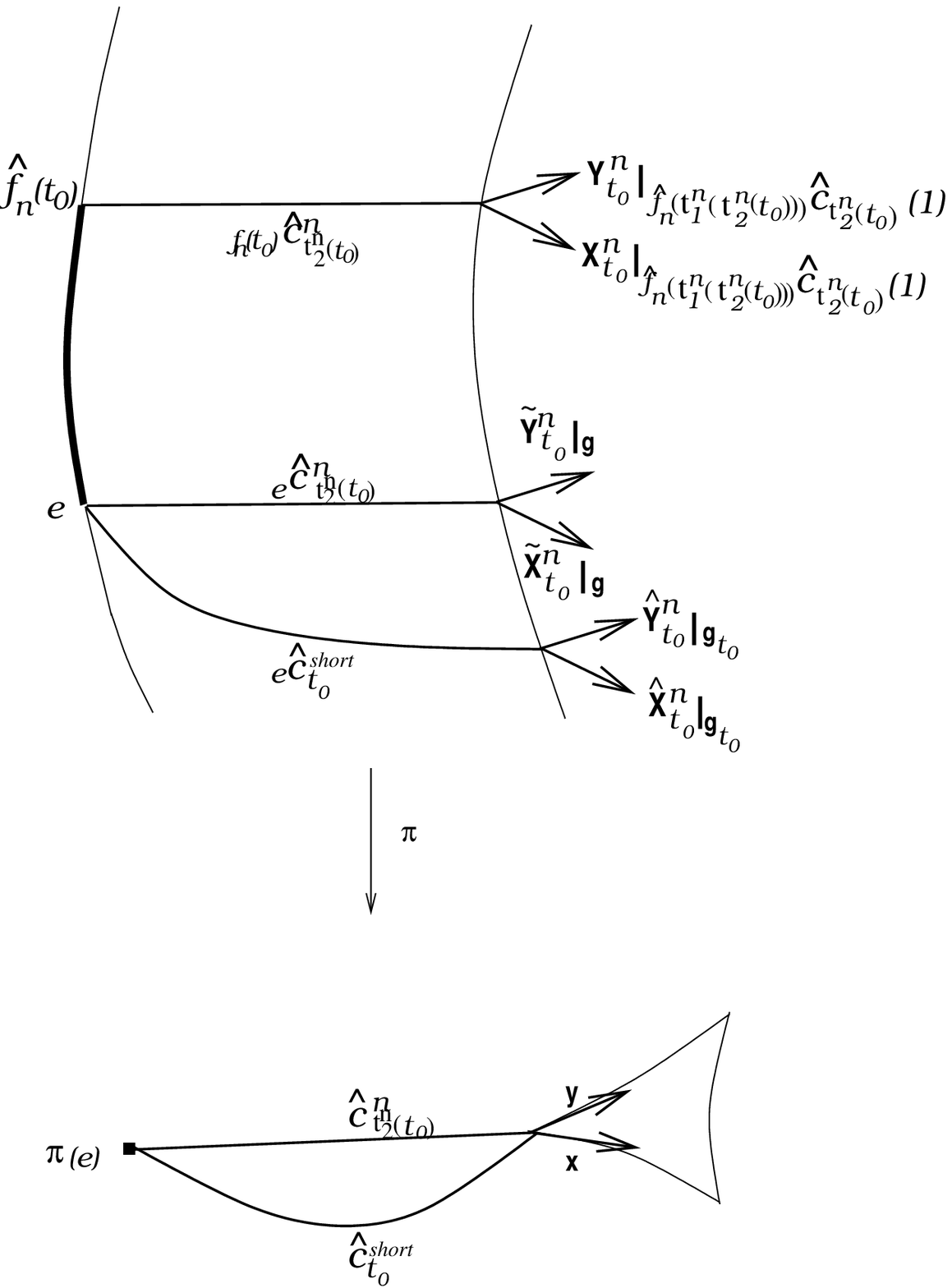}}
  \end{figure}

  Consider the $ (j_n (t^n_2 (t_0)) = j_n (t_0) + 1) $-th triangle in
  $ \hat{A} _n $, its starting point and the horizontal lifting of
  $$
    x :=
    \lim _{t \rightarrow 0^{+}}
    \frac{1}{\mid \dot{ \hat{\varphi} } ^n_{t^n_2 (t_0)} (t) \mid }
    \cdot
    \dot{ \hat{\varphi} } ^n_{t^n_2 (t_0)} (t)
  $$
  and
  $$
    y :=
    - \lim _{t \rightarrow 1^{-}}
    \frac{1}{\mid \dot{ \hat{\varphi} } ^n_{t^n_2 (t_0)} (t) \mid }
    \cdot
    \dot{ \hat{\varphi} } ^n_{t^n_2 (t_0) } (t)
  $$
  at ${ _e \hat{c} ^n_{t^n_2 (t_0)} (1) } = : g $ , respectively
  and find
  $$
    \tilde{X} ^n_{t_0}  , \tilde{Y} ^n_{t_0}
    \in \mathfrak{k}^{ \bot}
  $$
  with
  $$ \pi _{*}  \tilde{X} ^n_{t_0} \mid _{g} = x  $$
  and
  $$  \pi _{*} \tilde{Y} ^n_{t_0} \mid _{g} = y  \: .$$
  Also consider the horizontal lifting of $ x $ and $ y $ at
  $ _e \hat{c} ^{short}_{t_0} (1) = : g_{t_0} $ and find
  $$
    \hat{X} ^n_{t_0}  , \hat{Y} ^n_{t_0}  \in \mathfrak{k}^{ \bot}
  $$
  with
  $$  \pi _{*} \hat{X} ^n_{t_0} \mid _{g_{t_0}} = x  $$
  and
  $$  \pi _{*} \hat{Y} ^n_{t_0} \mid _{g_{t_0}} = y  \: .$$
  Note
   $$
     _{\hat{f}_n (t_0)} \hat{c} ^n_{t^n_2 (t_0)}  \; =
     \; _{\hat{f}_n (t^n_1 (t^n_2 (t_0)))} \hat{c} ^n_{t^n_2 (t_0)} \: ,
   $$
  so

  \begin{align*}
    \pi _{*} \; X^n_{t_0} \mid_
    {_{\hat{f}_n (t_0)} \hat{c} ^n_{t^n_2 (t_0)} (1)}
    & = \pi _{*} \; X^n_{t^n_1 (t^n_2 (t_0))} \mid_
        {_{\hat{f}_n (t^n_1 (t^n_2 (t_0)))} \hat{c}^n_{t^n_2 (t_0)} (1)} \\
    & = \lim _{t \rightarrow 0^{+}}
        \frac{1}{\mid \dot{\hat{\varphi}}^n_{t^n_2 (t_0)} (t) \mid }
        \cdot
        \dot{ \hat{\varphi} } ^n_{t^n_2 (t_0)} (t) \\
    & = x \\
    & = \pi _{*} \tilde{X}^n_{t_0} \mid _g \\
    & = \pi _{*} \tilde{X}^n_{t_0} \mid_
        {_e \hat{c} ^n_{t^n_2 (t_0)} (1)} \quad ,
  \end{align*}

  \noindent
  which implies

  $$ X^n_{t_0} = Ad_{(\hat{f}_n (t_0))^{-1}} \tilde{X} ^n_{t_0} $$
  from
  $
    {_{\hat{f}_n (t_0)} \hat{c} ^n_{t^n_2 (t_0)}} =
    R_{{\hat{f}_n (t_0)}} \circ {_e \hat{c} ^n_{t^n_2 (t_0)}} .
  $

  Similarly,
  $$ Y^n_{t_0} = Ad_{(\hat{f}_n (t_0))^{-1}} \tilde{Y} ^n_{t_0} \: .$$
  And by considering a loop
  $$
    \bar{\hat{c}} ^{short} _{t_0} * \hat{c} ^n_{t^n_2 (t_0)} :
    [0,1] \rightarrow \mathbb{H} ^n \: ,
  $$
  where
  $ \bar{\hat{c}} ^{short}_{t_0} : [0,1] \rightarrow \mathbb{H}^3 $
  is given by
  $$
    \bar{\hat{c}} ^{short}_{t_0} (t) = \hat{c} ^{short}_{t_0} (1-t)
    \: ,
  $$
  and its horizontal lifting at $ _e \hat{c} ^{short}_{t_0} (1) $ ,
  we obtain
  $$
    \tilde{X} ^n_{t_0} =
    Ad
    _{
      (
        ({_e \hat{c} ^{short}_{t_0} (1)})^{-1} \cdot
        {_e \hat{c} ^n_{t^n_2 (t_0)} (1)}
      )^{-1}
    }
    \hat{X} ^n_{t_0}
  $$
  $$
    \tilde{Y} ^n_{t_0} =
    Ad
    _{
      (
        ({_e \hat{c} ^{short}_{t_0} (1)})^{-1} \cdot
        {_e \hat{c} ^n_{t^n_2 (t_0)} (1)}
      )^{-1}
    }
    \hat{Y} ^n_{t_0} \: .
    $$
  Then we get
  $$
    X^n_{t_0} =
    Ad
    _{
      (
        ({_e \hat{c} ^{short}_{t_0} (1)})^{-1} \cdot
        {_e \hat{c} ^n_{t^n_2 (t_0)} (1)} \cdot
        \hat{f}_n (t_0)
      )^{-1}
    }
    \hat{X} ^n_{t_0}
  $$
  $$
    Y^n_{t_0} =
    Ad
    _{
      (
        ({_e \hat{c} ^{short}_{t_0} (1)})^{-1} \cdot
        {_e \hat{c} ^n_{t^n_2 (t_0)} (1)} \cdot
        \hat{f}_n (t_0)
      )^{-1}
    }
    \hat{Y} ^n_{t_0} \: .$$
  Since both
  $
    ({_e \hat{c} ^{short}_{t_0} (1)})^{-1} \cdot
    {_e \hat{c} ^n_{t^n_2 (t_0)} (1)}
  $
  and $ \hat{f}_n (t_0) $ are elements in $ K = \text{SO}(n) $, we get
  $$
    [X^n_{t_0}, Y^n_{t_0}] =
    Ad
    _{
      (
        ({_e \hat{c} ^{short}_{t_0} (1)})^{-1} \cdot
        {_e \hat{c} ^n_{t^n_2 (t_0)} (1)} \cdot
        \hat{f}_n (t_0)
      )^{-1}
    }
    [\hat{X} ^n_{t_0}, \hat{Y}^n_{t_0}] \: .
  $$

  Note $ \hat{c} ^n_{t^n_2 (t_0) } = {_1 \hat{c} ^n_{t_0} } $ .
  So we can rewrite $ X^n_{t_0} $ and $ Y^n_{t_0} $ as
  $$
    X^n_{t_0} =
    Ad
    _{
      (
        ({_e \hat{c} ^{short}_{t_0} (1)})^{-1} \cdot
        {_e (_1 \hat{c}) ^n_{t_0} (1)} \cdot
        \hat{f}_n (t_0)
      ) ^{-1}
    }
    \hat{X} ^n_{t_0}
  $$
  $$
    Y^n_{t_0} =
    Ad
    _{
      (
        ({_e \hat{c} ^{short}_{t_0} (1)})^{-1} \cdot
        {_e (_1 \hat{c}) ^n_{t_0} (1)} \cdot
        \hat{f}_n (t_0)
      ) ^{-1}
    }
    \hat{Y} ^n_{t_0} \: .
  $$

  Since ${_1 \hat{c} ^n_{t_0}} $ converges to
  $ {_1 \tilde{c} _{t_0} } $ in $\mathbb{H}^n$ as $n$,
  the number of steps (not the dimension of $\mathbb{H}^n$) goes to
  $\infty$, ${_e (_1 \hat{c}) ^n_{t_0}} $ converges to
  ${_e (_1 \tilde{c}) _{t_0} } $ in $ \text{SO}_0(1,n) \: .$

  Since
  $$
    \hat{f}_n (t_0) =
    \text{ the value , at } t=1 ,
    \text{ of the horizontal lifting of } \hat{\gamma} ^n_{t_0}
    \text{ at } e
  $$
  and
  $$
    \hat{\gamma} ^n_{t_0} \text{ converges to }
    \tilde{\gamma} ^{n_0}_{t_0}
    = \tilde{\gamma} _{t_0} \text{ as n goes to } \infty \: ,
  $$
  we get, from Proposition \ref{converge-lift},

 \begin{align*}
   \bar{f}(t_0)
   & = \text{ the value , at } t=1 ,
       \text{ of the horizontal lifting of } \tilde{\gamma} _{t_0}
       \text{ at } e  \\
   & = \lim _{n \rightarrow \infty} \hat{f}_n (t_0) .
 \end{align*}

  Then we also get

  $$
    _{\hat{f}_n (t_0)} (_1 \hat{c}) ^n_{t_0}
    = R_{\hat{f}_n (t_0)} \circ {_e (_1 \hat{c}) ^n_{t_0}}
    = {_e (_1 \hat{c}) ^n_{t_0}} \cdot \hat{f}_n(t_0),
  $$
  which will converge to
  $$
    _e (_1 \tilde{c}) _{t_0} \cdot \bar{f}(t_0)
    = R_{\bar{f}(t_0)} \circ {_e (_1 \tilde{c}) _{t}}
    = {_{\bar{f}(t_0)} (_1 \tilde{c}) _{t_0}} \: .
  $$

  \subsubsection{ \textbf{
                  the comparison of
                  $\mathbf{(j_n (t^n_2 (t_0)) = j_n (t_0) + 1)}$-th
                  triangles on the given surface $\mathbf{S}$
                 }
                }

  \begin{figure}[h]
       \centering{\includegraphics[width=1.93in]{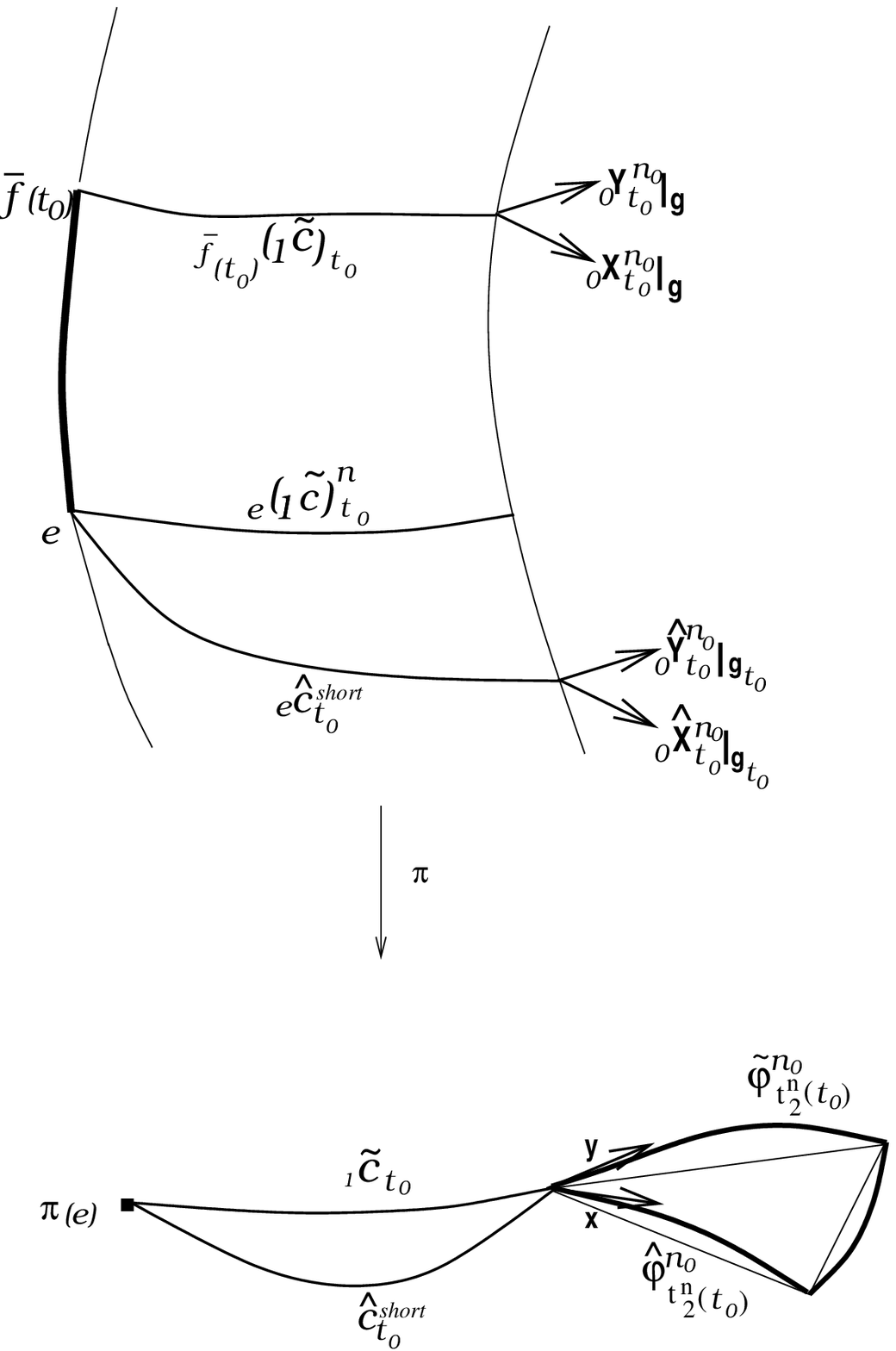}}
  \end{figure}
  Now, consider the
  $(j_{n_0} (t^{n_0}_2 (t_0)) = j_{n_0} (t_0) + 1) $-th triangle,
  lying on $S$, in $\tilde{A}_{n_0}$, its starting point and the
  horizontal lifting of
  $$
    x :=
    \lim _{t \rightarrow 0^{+}}
    \frac{1}
         {\mid \dot{\tilde{\varphi}}^{n_0}_{t^{n_0}_2 (t_0)}(t)\mid}
    \cdot
    \dot{\tilde{\varphi}}^{n_0}_{t^{n_0}_2 (t_0)} (t)
  $$
  and
  $$
    y := - \lim _{t \rightarrow 1^{-}}
    \frac{1}
         {\mid \dot{\tilde{\varphi} }^{n_0}_{t^{n_0}_2 (t_0)}(t)\mid}
    \cdot
    \dot{ \tilde{\varphi} } ^{n_0}_{t^{n_0}_2 (t_0)} (t)
  $$
  at $ g_{t_0} = {_e \hat{c} ^{short}_{t_0} (1)} $ and
  at $ g := {_{\bar{f}(t_0) } (_1  \tilde{c})_{t_0} (1)} $, respectively,
  and find
  $$
    _0 \hat{X}^{n_0}_{t_0} , \; {_0 \hat{Y}^{n_0}_{t_0}}, \; ,
    {_0 X^{n_0}_{t_0}}     , \; {_0 Y^{n_0}_{t_0}}
    \in \mathfrak{k}^{ \bot}
  $$
  with
  $$
    \pi _{*} \: {{_0 \hat{X}^{n_0}_{t_0}} \mid _{ g_{t_0} }}
    = x
    = \pi _{*} \text{ } {_0 X^{n_0}_{t_0}} \mid _g
  $$
  and
  $$
    \pi _{*} \: {{_0 \hat{Y}^{n_0}_{t_0}} \mid _{ g_{t_0} }}
    = y
    = \pi _{*} \text{ } {_0 Y^{n_0}_{t_0}} \mid _g  \: .
  $$
  Then,
  $$
    g = {_{\bar{f}(t_0) } (_1  \tilde{c})_{t_0} (1)}
    = {_e (_1  \tilde{c})_{t_0} (1)} \cdot \bar{f}(t_0)
    = g_{t_0} \cdot ({_e \hat{c} ^{short}_{t_0} (1)})^{-1 }
              \cdot {_e (_1  \tilde{c})_{t_0} (1)} \cdot \bar{f}(t_0)
  $$
  implies that
  $$
    _0 X^{n_0}_{t_0} =
    Ad
    _{
      (
        ({_e \hat{c} ^{short}_{t_0} (1)})^{-1}
        \cdot {_e (_1 \tilde{c}) _{t_0} (1)}
        \cdot \bar{f}(t_0)
      ) ^{-1}
    }
    \;
    {_0 \hat{X}} ^{n_0}_{t_0}
  $$

  $$
    _0 Y^{n_0}_{t_0} =
    Ad
    _{
      (
        ({_e \hat{c} ^{short}_{t_0} (1)})^{-1}
        \cdot {_e (_1 \tilde{c}) _{t_0} (1)}
        \cdot \bar{f}(t_0)
      ) ^{-1}
    }
    \;
    {_0 \hat{Y}} ^{n_0}_{t_0}
    \: .
  $$

  \subsubsection{\textbf{
      The convergence of tangent planes induced by
      $\mathbf{(j_n (t^n_2 (t_0)) = j_n (t_0) + 1)}$-th triangles at
      $\mathbf{t=t_0}$
      and the convergence of $\mathbf{\hat{f}_n}$ under
      $\mathbf{
                \displaystyle{
                    \lim _{n \rightarrow \infty}
                    \lim_{t \rightarrow {t_0}^{+}}}
         }
      $
   }
  }
  \begin{figure}[h]
       \centering{\includegraphics[width=1.78in]{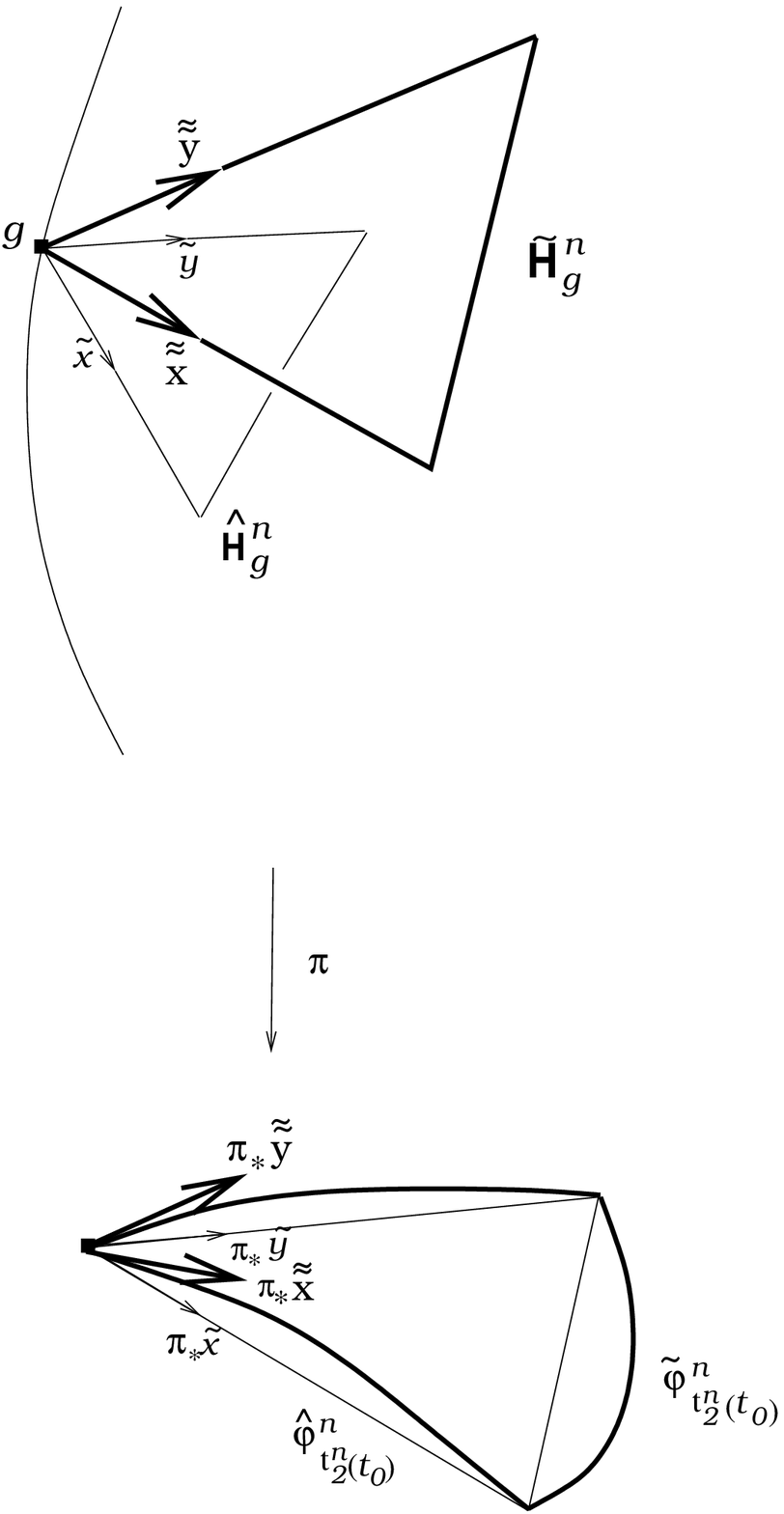}}
  \end{figure}
  Now, for
  $ g \in \pi ^{-1} (\pi (_{\hat{f}_n (t_0)} (_1 \hat{c}) ^n_{t_0} (1))) ,$
  let
  $$  \hat{H}^n_g := \text{Span} \{ \tilde{x}, \tilde{y} \} \: ,$$
  where $\tilde{x}, \tilde{y}$ are horizontal vectors at $g$ satisfying
  $$
    \pi_{*} \; \tilde{x} =
    \lim _{t \rightarrow 0^{+}}
    \frac{1}{\mid \dot{ \hat{\varphi} } ^n_{t^n_2(t_0)} (t) \mid }
    \cdot
    \dot{ \hat{\varphi} } ^n_{t^n_2(t_0)} (t)
  $$
  $$
    \pi_{*} \; \tilde{y} =
    - \lim _{t \rightarrow 1^{-}}
    \frac{1}{\mid \dot{ \hat{\varphi} } ^n_{t^n_2(t_0)} (t) \mid }
    \cdot
    \dot{ \hat{\varphi} } ^n_{t^n_2(t_0)} (t)
    \: .
  $$
  Also, for
  $ g \in \pi ^{-1} (\pi (_{\hat{f}_n (t_0)} (_1 \hat{c}) ^n_{t_0} (1))) ,$
  let
  $$
    \tilde{H}^n_g :=
    \text{Span} \{ \tilde{\tilde{x}}, \tilde{\tilde{y}} \} \: ,
  $$
  where $\tilde{\tilde{x}}, \tilde{\tilde{y}}$ are horizontal vectors
  at $g$ satisfying
  $$
    \pi_{*} \; \tilde{\tilde{x}} =
    \lim _{t \rightarrow 0^{+}}
    \frac{1}{\mid \dot{ \tilde{\varphi} } ^n_{t^n_2(t_0)} (t) \mid }
    \cdot
    \dot{ \tilde{\varphi} } ^n_{t^n_2(t_0)} (t)
  $$
  $$
    \pi_{*} \; \tilde{\tilde{y}} =
    - \lim _{t \rightarrow 1^{-}}
    \frac{1}{\mid \dot{ \tilde{\varphi} } ^n_{t^n_2(t_0)} (t) \mid }
    \cdot
    \dot{ \tilde{\varphi} } ^n_{t^n_2(t_0)} (t) .
  $$

 \bigskip

  Note, for
  $
    t \in
   (t_0, \; t^n_2 (t_0)) = (t^n_1 (t^n_2 (t_0)), \; t^n_2 (t_0)) ,
  $

  \medskip
  \noindent
  $
    \omega \;
    (\; \tfrac{1}{\mid \dot{\hat{f}}_n (t) \mid} \cdot \dot{\hat{f}} _n (t) \;)
  $
  \begin{align*}
    = \: (-1) \cdot
      \big(\;
        & \text
            {
              the unit curvature of the 2-dimensional horizontal
              oriented tangent
            }
        \\
        & \text{plane,}
        \\
        & \hat{H}^n_{_{\hat{f}_n (t_0)} (_1 \hat{c}) ^n_{t_0} (1)} \: = \:
          \hat{H}^n_{
                      _{\hat{f}_n (t^n_1 (t^n_2(t_0)))}
                      \hat{c} ^n_{t^n_2(t_0)}
                      (1)
                    }
        \\
        & =
          \text{Span}
          \{
            X^n_{t^n_1 (t^n_2 (t_0))} \mid
            _{_{\hat{f}_n (t^n_1 (t^n_2(t_0)))}\hat{c}^n_{t^n_2(t_0)} (1)},
            \;
            Y^n_{t^n_1 (t^n_2(t_0))} \mid
            _{_{\hat{f}_n (t^n_1 (t^n_2(t_0)))}\hat{c}^n_{t^n_2(t_0)} (1)}
          \}
        \\
        & =
          \text{Span}
          \{
            X^n_{t_0} \mid
            _ {_{\hat{f}_n (t_0)} ({_1 \hat{c}}) ^n_{t_0} (1)} , \;
            Y^n_{t_0} \mid
            _ {_{\hat{f}_n (t_0)} ({_1 \hat{c}}) ^n_{t_0} (1)}
          \} \:
        \\
        & \text{at }
          {_{\hat{f}_n (t^n_1 (t^n_2(t_0)))} \hat{c} ^n_{t^n_2(t_0)} (1)}
          = {_{\hat{f}_n (t_0)} (_1 \hat{c}) ^n_{t_0} (1)} \; ,
        \\
        & \text{which projects to the tangent plane}
        \\
        & \text{ \: \; of the } (j_n(t^n_2 (t_0)) = j_n (t_0) +1)
          \text{-th triangle in } \hat{A}_n
        \\
        & \text{at }  \pi ({_{\hat{f}_n (t_0)} (_1 \hat{c}) ^n_{t_0} (1)})
          = {_1 \hat{c}^n_{t_0} (1)} = \hat{c}^n_{t^n_2 (t_0)} (1)
        \\
        & \text{ \: \; = the starting point of the }
          (j_n(t^n_2 (t_0)) = j_n (t_0) +1)
          \text{-th triangle in } \hat{A}_n
        \\
        & \text{
                 with respect to the connection of the principal
                 bundle
                }
           \pi : \text{SO}_0 (1,n) \rightarrow \mathbb{H} ^n
      \big)
  \end{align*}
  $
    = \tfrac{1}
            {
              \mid
                [
                  X^n_{t^n_1 (t^n_2(t_0))}, \;
                  Y^n_{t^n_1 (t^n_2(t_0))}
                ]
              \mid
            }
      \cdot  \;
      [ X^n_{t^n_1 (t^n_2(t_0))}, \; Y^n_{t^n_1 (t^n_2(t_0))} ]
    \\
    = \tfrac{1}
            { \mid [ X^n_{t_0}, \; Y^n_{t_0} ]  \mid }
      \cdot  \;
      [ X^n_{t_0}, \text{ } Y^n_{t_0} ]
    \\
    = \frac
        {
          Ad
            _{
               (
                 ({_e \hat{c} ^{short}_{t_0} (1)})^{-1}
                 \cdot {_e (_1 \hat{c}) ^n_{t_0} (1)}
                 \cdot \hat{f}_n (t_0)
               )^{-1}
             }
            [\hat{X} ^n_{t_0}, \hat{Y}^n_{t_0}]
        }
        {
          \mid
            Ad
              _{
                 (
                   ({_e \hat{c} ^{short}_{t_0} (1)})^{-1}
                   \cdot {_e  (_1 \hat{c}) ^n_{t_0} (1)}
                   \cdot \hat{f}_n (t_0)
                 )^{-1}
               }
              [\hat{X} ^n_{t_0}, \hat{Y}^n_{t_0}]
          \mid
        }
    \\
    = Ad
        _{
           (
             ({_e \hat{c} ^{short}_{t_0} (1)})^{-1}
             \cdot {_e (_1 \hat{c}) ^n_{t_0} (1)}
             \cdot \hat{f}_n (t_0)
           )^{-1}
         }
        \Big(
          \tfrac{1}
                {\mid [\hat{X} ^n_{t_0}, \hat{Y}^n_{t_0}] \mid}
          \cdot
          [\hat{X} ^n_{t_0}, \hat{Y}^n_{t_0}]
        \Big)
  $
  \begin{align*}
    = (-1) \cdot
      Ad
        _{
           (
             ({_e \hat{c} ^{short}_{t_0} (1)})^{-1}
             \cdot {_e (_1 \hat{c}) ^n_{t_0} (1)}
             \cdot \hat{f}_n (t_0)
           )^{-1}
         }
        \big( \;
          &\text{the unit curvature of the 2-dimensional }
          \\
          & \text{horizontal oriented tangent plane, }
  \end{align*}
          \begin{flushright}
             $
               \hat{H}^n_{_e \hat{c} ^{short}_{t_0} (1)}
               = \text{Span}
                 \{
                   \hat{X}^n_{t_0}\mid _{_e \hat{c}^{short}_{t_0}(1)}
                   ,
                   \hat{Y}^n_{t_0}\mid _{_e \hat{c}^{short}_{t_0}(1)}
                 \}
               \big)
             $
          \end{flushright}

  Note the tangent plane of the
  $(j_n(t^n_2 (t_0)) = j_n(t_0) + 1)$-th triangle in $\hat{A}_n  $ at
  $
    _ 1 \hat{c} ^n_{t_0} (1)
    = \hat{c} ^n_{t^n_2 (t_0)} (1)
    = \hat{c} ^{n_0}_{t^n_2 (t_0)} (1)
    = {_1 \hat{c} ^{n_0}_{t_0} (1)}
    = \tilde{c} ^{n_0}_{t^n_2 (t_0)} (1)
    = {_1 \tilde{c} ^{n_0}_{t_0} (1)}
  $
  for all $ n \geq n_0 \; ,$    the starting point of the
  $ (j_n (t^n_2 (t_0)) = j_n (t_0) + 1) $-th triangle in
  $ \hat{A} _n $ and the ending point of the $ j_n (t_0) $-th
  triangle in $ \hat{A} _n $ for all $ n \geq n_0 $ at the same time,
  which is also the starting point of the
  $ (j_{n_0} (t^{n_0}_2 (t_0)) = j_{n_0} (t_0) + 1) $-th triangle,
  lying on S,  in $ \tilde{A} _{n_0} $ and the ending point of the
  $ j_{n_0} (t_0) $-th triangle in $ \tilde{A} _{n_0} $ at the same
  time, will converge to the tangent plane of $S$ at
  $
    \tilde{c} ^{n_0}_{t^n_2 (t_0)} (1)
    = {_1 \tilde{c} ^{n_0}_{t_0} (1)} \: .
  $

  \begin{figure}[h]
       \centering{\includegraphics[width=3in]{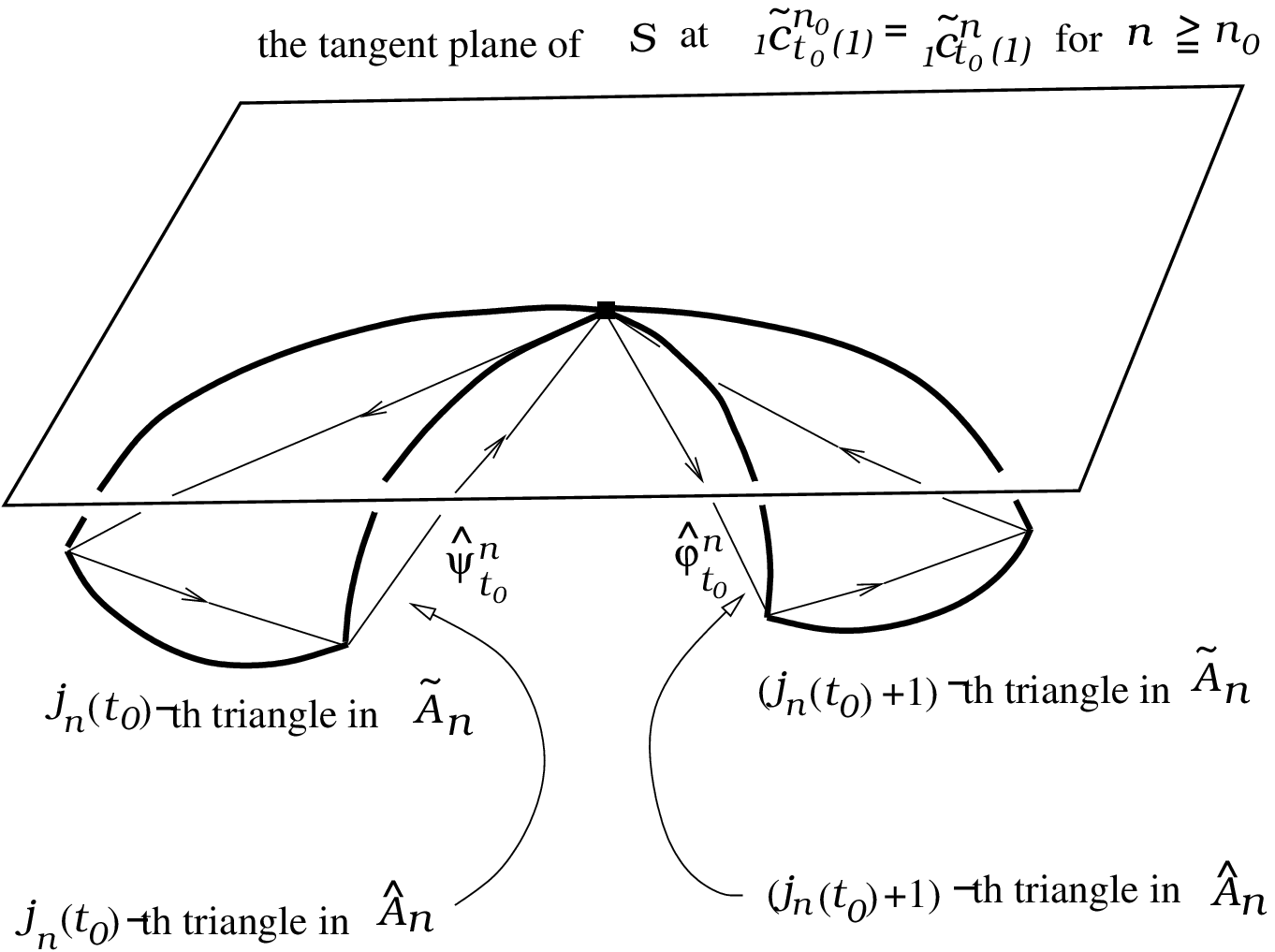}}
  \end{figure}

  And, note, in general, if $ \lim_{n \rightarrow \infty} g_n = g_{_0} $
  in $G$ and $ \lim_{n \rightarrow \infty} X_n = X_0 $ in
  $ \mathfrak{g} \: ,$ then

  \begin{align*}
    \lim_{n \rightarrow \infty} t \cdot Ad_{g_n} X_n
    & = \lim_{n \rightarrow \infty}
      \text{exp} ^{-1} (\text{exp} (t \cdot Ad_{g_n} x_n))
    \\
    & = \lim_{n \rightarrow \infty}
        \text{exp} ^{-1}
        (g_n \cdot \text{exp}(t \cdot X_n ) \cdot {g_{n}} ^{-1} )
    \\
    & = \text{exp} ^{-1}
        (g_{_0} \cdot \text{exp}(t \cdot X_0 ) \cdot {g_{_0}} ^{-1})
    \\
    & = \text{exp} ^{-1} (\text{exp} (t \cdot Ad_{g_{_0}} X_0))
    \\
    & = t \cdot Ad_{g_{_0}} X_0 \: .
  \end{align*}

  Now, refer to previous three pictures. Then we get, for
  $
    t \in \big( t_0, t^n_2 (t_0) \big) \quad ,
  $

  \medskip
  \noindent
  $
    \omega \;
    (\; \tfrac{1}{\mid \dot{\hat{f}}_n (t) \mid} \cdot \dot{\hat{f}} _n (t) \;)
    \\
    = \tfrac{1}
            { \mid [ X^n_{t_0}, \; Y^n_{t_0} ]  \mid }
      \cdot  \;
      [ X^n_{t_0}, \; Y^n_{t_0} ]
    \\
    = Ad
        _{
           (
             ({_e \hat{c} ^{short}_{t_0} (1)})^{-1}
             \cdot {_e (_1 \hat{c}) ^n_{t_0} (1)}
             \cdot \hat{f}_n (t_0)
           )^{-1}
         }
        \Big(
          \tfrac{1}
                {\mid [\hat{X} ^n_{t_0}, \hat{Y}^n_{t_0}] \mid}
          \cdot
          [\hat{X} ^n_{t_0}, \hat{Y}^n_{t_0}]
        \Big)
  $
  \begin{align*}
    = (-1) \cdot
      Ad
        _{
           (
             ({_e \hat{c} ^{short}_{t_0} (1)})^{-1}
             \cdot {_e (_1 \hat{c}) ^n_{t_0} (1)}
             \cdot \hat{f}_n (t_0)
           )^{-1}
         }
        \big( \;
          &\text{the unit curvature of the 2-dimensional }
          \\
          & \text{horizontal oriented tangent plane, }
  \end{align*}
          \begin{flushright}
             $
               \hat{H}^n_{_e \hat{c} ^{short}_{t_0} (1)}
               = \text{Span}
                 \{
                   \hat{X}^n_{t_0}\mid _{_e \hat{c}^{short}_{t_0}(1)}
                   ,
                   \hat{Y}^n_{t_0}\mid _{_e \hat{c}^{short}_{t_0}(1)}
                 \}
               \big),
             $
          \end{flushright}

  \medskip
  \noindent
  which will converge to

  \begin{align*}
    (-1) \cdot
    Ad
      _{
         (
           ({_e \hat{c} ^{short}_{t_0} (1)})^{-1}
           \cdot {_e (_1 \tilde{c})_{t_0} (1)}
           \cdot \bar{f}(t_0)
         )^{-1}
       }
      \big( \;
        &\text{the unit curvature of the 2-dimensional }
        \\
        & \text{horizontal oriented tangent plane, }
  \end{align*}
        \begin{flushright}
           $
             \tilde{H}^{n_0}_{_e \hat{c} ^{short}_{t_0} (1)}
               = \text{Span}
               \{
                  {_0 \hat{X}}^{n_0}_{t_0}\mid
                     _{_e \hat{c}^{short}_{t_0}(1)}
                  , \;
                  {_0 \hat{Y}}^{n_0}_{t_0}\mid
                     _{_e \hat{c}^{short}_{t_0}(1)}
               \}
             \big)
           $
        \end{flushright}
  \noindent
  $
    = Ad
        _{
           (
             ({_e \hat{c} ^{short}_{t_0} (1)})^{-1}
             \cdot {_e (_1 \tilde{c})_{t_0} (1)}
             \cdot \bar{f}(t_0)
           )^{-1}
         }
        \Big(
          \tfrac{1}
                {\mid
                  [
                    {_0 \hat{X}}^{n_0}_{t_0}
                    , \;
                    {_0 \hat{Y}}^{n_0}_{t_0}
                  ]
                \mid}
          \cdot
          [{_0 \hat{X}}^{n_0}_{t_0},\; {_0 \hat{Y}}^{n_0}_{t_0}]
        \Big)
    \\
    = \tfrac{1}
            { \mid [{_0 X}^{n_0}_{t_0},\; {_0 Y}^{n_0}_{t_0}]  \mid }
      \cdot  \;
      [{_0 X}^{n_0}_{t_0},\; {_0 Y}^{n_0}_{t_0}]
  $
  \noindent
  \begin{align*}
    = \: (-1) \cdot
      \big(\;
        & \text
            {
              the unit curvature of the 2-dimensional horizontal
              oriented tangent
            }
        \\
        & \text{plane, } \:
          \tilde{H}^{n_0}_{_{\bar{f}(t_0)} (_1 \tilde{c})_{t_0} (1)}
          \: = \:
          \text{Span}
          \{
            {_0 X}^{n_0}_{t_0} \mid
            _ {_{\bar{f}(t_0)} ({_1 \tilde{c}})_{t_0} (1)} , \;
            {_0 Y}^{n_0}_{t_0} \mid
            _ {_{\bar{f}(t_0)} ({_1 \tilde{c}})_{t_0} (1)}
          \}
        \\
        & \text{at } {_{\bar{f}(t_0)} ({_1 \tilde{c}})_{t_0} (1)} \; ,
        \\
        & \text{which projects to the tangent plane}
        \\
        & \text{ \quad of the }
          (j_{n_0}(t^{n_0}_2 (t_0)) = j_{n_0} (t_0) +1)
          \text{-th triangle in } \tilde{A}_{n_0} \; ,
        \\
        & \text{\hspace{1cm}  where } \quad
          n_0  =
          \text{ min }
          \{ n_1   \mid n+1 \geq n_1 \Rightarrow  t_0 \in D_n \}
          \quad ,
        \\
        & \text{ \quad - so tangent to the given disk $S$ -}
        \\
        & \text{at }
          \pi (_{\bar{f}(t_0)} ({_1 \tilde{c}})_{t_0} (1))
          = {_1 \tilde{c}_{t_0} (1)}
          = \tilde{c}^{n_0}_{t^{n_0}_2 (t_0)} (1)
        \\
        & \text{ \: \; = the starting point of the }
          (j_{n_0}(t^{n_0}_2 (t_0)) = j_{n_0} (t_0) +1)
          \text{-th triangle in } \tilde{A}_{n_0}
        \\
        & \text{\noindent
                 with respect to the connection of the principal
                 bundle
                }
           \pi : \text{SO}_0 (1,n) \rightarrow \mathbb{H} ^n
      \big)
  \end{align*}

  \noindent
  under
  $
    \displaystyle{
                   \lim _{n \rightarrow \infty}
                   \lim_{t \rightarrow {t_0}^{+}}} \: .
  $

  So, under the identification of $ T_e K $ with $ \mathfrak{k} \; ,$
  \\
  $
    \displaystyle{
                   \lim_{n \rightarrow \infty}
                   \lim_{t \rightarrow {t_0}^{+}}
                 }
       {L_{(\hat{f}_n(t)^{-1})}}_{*} \;
       \tfrac{1}{\mid \dot{\hat{f}} _n (t)  \mid}
       \cdot \dot{\hat{f}} _n (t)
    \\
    = \: \displaystyle{\lim_{n \rightarrow \infty}}
            \frac{1}{\mid  [X^n_{t_0} ,  Y^n_{t_0} ] \mid}
            \cdot [ X^n_{t_0} ,  Y^n_{t_0} ]
    \\
    = \: \frac{1}{
                   \mid
                     [{_0 X^{n_0}_{t_0}} , {_0 Y^{n_0}_{t_0} }]
                   \mid
                 }
         \cdot  [{_0 X^{n_0}_{t_0}} , {_0 Y^{n_0}_{t_0} }]
  $
  \noindent
  \begin{align*}
    = \: (-1) \cdot
      \big(\;
        & \text
            {
              the unit curvature of the 2-dimensional horizontal
              oriented tangent
            }
        \\
        & \text{plane, } \:
          \tilde{H}^{n_0}_{_{\bar{f}(t_0)} (_1 \tilde{c})_{t_0} (1)}
          \text{ at } {_{\bar{f}(t_0)} ({_1 \tilde{c}})_{t_0} (1)} \; ,
        \\
        & \text{which projects to the tangent plane of } S
          \text{ at }
          \pi (_{\bar{f}(t_0)} ({_1 \tilde{c}})_{t_0} (1))
          = {_1 \tilde{c}_{t_0} (1)}
        \\
        & \text{\noindent
                 with respect to the connection of the principal
                 bundle
                }
           \pi : \text{SO}_0 (1,n) \rightarrow \mathbb{H} ^n
      \big)
  \end{align*}

  \bigskip

  \subsubsection{
    \textbf{
    $
     \mathbf{
      \lim_{n \rightarrow \infty} \lim_{t \rightarrow {t_0}^{+}}
        {L_{(\hat{f}_n(t)^{-1})}}_{*} \;
        \tfrac{1}{\mid \dot{\hat{f}} _n (t)  \mid}
        \cdot \dot{\hat{f}} _n (t) \:
      = \: \lim_{n \rightarrow \infty} \lim_{t \rightarrow {t_0}^{-}}
              {L_{(\hat{f}_n(t)^{-1})}}_{*} \;
              \tfrac{1}{\mid \dot{\hat{f}} _n (t)  \mid}
              \cdot \dot{\hat{f}} _n (t)
      }
    $
   }
   }
  To show
  $$
    \lim_{n \rightarrow \infty} \lim_{t \rightarrow {t_0}^{+}}
      {L_{(\hat{f}_n(t)^{-1})}}_{*} \;
      \tfrac{1}{\mid \dot{\hat{f}} _n (t)  \mid}
      \cdot \dot{\hat{f}} _n (t) \:
    = \: \lim_{n \rightarrow \infty} \lim_{t \rightarrow {t_0}^{-}}
            {L_{(\hat{f}_n(t)^{-1})}}_{*} \;
            \tfrac{1}{\mid \dot{\hat{f}} _n (t)  \mid}
            \cdot \dot{\hat{f}} _n (t) \: ,
  $$
  for
  $
    g \in \pi ^{-1}
          (\pi (_{\hat{f}_n (t^n_1(t_0))} (_1 \hat{c}) ^n_{t_0} (1)))
    = \pi ^{-1} (\pi (_{\hat{f}_n (t_0)} (_1 \hat{c}) ^n_{t_0} (1))) \; ,
  $
  let
  $$  _1 \hat{H}^n_g := \text{Span} \{ \tilde{x}, \tilde{y} \} \: ,$$
  where $\tilde{x}, \tilde{y}$ are horizontal vectors at $g$
  satisfying
  $$
    \pi_{*} \; \tilde{x} =
    \lim _{t \rightarrow 0^{+}}
       \frac{1}{\mid \dot{ \hat{\psi} } ^n_{t_0} (t) \mid }
       \cdot
       \dot{ \hat{\psi} } ^n_{t_0} (t)
  $$
  $$
    \pi_{*} \; \tilde{y} =
    - \lim _{t \rightarrow 1^{-}}
         \frac{1}{\mid \dot{ \hat{\psi} } ^n_{t_0} (t) \mid }
         \cdot
         \dot{ \hat{\psi} } ^n_{t_0} (t) \: .
  $$
  Also, for
  $
    g \in \pi ^{-1}
          (\pi (_{\hat{f}_n (t^n_1(t_0))} (_1 \hat{c}) ^n_{t_0} (1))) \: ,
  $
  let
  $$
    _1 \tilde{H}^n_g :=
    \text{Span} \{ \tilde{\tilde{x}}, \tilde{\tilde{y}} \} \: ,
  $$
  where $\tilde{\tilde{x}}, \tilde{\tilde{y}}$ are horizontal vectors
  at $g$ satisfying
  $$
    \pi_{*} \; \tilde{\tilde{x}} =
    \lim _{t \rightarrow 0^{+}}
       \frac{1}{\mid \dot{ \tilde{\psi} } ^n_{t_0} (t) \mid }
       \cdot
       \dot{ \tilde{\psi} } ^n_{t_0} (t)
  $$
  $$
    \pi_{*} \; \tilde{\tilde{y}} =
    - \lim _{t \rightarrow 1^{-}}
         \frac{1}{\mid \dot{ \tilde{\psi} } ^n_{t_0} (t) \mid }
         \cdot
         \dot{ \tilde{\psi} } ^n_{t_0} (t) \: .
  $$

  Now,  consider the horizontal lifting of
  $$
    z :=
    \lim _{t \rightarrow 0^{+}}
       \frac{1}{\mid \dot{ \hat{\psi} } ^n_{t_0} (t) \mid }
       \cdot
       \dot{ \hat{\psi} } ^n_{t_0} (t)
  $$
  and
  $$
    w :=
    - \lim _{t \rightarrow 1^{-}}
         \frac{1}{\mid \dot{ \hat{\psi} } ^n_{t_0} (t) \mid }
         \cdot
         \dot{ \hat{\psi} } ^n_{t_0} (t)
  $$
  at
  $ g := {_{\hat{f}_n (t^n_1 (t_0)) } {(_1 \hat{c})} ^n_{t_0} (1)}
       = { _{\hat{f}_n (t^n_1 (t_0)) } \hat{c} ^n_{t^n_2 (t_0)} (1) } \; ,
  $
  respectively and find
  $$ Z^n_{t_0}  , \text{ } W^n_{t_0}  \in \mathfrak{k}^{ \bot} $$
  with
  $$ \pi _{*} \text{ } Z^n_{t_0} \mid _{ g } = z  $$
  and
  $$ \pi _{*} \text{ } W^n_{t_0} \mid _{ g } = w  \: .$$

  Also consider the horizontal lifting of $ z $ and $ w $
  at
  $
    {_e (_1 \hat{c}) ^n_{t_0}} (1)
    = { _e \hat{c} ^n_{t^n_2 (t_0)} (1) }
    = : g \; ,
  $
  respectively
  and find
  $$
    \tilde{Z} ^n_{t_0}  , \; {\tilde{W}} ^n_{t_0}
    \in \mathfrak{k}^{ \bot}
  $$
  with
  $$ \pi _{*} \; \tilde{Z} ^n_{t_0} \mid _{g} = z  $$
  and
  $$  \pi _{*} \; {\tilde{W}} ^n_{t_0} \mid _{g} = w  \: .$$
  And consider the horizontal lifting of $ z $ and $ w $ at
  $ _e \hat{c} ^{short}_{t_0} (1) = : g_{t_0} $ and find
  $$
    \hat{Z} ^n_{t_0}  , \; {\hat{W}} ^n_{t_0}
    \in \mathfrak{k}^{ \bot}
  $$
  with
  $$ \pi _{*} \; {\hat{Z}} ^n_{t_0} \mid _{g_{t_0}} = z  $$
  and
  $$ \pi _{*} \; {\hat{W}} ^n_{t_0} \mid _{g_{t_0}} = w  \: .$$
  Note $ \text{Im} \varphi^n_{t_0} = \text{ Im} \psi^n_{t_0} $ is the
  boundary of a geodesic triangle in $ \mathbb{H}^n \: .$

 Then, from Facts, mentioned earlier in this section, and from the
 property in Section \ref{sec:n=2}, we get

  \begin{align*}
    \hat{f}_n (t)
    & = \hat{f}_n(t^n_1 (t_0)) \cdot
        \text{ exp }
          \Big(
             (t-t^n_1 (t_0))
             \cdot
             \frac
                {
                  (
                    \text{ Area of } j_n (t_0)
                    \text{-th triangle in } \hat{A} _n
                  )
                }
                {
                  ( (t_0 - t^n_1 (t_0)) \cdot
                    \mid
                      [X^n_{t^n_1(t_0)} , Y^n_{t^n_1 (t_0)} ]
                    \mid
                  )
                }
             \cdot
             [X^n_{t^n_1 (t_0)} , Y^n_{t^n_1 (t_0)} ]
          \Big) \\
    & = \hat{f}_n(t^n_1 (t_0)) \cdot
        \text{ exp }
          \Big(
             (t -t^n_1 (t_0)) \cdot
             \frac
               {
               (
                 \text{ Area of } j_n (t_0)
                 \text{-th triangle in } \hat{A} _n )
               }
               {
                 (
                   (t_0 - t^n_1 (t_0)) \cdot
                   \mid  [Z^n_{t_0} , W^n_{t_0} ] \mid
                 )
               }
             \cdot
             [Z^n_{t_0} , W^n_{t_0} ]
          \Big) \\
  \end{align*}

  \noindent
  for $ t \in [t^n_1 (t_0), t_0 ] \: .$

  \bigskip

  Note

  \begin{align*}
    \pi _{*} \; Z^n_{t_0} \mid
      _{_{\hat{f}_n (t^n_1  (t_0))} (_1 \hat{c}) ^n_{t_0} (1)}
   & = \lim _{t \rightarrow 0^{+}}
          \frac{1}{\mid \dot{ \hat{\psi} } ^n_{t_0} (t) \mid }
          \cdot
          \dot{ \hat{\psi} } ^n_{t_0} (t)
   \\
   & = z
   \\
   & = \pi _{*} \;
         \tilde{Z}^n_{t_0} \mid _{_e {(_1 \hat{c})} ^n_{t_0} (1)}
   \; ,
  \end{align*}

  \noindent
  which implies

  $$ Z^n_{t_0} = Ad_{(\hat{f}_n (t^n_1 (t_0)))^{-1}} \tilde{Z} ^n_{t_0} $$
  from
  $
    {_{\hat{f}_n (t^n_1 (t_0))} {(_1 \hat{c})} ^n_{t_0}} =
    R_{{\hat{f}_n (t^n_1 (t_0))}} \circ {_e (_1 \hat{c}) ^n_{t_0}} \: .
  $

  Similarly,
  $$
    W^n_{t_0} = Ad_{(\hat{f}_n {t^n_1 (t_0)})^{-1}} \tilde{W} ^n_{t_0} \: .
  $$

  And by considering a loop
  $$
    \bar{\hat{c}} ^{short} _{t_0} * {_1 \hat{c}} ^n_{t_0} :
    [0,1] \rightarrow \mathbb{H} ^n \: ,
  $$
  where
  $ \bar{\hat{c}} ^{short}_{t_0} : [0,1] \rightarrow \mathbb{H}^n $
  is given by
  $$
    \bar{\hat{c}} ^{short}_{t_0} (t) = \hat{c} ^{short}_{t_0} (1-t)
    \: ,
  $$
  and its horizontal lifting at $ _e \hat{c} ^{short}_{t_0} (1) \: ,$
  we obtain
  $$
    \tilde{Z} ^n_{t_0} =
    Ad
      _{
         (
           ({_e \hat{c} ^{short}_{t_0} (1)})^{-1} \cdot
           {_e (_1 \hat{c}) ^n_{t_0} (1) }
         )^{-1}
       }
      \hat{Z} ^n_{t_0}
  $$
  $$
    \tilde{W} ^n_{t_0} =
    Ad
      _{
         (
           ({_e \hat{c} ^{short}_{t_0} (1)})^{-1} \cdot
           {_e (_1 \hat{c}) ^n_{t_0} (1)}
         )^{-1}
       }
      \hat{W} ^n_{t_0} \: .$$
  Then we get
  $$
    Z^n_{t_0} =
    Ad
      _{
         (
           ({_e \hat{c} ^{short}_{t_0} (1)})^{-1}
           \cdot
           {_e (_1 \hat{c}) ^n_{t_0} (1)}
           \cdot
           \hat{f}_n (t^n_1 (t_0))
         )^{-1}
       }
      \hat{Z} ^n_{t_0}
  $$
  $$
    W^n_{t_0} =
      Ad
        _{
           (
             ({_e \hat{c} ^{short}_{t_0} (1)})^{-1}
             \cdot
             {_e (_1 \hat{c}) ^n_{t_0} (1)}
             \cdot
             \hat{f}_n (t^n_1 (t_0))
           )^{-1}
         }
        \hat{W} ^n_{t_0} \: .
  $$

  Since both
  $
    ({_e \hat{c} ^{short}_{t_0} (1)})^{-1} \cdot
    {_e (_1 \hat{c}) ^n_{t_0} (1)}$
  and $ \hat{f}_n (t^n_1 (t_0)) $ are elements in $ K = \text{SO}(n) \; ,$
  we get
  $$
    [Z^n_{t_0}, W^n_{t_0}] =
     Ad
       _{
         (
           ({_e \hat{c} ^{short}_{t_0} (1)})^{-1}
           \cdot
           {_e (_1 \hat{c}) ^n_{t_0} (1)}
           \cdot
           \hat{f}_n (t^n_1 (t_0))
         )^{-1}
        }
       [\hat{Z} ^n_{t_0}, \hat{W}^n_{t_0}] \: .
  $$

  \bigskip

  Note the tangent plane of the $j_n(t_0) $-th triangle in $\hat{A}_n
  $ at $ _1 \hat{c} ^n_{t_0} (1) = \hat{c} ^n_{t^n_2 (t_0)} (1) =
  \hat{c} ^{n_0}_{t^n_2 (t_0)} (1) = {_1 \hat{c} ^{n_0}_{t_0} (1)} =
  \tilde{c} ^{n_0}_{t^n_2 (t_0)} (1) = {_1 \tilde{c} ^{n_0}_{t_0}
  (1)} $ for all $ n \geq n_0 $ ,   the starting point of the $ (j_n
  (t^n_2 (t_0)) = j_n (t_0) + 1) $-th triangle in $ \hat{A} _n $ and
  also  the ending point of the $ j_n (t_0) $-th triangle in $
  \hat{A} _n $ for all $ n \geq n_0 $, which is also the starting
  point of the $ (j_{n_0} (t^{n_0}_2 (t_0)) = j_{n_0} (t_0) + 1) $-th
  triangle, lying on S,  in $ \tilde{A} _{n_0} $ and also  the ending
  point of the $ j_{n_0} (t_0) $-th triangle in $ \tilde{A} _{n_0} $,
  will converge to the tangent plane of $S$ at $ \tilde{c}
  ^{n_0}_{t^n_2 (t_0)} (1) = {_1 \tilde{c} ^{n_0}_{t_0} (1)} $, which
  implies that for $t \in \big( t^n_1(t_0), t_0 \big)$

  \medskip
  \noindent
  $
    \omega \;
    (\; \tfrac{1}{\mid \dot{\hat{f}}_n (t) \mid} \cdot \dot{\hat{f}} _n (t) \;)
    \\
    = \tfrac{1}
            { \mid [ Z^n_{t_0}, \; W^n_{t_0} ]  \mid }
      \cdot  \;
      [ Z^n_{t_0}, \; W^n_{t_0} ]
    \\
    = Ad
        _{
           (
             ({_e \hat{c} ^{short}_{t_0} (1)})^{-1}
             \cdot {_e (_1 \hat{c}) ^n_{t_0} (1)}
             \cdot \hat{f}_n (t^n_1(t_0))
           )^{-1}
         }
        \Big(
          \tfrac{1}
                {\mid [\hat{Z} ^n_{t_0}, \hat{W}^n_{t_0}] \mid}
          \cdot
          [\hat{Z} ^n_{t_0}, \hat{W}^n_{t_0}]
        \Big)
  $
  \begin{align*}
    = (-1) \cdot
      Ad
        _{
           (
             ({_e \hat{c} ^{short}_{t_0} (1)})^{-1}
             \cdot {_e (_1 \hat{c}) ^n_{t_0} (1)}
             \cdot \hat{f}_n (t^n_1(t_0))
           )^{-1}
         }
        \big( \;
          &\text{the unit curvature of the 2-dimensional }
          \\
          & \text{horizontal oriented tangent plane, }
  \end{align*}
          \begin{flushright}
             $
               {_1 \hat{H}}^n_{_e \hat{c} ^{short}_{t_0} (1)}
               = \text{Span}
                 \{
                   \hat{Z}^n_{t_0}\mid _{_e \hat{c}^{short}_{t_0}(1)}
                   ,
                   \hat{W}^n_{t_0}\mid _{_e \hat{c}^{short}_{t_0}(1)}
                 \}
               \big)
             $
          \end{flushright}

  \medskip
  \noindent
  will converge to

  \begin{align*}
    (-1) \cdot
    Ad
      _{
         (
           ({_e \hat{c} ^{short}_{t_0} (1)})^{-1}
           \cdot {_e (_1 \tilde{c})_{t_0} (1)}
           \cdot \bar{f}(t_0)
         )^{-1}
       }
      \big( \;
        &\text{the unit curvature of the 2-dimensional }
        \\
        & \text{horizontal oriented tangent plane, }
  \end{align*}
        \begin{flushright}
           $
             \tilde{H}^{n_0}_{_e \hat{c} ^{short}_{t_0} (1)}
             = \text{Span}
               \{
                  {_0 \hat{X}}^{n_0}_{t_0}\mid
                     _{_e \hat{c}^{short}_{t_0}(1)}
                  , \;
                  {_0 \hat{Y}}^{n_0}_{t_0}\mid
                     _{_e \hat{c}^{short}_{t_0}(1)}
               \}
             \big)
           $
        \end{flushright}
  \noindent
  $
    = Ad
        _{
           (
             ({_e \hat{c} ^{short}_{t_0} (1)})^{-1}
             \cdot {_e (_1 \tilde{c})_{t_0} (1)}
             \cdot \bar{f}(t_0)
           )^{-1}
         }
        \Big(
          \tfrac{1}
                {\mid
                  [
                    {_0 \hat{X}}^{n_0}_{t_0}
                    , \;
                    {_0 \hat{Y}}^{n_0}_{t_0}
                  ]
                \mid}
          \cdot
          [{_0 \hat{X}}^{n_0}_{t_0},\; {_0 \hat{Y}}^{n_0}_{t_0}]
        \Big)
    \\
    = \tfrac{1}
            { \mid [{_0 X}^{n_0}_{t_0},\; {_0 Y}^{n_0}_{t_0}]  \mid }
      \cdot  \;
      [{_0 X}^{n_0}_{t_0},\; {_0 Y}^{n_0}_{t_0}]
  $
  \noindent
  \begin{align*}
    = \: (-1) \cdot
      \big(\;
        & \text
            {
              the unit curvature of the 2-dimensional horizontal
              oriented tangent
            }
        \\
        & \text{plane, } \:
          \tilde{H}^{n_0}_{_{\bar{f}(t_0)} (_1 \tilde{c})_{t_0} (1)}
          \: = \:
          \text{Span}
          \{
            {_0 X}^{n_0}_{t_0} \mid
            _ {_{\bar{f}(t_0)} ({_1 \tilde{c}})_{t_0} (1)} , \;
            {_0 Y}^{n_0}_{t_0} \mid
            _ {_{\bar{f}(t_0)} ({_1 \tilde{c}})_{t_0} (1)}
          \}
        \\
        & \text{at } {_{\bar{f}(t_0)} ({_1 \tilde{c}})_{t_0} (1)} \; ,
        \\
        & \text{which projects to the tangent plane}
        \\
        & \text{ \quad of the }
          (j_{n_0}(t^{n_0}_2 (t_0)) = j_{n_0} (t_0) +1)
          \text{-th triangle in } \tilde{A}_{n_0} \; ,
        \\
        & \text{\hspace{1cm}  where } \quad
          n_0  =
          \text{ min }
          \{ n_1   \mid n+1 \geq n_1 \Rightarrow  t_0 \in D_n \}
          \quad ,
        \\
        & \text{ \quad - so tangent to the given disk $S$ -}
        \\
        & \text{at }
          \pi (_{\bar{f}(t_0)} ({_1 \tilde{c}})_{t_0} (1))
          = {_1 \tilde{c}_{t_0} (1)}
          = \tilde{c}^{n_0}_{t^{n_0}_2 (t_0)} (1)
        \\
        & \text{ \: \; = the starting point of the }
          (j_{n_0}(t^{n_0}_2 (t_0)) = j_{n_0} (t_0) +1)
          \text{-th triangle in } \tilde{A}_{n_0}
        \\
        & \text{\noindent
                 with respect to the connection of the principal
                 bundle
                }
           \pi : \text{SO}_0 (1,n) \rightarrow \mathbb{H} ^n
      \big)
  \end{align*}

  Thus we get

  \begin{align*}
    \lim_{n \rightarrow \infty} \lim_{t \rightarrow {t_0}^{-}}
       {L_{(\hat{f}_n(t)^{-1})}}_{*} \;
       \Big(
         \tfrac{1}{\mid \dot{\hat{f}} _n (t)  \mid} \cdot \dot{\hat{f}} _n (t)
       \Big)
    & = \lim_{n \rightarrow \infty}
           \frac{1}{\mid  [Z^n_{t_0} , W^n_{t_0} ] \mid} \cdot
           [Z^n_{t_0} , W^n_{t_0} ]
    \\
    & = \; \frac{1}
                {
                  \mid
                    [{_0 X^{n_0}_{t_0}} , {_0 Y^{n_0}_{t_0} }]
                  \mid
                }
           \cdot
           [{_0 X^{n_0}_{t_0}} , {_0 Y^{n_0}_{t_0} }]
    \\
    & = \lim_{n \rightarrow \infty} \lim_{t \rightarrow {t_0}^{+}}
          {L_{(\hat{f}_n(t)^{-1})}}_{*} \;
          \Big(
             \tfrac{1}{\mid \dot{\hat{f}} _n (t)  \mid} \cdot \dot{\hat{f}} _n (t)
          \Big)
    \quad .
  \end{align*}

  \vspace{1cm}
  \subsubsection{\textbf{Main Part}}

  Define a function
  $ s_n : D_n - \{0\} \rightarrow (0, \infty) $  as follows : \\
  Given $ t \in D_n - \{0\} $, assume $t$ is the
  $j$-th element in $D_n$, $i.e.$, $ j=j_n (t) \: .$ Then,

  \begin{align*}
    s_n (t) :
    & = \sum^{j=j_n(t)}_{i=1}
           ( \text{ the area of i-th triangle in } \hat{A}_n )
    \\
    & = \text{ the area of the region surrounded by }
        \hat{\gamma}^n_t
        \text { in the } n \text{-th step pleated surface} .
  \end{align*}

  Note in $S$, for $n \geq n_0$,
  $$
    \text{the region surrounded by } \tilde{\gamma} ^n_{t_0}
    \text{ in } S
    =
    \text{ the region surrounded by } \tilde{\gamma} ^{n_0}_{t_0}
    \text{ in } S \; ,
  $$
  so we get
  $$
    \lim_{n \rightarrow \infty} s_n (t_0)  =
    \text{ the area of the region surrounded by }
    \tilde{\gamma} ^{n_0}_{t_0}
    \text{ in } S  =: s(t_0) \: .
  $$
  Thus, we obtain a function
  $$
    s : \cup ^{\infty}_{n=1} D_n - \{ 0 \} \rightarrow (0, \infty)
    \: .
  $$
  Now, induce a function
  $$
    {f_n} \; : \;
    [ 0, \text{ the area of the } n \text{-step pleated surface } ]
    \rightarrow K \; ,
  $$
  which is the reparametrization of $ \hat{f}_n $ with
  $ \mid \! {\dot{f}_n} (t) \! \; \mid = 1 $ on \\
  $
    [0, \text{ the area of the } n \text{-step pleated surface } ]
    \; - \;
    \\
    \Big\{ \:
       \displaystyle{\sum ^{j}_{i=1}}
       (\text{ the area of the } i \text{-th triangle in } \hat{A}_n )
       \: \mid \:
       j=1,2,\cdots, \: \mid \! \hat{A}_n  \! \mid
    \Big\}
    \: .
  $
  \\
  Then we get
  $$
    {f_n} (s_n (t)) \; = \; \hat{f}_n (t) = \bar{f}_n (t) \qquad
    \text{ for } t \in D_n - \{ 0 \} \: .
  $$
  Define a function
  $$
    {f} :
    \{ s(t) \mid t \in \cup ^{\infty}_{n=1} D_n - \{ 0 \} \}
    \rightarrow K
  $$
  by
  $$ {f} (s(t)) = \bar{f}(t) \: .$$
  Then we get \\
  $$
    {f} (s(t_0))
    = \bar{f}(t_0)
    = \lim _{n \rightarrow \infty} \bar{f}_n (t_0)
    = \lim _{n \rightarrow \infty} {f_n} (s_n (t_0)) \: .
  $$
  Note, for $ t_1 \in \cup^{\infty}_{n=1} D_n - \{ 0 \} \: ,$

  \begin{align*}
    \bar{f}(t_1)
    & = \lim _{n \rightarrow \infty} \bar{f}_n (t_1)
    \\
    & = \lim _{n \rightarrow \infty}
          (
            \text{the value, at } t=1,
            \text{ of the horizontal lifting of  }
            \hat{\gamma}^n_{t_1} \text{ at } e
          )
    \\
    & = \text{ the value, at } t=1 ,
        \text{ of the horizontal lifting of }
        \tilde{\gamma}_{t_1}  \text{ at } e \: .
  \end{align*}

  Since $\tilde{\gamma}_{t_1}$ converges to $\tilde{\gamma}_{t_0}$
  as $t_1$ approaches $t_0$ in $ \cup^{\infty}_{n=1} D_n$, Proposition \ref{converge-lift}
  implies that $\bar{f}$ will be continuous on
  $ \displaystyle{\bigcup ^{\infty}_{n=1}} D_n - \{ 0 \} $
  and we can extend $\bar{f}$ on [0,1]. And from
  $ {f} (s(t_0)) = \bar{f}(t_0) \; ,$
  ${f}$ will be continuous on
  $
    \big\{
       s(t)
       \: \mid \:
       t \in  \displaystyle{\bigcup^{\infty}_{n=1}} D_n - \{ 0 \}
    \big\}
    \: .
  $

  Note $s$ is continuous on
  $ \displaystyle{\bigcup ^{\infty}_{n=1}} D_n - \{ 0 \} $
  and so it can be extended on [0,1].

  Since
  $
    \big\{
      s(t)
      \: \mid t \:
      \in \displaystyle{\bigcup ^{\infty}_{n=1}} D_n - \{0\}
    \big\}
  $
  is a dense subset of $[0,$ the area of $ S] $, we can
  extend $f$ on $[0,$ the area of $ S] $ continuously. Call it
  $ f $ as well. Then we get
  $$  {f} \circ s = \bar{f} \text{ is continuous on } [0,1]  $$
  and
  $$
    f(\text{ the area of } S)
    = \bar{f}(1)
    = \lim_{t \rightarrow 1, t \in \bigcup ^{\infty}_{n=1} D_n }
       {_e \tilde{\gamma}} _t (1)
    = {_e \tilde{\gamma} } (1),
  $$
  \noindent
  where $ \tilde{\gamma} : [0,1] \rightarrow S $ is the boundary
  curve of S and $ _e \tilde{\gamma} $ is its horizontal lifting at
  $e$.

  \bigskip

  Let's show $ f $ is a $ C^1 $ curve.

  Define a function      $  F_n $ from \\
  $
    [0, \text{ the area of the } n \text{-step pleated surface } ]
    \; - \;
    \\
    \Big\{ \;
       \displaystyle{\sum ^{j}_{i=1}}
       (\text{ the area of the } i \text{-th triangle in } \hat{A}_n )
       \; \mid \;
       j=1,2,\cdots, {\mid \! \hat{A}_n \! \mid \;}
    \Big\}
  $  \\
  to
  $$ \text {the unit sphere in } \mathfrak{k} $$
  by
  $$
    F_n (t)  \: = \:
    {L_{(f_n (t)) ^{-1}}}_{*} \; {\dot{f} _n} (t)  \: .
  $$

  And define a function
  $$
    F :
    \big\{
       s(t)
       \: \mid \:
       t \in \displaystyle{\bigcup ^{\infty}_{n=1}} D_n - \{ 0 \}
    \big\}
    \rightarrow \text{ the unit sphere in } \mathfrak{k}
  $$
  by
  $$
   F (s(t_0)) \:
   = \: \lim_{n \rightarrow \infty } \lim_{t \rightarrow {t_0} ^{-}}
       {L_{(\hat{f}_n(t)^{-1})}}_{*} \;
         \Big(
           \tfrac{1}{\mid \dot{\hat{f}} _n (t) \mid} \cdot \dot{\hat{f}} _n (t)
         \Big)
   = \: \lim_{n \rightarrow \infty } \lim_{t \rightarrow {t_0} ^{+}}
       {L_{(\hat{f}_n(t)^{-1})}}_{*} \;
         \Big(
           \tfrac{1}{\mid \dot{\hat{f}} _n (t) \mid} \cdot \dot{\hat{f}} _n (t)
         \Big)
   \: .
  $$
  Then, $ F_n $ is constant on the interval
  $$(0, \text{ the area of the first triangle in } \hat{A}_n)$$
  and on the interval
  $$
    \bigg(
       \displaystyle{\sum ^{j}_{i=1}}
         (
           \text{\small the area of the }
           i \text{\small -th triangle in } \hat{A}_n
         )
       \; , \:
       \displaystyle{\sum ^{j+1}_{i=1}}
         (
           \text{\small the area of the }
           i \text{\small -th triangle in } \hat{A}_n
         )
    \bigg)
  $$
  for each $ j=1,2,\cdots, {\mid \! \hat{A}_n \! \mid} \; , \:  $  and

  \begin{align*}
    F(s(t_0))
    & = \; \lim_{n \rightarrow \infty }
           \lim_{t \rightarrow {t_0} ^{-}}
              {L_{(\hat{f}_n(t)^{-1})}}_{*} \;
                \Big(
                    \tfrac{1}{\mid \dot{\hat{f}} _n (t) \mid}
                    \cdot
                    \dot{\hat{f}} _n (t)
                \Big)
           = \;
             \lim_{n \rightarrow \infty }
             \lim_{t \rightarrow {t_0} ^{+}}
               {L_{(\hat{f}_n(t)^{-1})}}_{*} \;
                 \Big(
                    \tfrac{1}{\mid \dot{\hat{f}} _n (t) \mid}
                    \cdot
                    \dot{\hat{f}} _n (t)
                 \Big)
    \\
    & = \; \lim_{n \rightarrow \infty }
           \lim_{t \rightarrow {s_n(t_0)} ^{-}}
             {L_{(f_n(t))^{-1}}}_{*} \; {\dot{f}_n (t)} \;
           = \; \lim_{n \rightarrow \infty }
                \lim_{t \rightarrow {s_n(t_0)} ^{+}}
                  {L_{(f_n(t))^{-1}}}_{*} \; {\dot{f}_n (t)}
    \\
    & = \; \lim_{n \rightarrow \infty}
           \lim_{t \rightarrow {s_n(t_0)}} \;
             {L_{(f_n(t))^{-1}}}_{*} \; {\dot{f}_n (t)}
    \\
    & = \; \lim_{n \rightarrow \infty}
           \lim_{t \rightarrow {s_n(t_0)}} \; {F_n} (t)
    \: .
  \end{align*}

  \noindent
  Also

  \medskip
  \noindent
  $
    F(s(t_0))        \\
    = \; \displaystyle{
                       \lim_{n \rightarrow \infty }
                       \lim_{t \rightarrow {t_0} }
                     }
      \;
      \omega \;
         \big(
           \;
           \tfrac{1}{\mid \dot{\hat{f}} _n (t) \mid } \cdot \dot{\hat{f}} _n (t)
           \;
         \big)
  $
  \begin{align*}
    = \; (-1) \cdot
      \big( \;
         \text{
               the unit curvature of the 2-dimensional
               horizontal oriented tangent plane,
              }
  \end{align*}
        \begin{flushright}
           $
             \tilde{H}^{n_0}
                 _{_e (_1 \tilde{c})_{t_0} (1) \cdot \bar{f}(t_0)}
             =
             \tilde{H}^{n_0}
                 _{ _{\bar{f}(t_0)} (_1 \tilde{c})_{t_0} (1) }
             \: ,
           $
           \hspace*{1cm}
           \\
           which projects to the tangent plane of $S$ at
           $(_1 \tilde{c}) _{t_0} (1) $  .
           \big)
        \end{flushright}

  Note paths $ _1 \tilde{c} _t  $ on $S$ gives us
  $$
    \lim_{t \rightarrow t_0 , t \in \cup ^{\infty}_{n=1} D_n}
       {_1 \tilde{c} _t} (1) \; = \; {_1 \tilde{c} _{t_0}} (1)
  $$
  and

  \begin{align*}
    \lim_{t \rightarrow t_0 , t \in \cup ^{\infty}_{n=1} D_n}
       {_{\bar{f}(t)} (_1 \tilde{c}) _t} (1)
    & = \lim_{t \rightarrow t_0 , t \in \cup ^{\infty}_{n=1} D_n}
            {_e (_1 \tilde{c}) _t} (1) \cdot \bar{f}(t)
    \\
    & = {_e (_1 \tilde{c}) _{t_0}} (1) \cdot \bar{f}(t_0)
    \\
    & = {_{\bar{f}(t_0)} (_1 \tilde{c}) _{t_0}} (1) \: .
  \end{align*}

  Then we get \\
  $
    \displaystyle{
                    \lim_{t \rightarrow t_0 , \;
                    t \in {\cup ^{\infty}_{n=1}} D_n}
                  }
           {F} (s(t))
  $
  \begin{align*}
    = \displaystyle{
                    \lim_{t \rightarrow t_0 , \;
                    t \in {\cup ^{\infty}_{n=1}} D_n}
                   }
      \; (-1) \cdot
      \big( \;
         &\text{the unit curvature of the 2-dimensional horizontal}
         \\
         &\text{oriented tangent plane, } \quad
             \tilde{H}^{n_1(t)}_{_{\bar{f}(t)} (_1 \tilde{c}) _{t} (1)}
             \quad ,
         \\
         &\text{for some }  n_1(t) \in \mathbb{N} \text{ depending on } t,
         \\
         &\text{whose projection is the tangent plane of the}
         \\
         &\text{surface } S \text{ at } _1 \tilde{c} _t (1)
         \:
       \big)
  \end{align*}

  \begin{align*}
    = \; (-1) \cdot
      \big( \;
         &\text{
                 the unit curvature of the 2-dimensional horizontal
                 oriented tangent
               }
         \\
         &\text{plane, } \quad
             \tilde{H}^{n_0}_{_{\bar{f}(t_0)} (_1 \tilde{c}) _{t_0} (1)}
           \quad ,
         \\
         &\text{whose projection is the tangent plane of the surface }
           S \text{ at } _1 \tilde{c} _{t_0} (1)
         \:
       \big)
  \end{align*}

  \noindent
  $ = {F(s(t_0))}  \: .  $

  So, we get
  $$
    F :
     \{ s(t) \mid t \in \cup ^{\infty}_{n=1} D_n - \{ 0 \} \}
     \rightarrow \text{ the unit sphere in } \mathfrak{k}
  $$
  is a continuous function. Since
  $ \{ s(t) \mid t \in \cup ^{\infty}_{n=1} D_n \} $
  is a dense subset of \\
  $[0,$ the area of $ S] $, we can extend $ F$ on
  $[0,$ the area of $ S] $ continuously. Call it also $F $.
  Consider the $C^1$ curve
  $$ \alpha : [0 , \text{ the area of } S] \rightarrow K $$
  satisfying

  \begin{align*}
  & \alpha (0) = e    \qquad  \text{and} \\
  & {L_{(\alpha(t)^{-1})}}_{*} \dot{\alpha} (t) \; = \; {F}(t)
  \: .
  \end{align*}

  Note the function
    $$
      {f_n} \; : \;
      [0, \text{ the area of the } n- \text{step pleated surface } ]
      \rightarrow K
    $$
  can be regarded as the piecewise integral curve of
    $$
      \dot{f}_n (t)
      = {L_{(f_n (t))}}_{*}
           ({L_{(f_n (t))^{-1}}}_{*} {\dot{f}_n (t)})
      = {L_{f_n (t)}}_{*} {F_n (t)} \: ,
    $$
 or equivalently the piecewise solution of the ODE
     $$
       {L_{(\alpha _n (t))^{-1}}}_{*} \dot{\alpha}_n (t)
       = {F_n (t)}  \: .
     $$
  Then
     $$
       F(s(t_0))
        =  \displaystyle{
                          \lim_{n \rightarrow \infty}
                          \lim_{ t \rightarrow s_n (t_0) }
                        }
           {L_{( f_n (t))^{-1}}}_{*} { \dot{f}_n (t)}
        =  \displaystyle{
                          \lim_{n \rightarrow \infty}
                          \lim_{ t \rightarrow s_n (t_0) }
                        }
           {F_n (t)}
     $$
  implies that
  $$ \alpha (s (t_0))
    = \displaystyle{
                     \lim_{n \rightarrow \infty}
                     \lim_{ t \rightarrow s_n (t_0) }
                   }
      {\alpha _n (t)}
    = \displaystyle{
                     \lim_{n \rightarrow \infty}
                     \lim_{ t \rightarrow s_n (t_0) }
                     }
      {f_n (t)}
    = \displaystyle{\lim_{n \rightarrow \infty} }
      {f_n (s_n (t_0))}
    = {f (s(t_0))} \: .
  $$

   Since $ \{ s(t) \mid t \in \cup ^{\infty}_{n=1} D_n \} $ is a
   dense subset of $[0,$ the area of $ S] $ and $F$ is continuous
   on [0, the area of $S$], we get
   $$
     f = \alpha
     \text { is a } C^1 \text{ curve on } [0, \text{ the area of } S]
     \: .
   $$

  Also, we obtain

  \begin{align*}
  \text{the length of the curve } f
  & = \text{the length of the curve } \alpha \\
  & = \int^{\text{the area of } S}_{0} \mid \dot{\alpha} (t) \mid dt \\
  & = \int^{\text{the area of } S}_{0} \mid {F(t)} \mid dt \\
  &= \text{the area of } S \: ,
  \end{align*}

  \noindent
  which proves Theorem \ref{thm}.

  \subsubsection{\textbf{Remarks on Factorization Lemma}}
  \label{factorization}

  `Factorization Lemma', introduced by Lichnerowicz,
  \emph{
        Theorie Globale des Connexions et des Groupes
        d'Holonomie
        },
  [3, vol 1, p.284],
  can give us another sequence of piecewise
  smooth loops
  $\mu_m : [0,1] \rightarrow \mathbb{H}^n, \; m=1,2, \cdots,$ with
  $\mu_m(0)= \pi(e)$
  such that it converges to $\partial S.$ And a similar way to make
  the sequence of curves $f_n : [0,1] \rightarrow K, \; n=1,2, \cdots,$
  can give us a sequence of curves
  $g_n : [0,1] \rightarrow K, \; n=1,2, \cdots,$ with $g_n(0)=e$
  such that
  $g_n(1)$ is the ending point of the horizontal lifting of
  $\mu_n$ at $e$ and that the length of $g_n $ is the area of the
  \emph{pleated surface}, the union of totally geodesic triangles
  obtained in the construction of $g_n$. Since the sequence of
  the areas converges to the area of $S,$
  Prop \ref{converge-lift} will say that $g_n(1)$ will converge
  to $_e \tilde{\gamma}(1)$ and that the distance from $e$ to
  $_e \tilde{\gamma}(1)$ is less that equal to the area of $S$.
  But the sequence $\{g_n\}$ may not converge to some curve from
  $e$ to $_e \tilde{\gamma}(1).$

\pagebreak

\begin{appendices}

\section{About Triangles}\label{sec:triangles}

  \medskip

  For each $n = 0, 1, 2, \cdots ,$ all triangles inside
  $2^n \cdot 3$-gon will consist of two kinds of triangles,
  \emph{interior} ones and \emph{exterior} ones.

  \medskip

  \subsection{The definition of interior triangles and the definition
  of their starting points and ending points}

  Consider a regular triangle whose vertices lie on the
  boundary of the given disk $ D^2 $ and one of whose vertices is the
  base point of the disk. Call the triangle $ T_0 $.  And the base
  point will be called its starting and ending point.

  Now let's define triangles $ T_{a_0 a_1 \cdots a_n} $ inductively
  as follows :

  \medskip

  Case 1)   $ n = 1 $ :

    The given orientation at the center of $ D^2 $ and  the base
    point, or equivalently the starting and ending point of $ T_0 $,
    will give the order $b_0$ of sides of $ T_0 $, where $b_0 =
    1.2.3$,  in the counter-clockwise or clockwise order. For the
    barycentric subdivision of $ T_0 $, thinking of the triangle
    with the base point as its vertex and with one side lying on the
    first side of $ T_0 $ as the first triangle will give the order
    of triangles in the counter-clockwise or clockwise order. The
    $i$-th triangle will be called $ T_{a_0 a_1} $, where $a_0 =0$
    and $a_1 =i$, for $i = 1, 2, \cdots , 6 .$

    \medskip

    \begin{figure}[h]
      \centering
        {
          \includegraphics[width=1.5in]{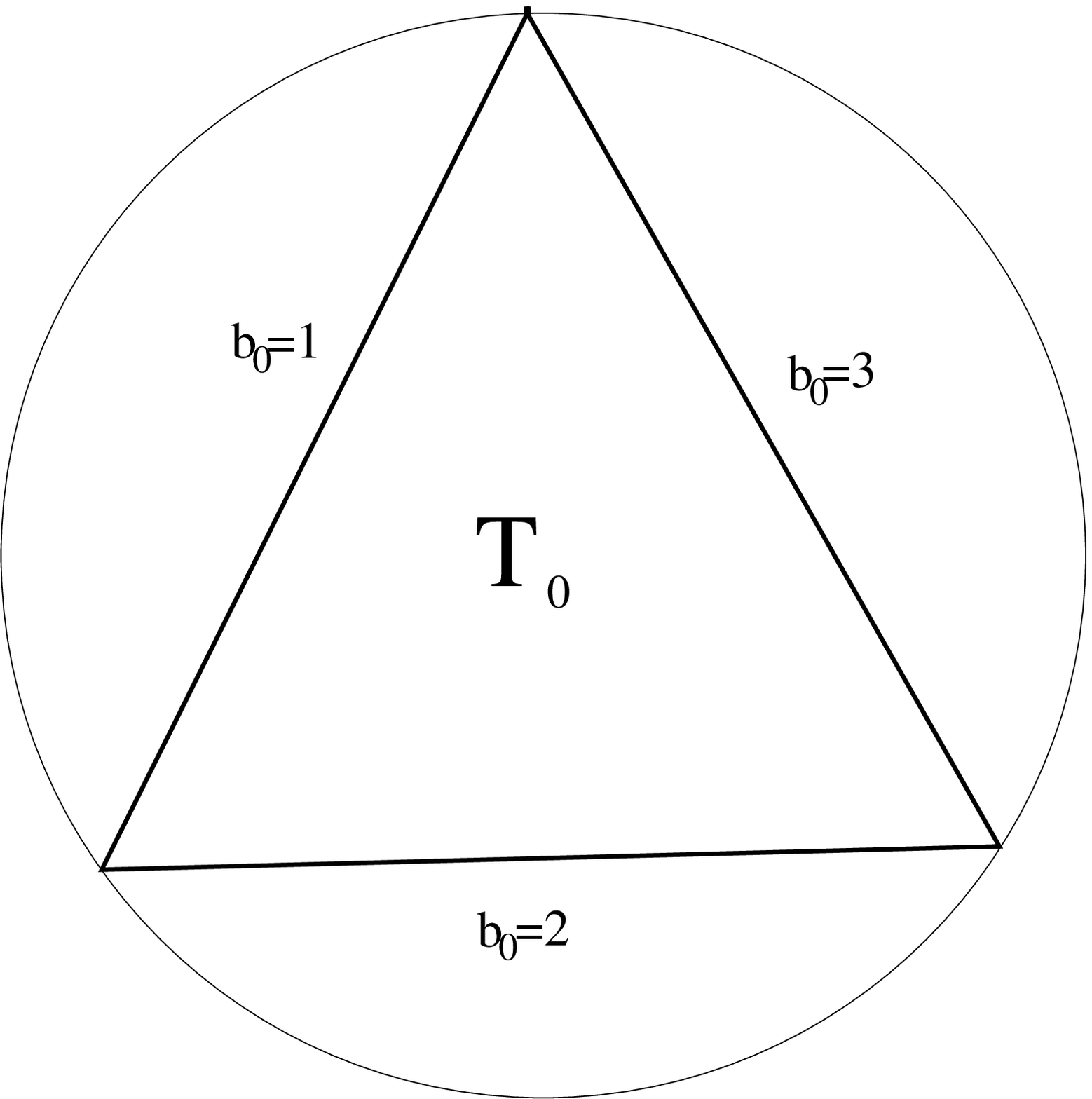}
          \hspace{1cm} \text{ or } \hspace{1cm}
          \includegraphics[width=1.5in]{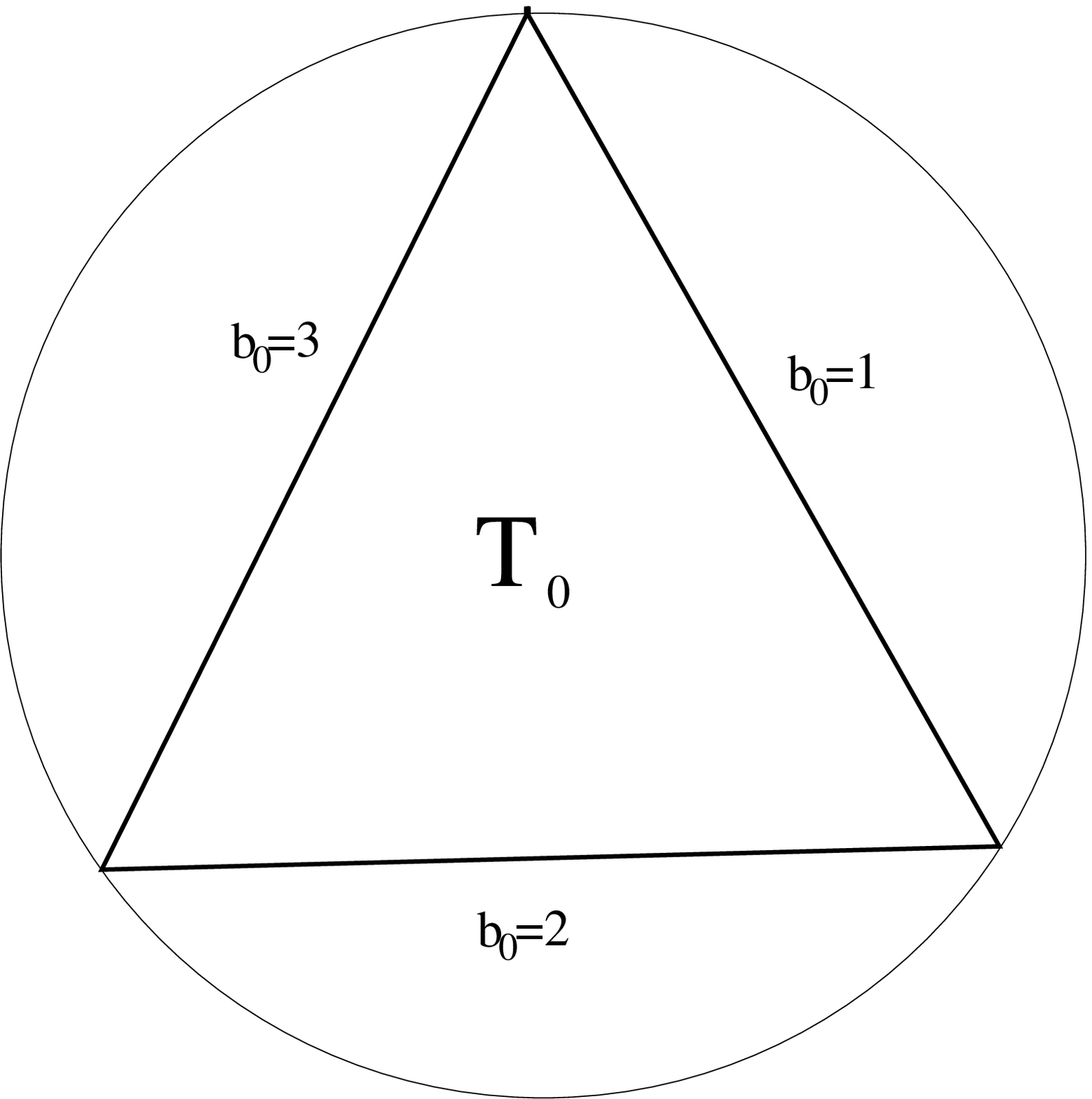}
        }
      \vspace{1cm}
      \centering
        {
          \includegraphics[width=1.5in]{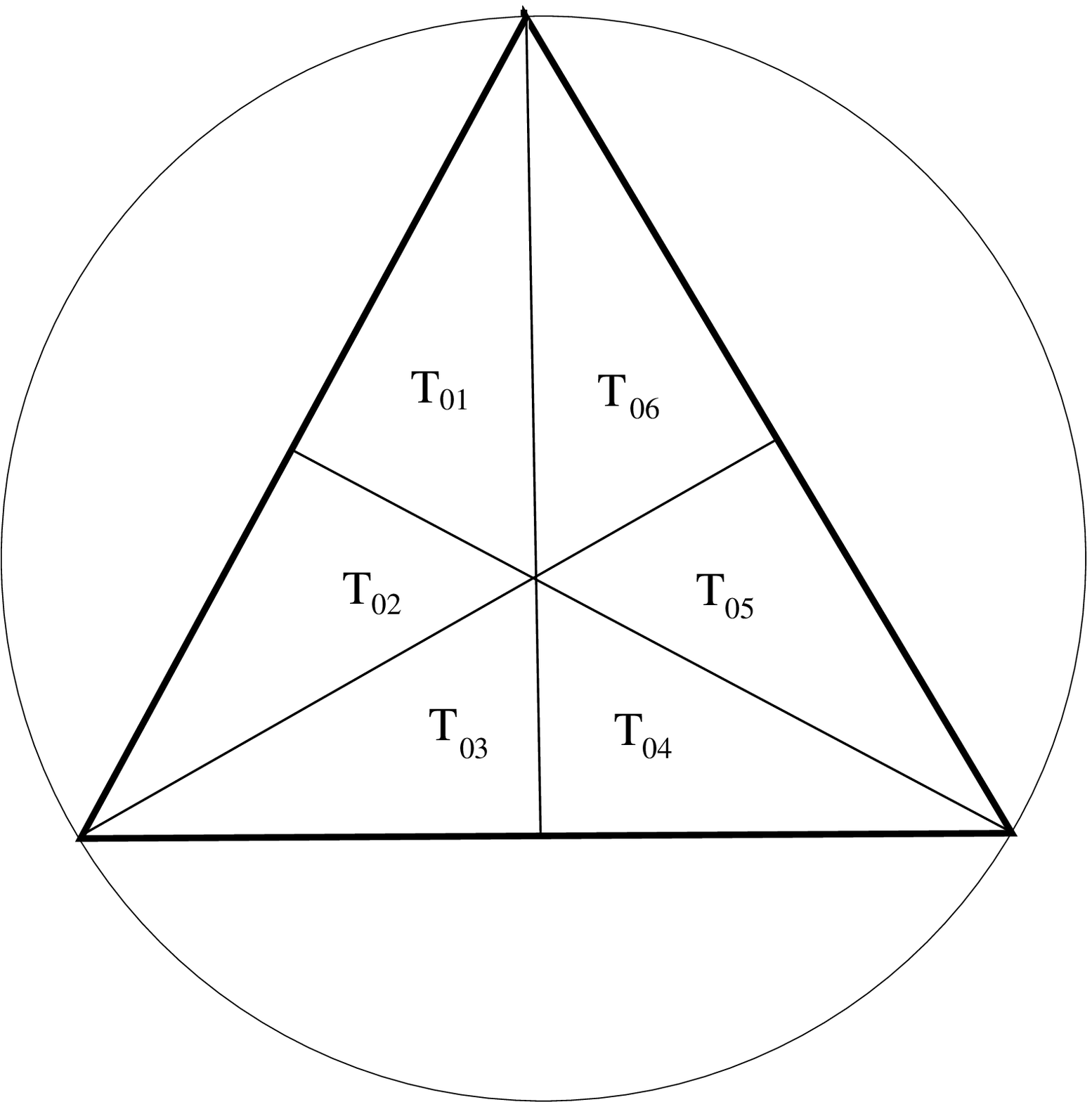}
          \hspace{3cm}
          \includegraphics[width=1.5in]{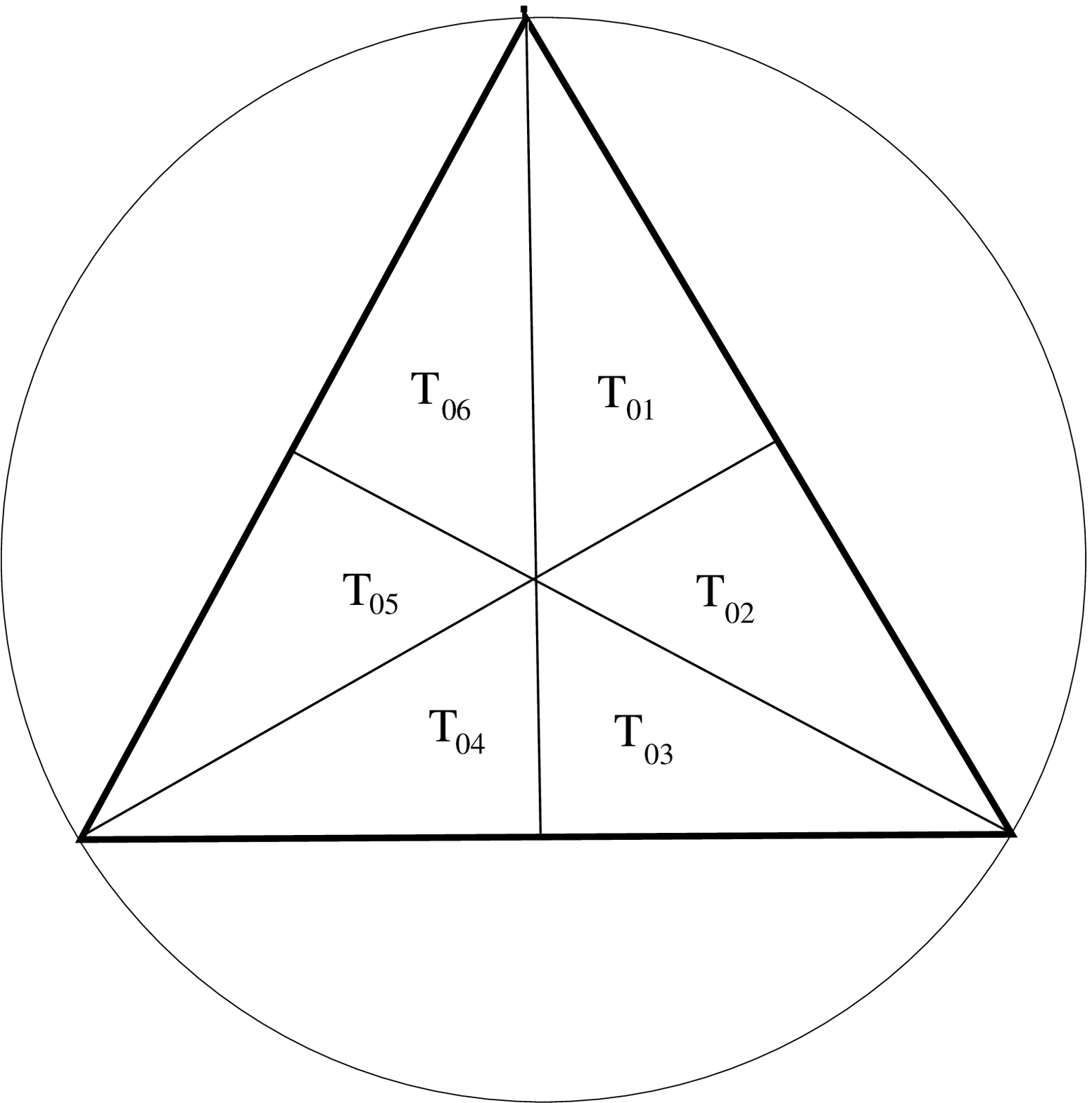}
        }
    \end{figure}

    For $ T_{01} $ , the base point, or equivalently the starting
    point of $ T_0 $, will be called the starting point of $T_{01} $
    and the barycenter of $ T_0 $ will be called the ending point of
    $ T_{01} $.

    For $ T_{0i} $ , where $ i=2,3,4,5 $ , the barycenter of $ T_0 $
    will be called the starting and ending point of $T_{0i} $ for
    $i = 2,3,4,5 $ .

    For $ T_{06} $ , the barycenter of $ T_0 $ will be called the
    starting point of $T_{06} $ and the base point , or equivalently
    the ending point of $ T_0 $ , will be called the ending point of
    $ T_{06}$

  \medskip

  Case 2)   $ n \geq 2 $  :

    Let $ L_{n-1} := T_{a_0 a_1 \cdots a_{n-2} a_{n-1}} $ be given,
    where $a_0=0$ and $a_1, \cdots , a_{n-1} \in \{1,2,3,4,5,6\}$.
    Let $L_{n-2} := T_{a_0 \cdots a_{n-2}}$ and assume the
    following properties:

    - $M_j := T_{a_0 a_1 \cdots a_{n-2} \, j}, \, j \in \{1,2,3,4,5,6 \},$
    consists one of six triangles obtained by the barycentric
    subdivision of $L_{n-2}$,

    - $L_{n-1}$ is also one of those, in other words,
    $$
       L_{n-1} = T_{a_0 a_1 \cdots a_{n-2} a_{n-1}}
       = T_{a_0 a_1 \cdots a_{n-2} \, j_0} = M_{j_0}
    $$
    $
       \text{for some } j_0 \in \{1,2,3,4,5,6 \}.
    $

    \begin{figure}[h]
    \centering
      {
        \includegraphics[width=1.5in]{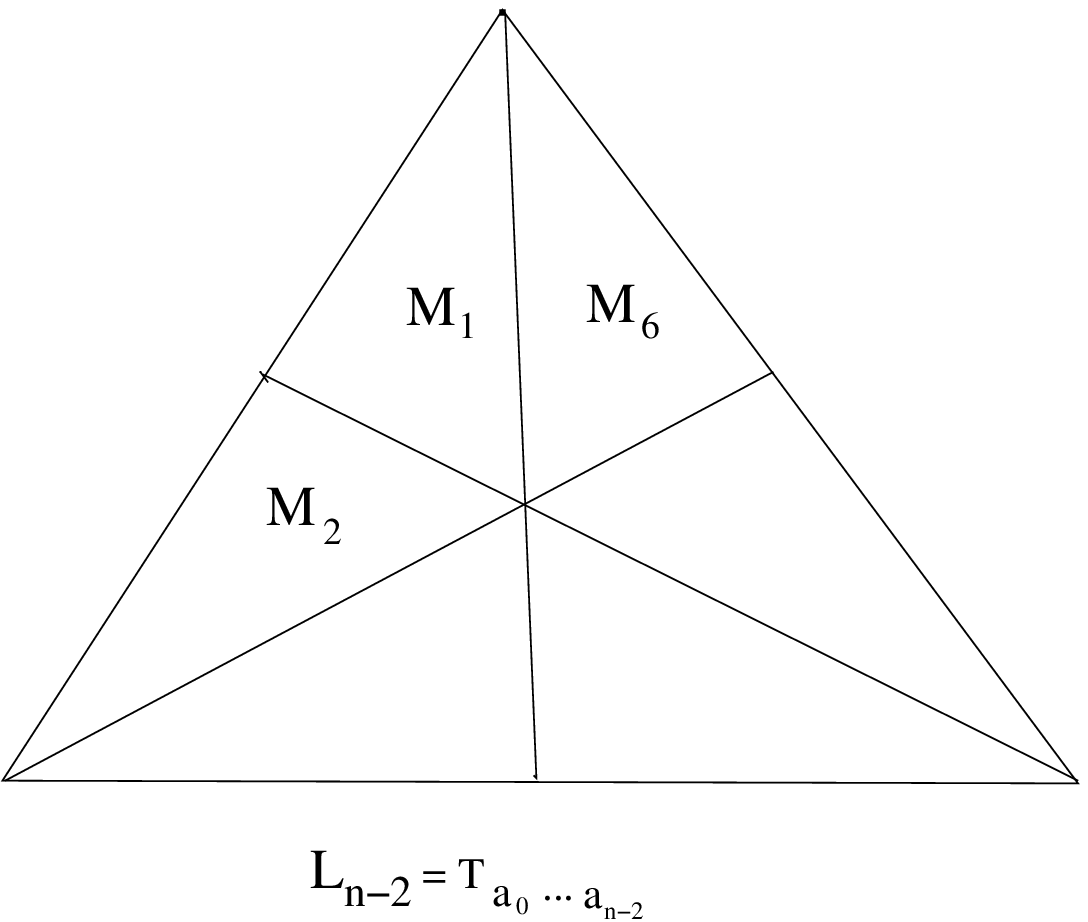}
        \hspace{2cm}
        \includegraphics[width=1.5in]{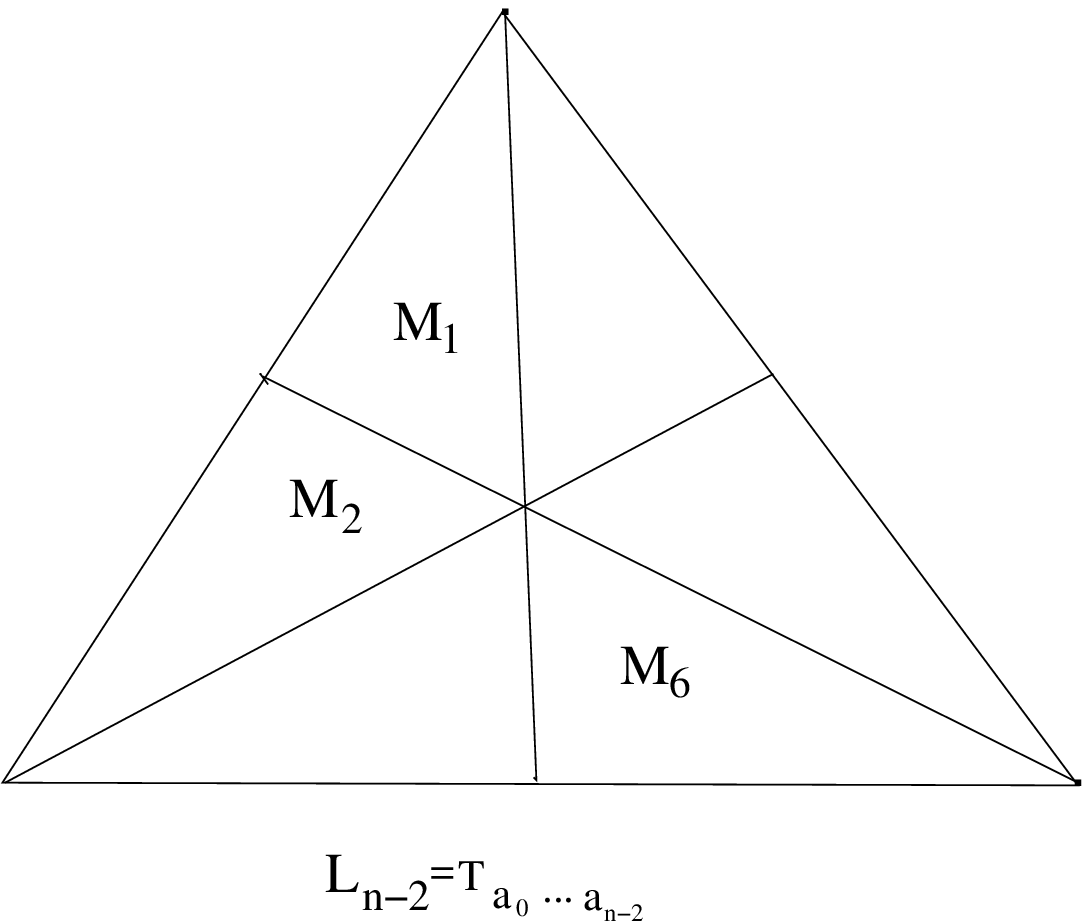}
      }
    \\
    \vspace{1cm}
    \centering
      {
        \includegraphics[width=1.5in]{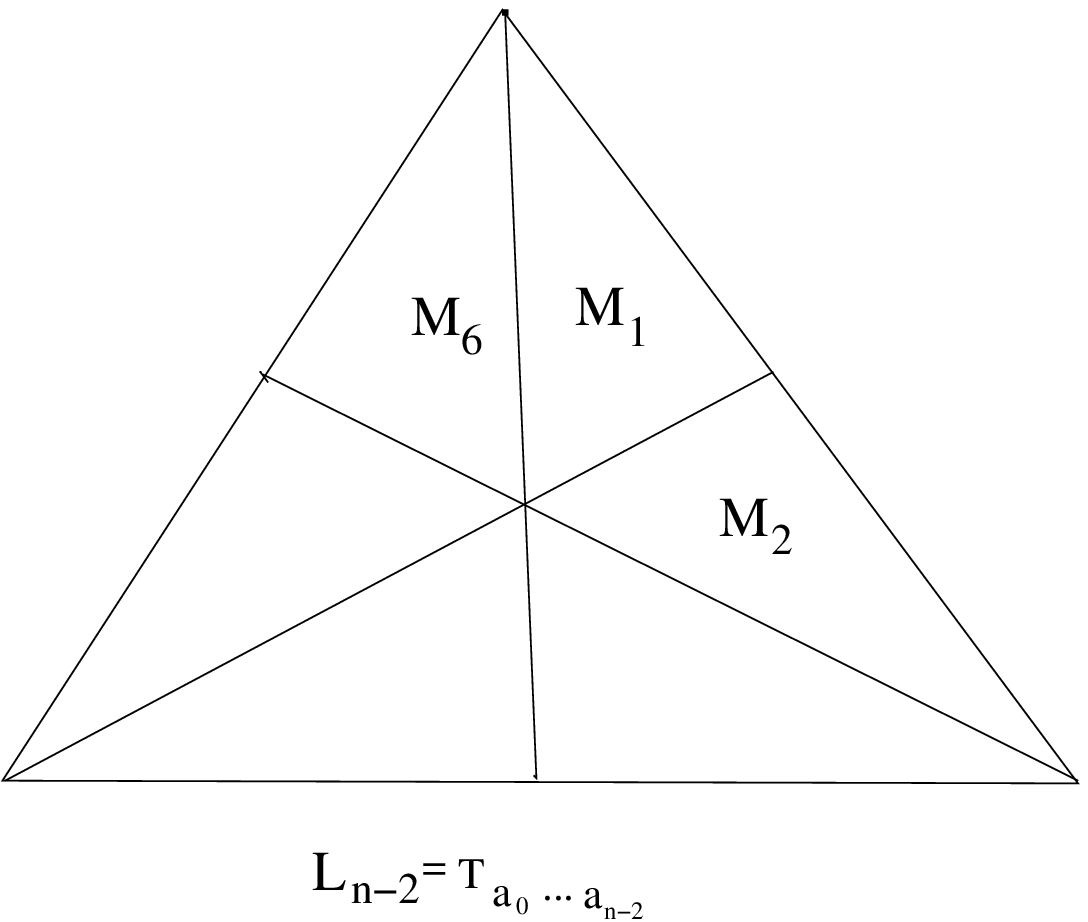}
        \hspace{2cm}
        \includegraphics[width=1.5in]{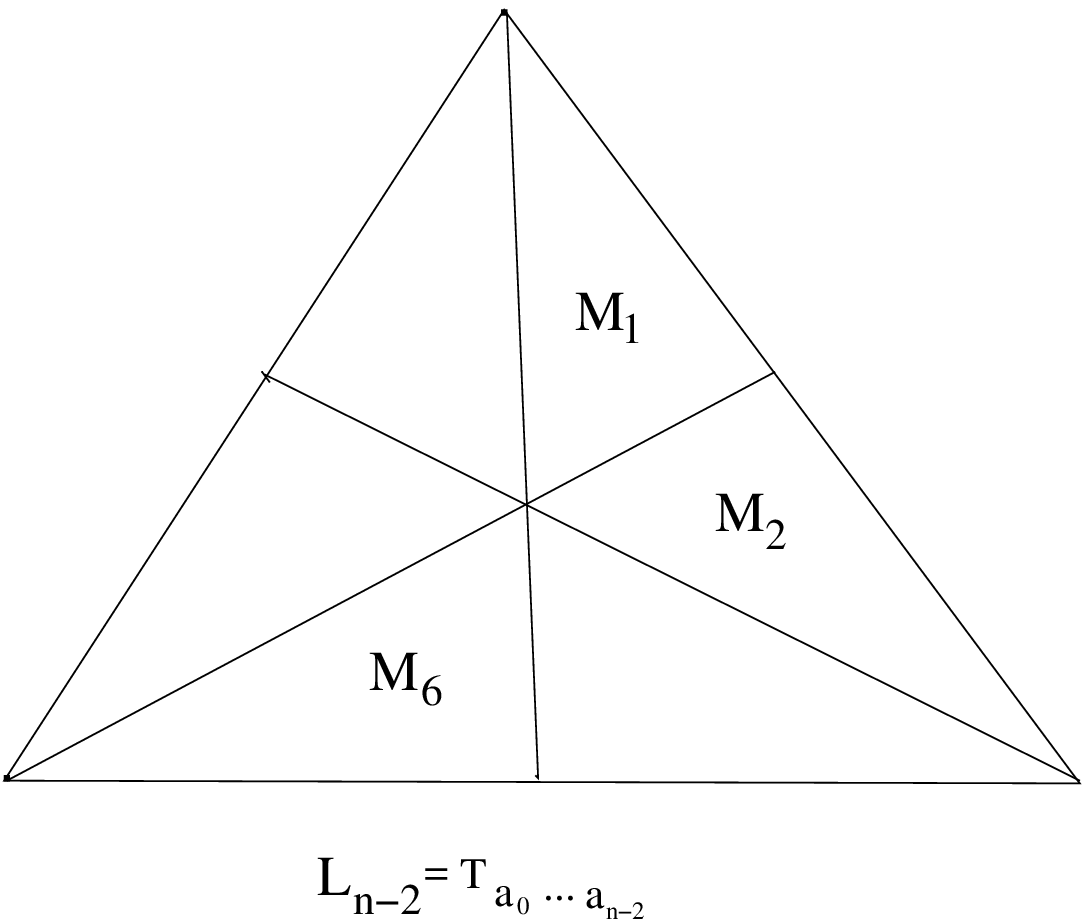}}
    \end{figure}

    - common vertex of $L_{n-2}$ and $ M_{1} $ is the
    starting point of each of them,

    - the barycenter of $ L_{n-2} $ is the starting point of $ M_{i} $
    for $i=2,3,4,5,6,$
    and the ending point of $ M_{i}$ for $i=1,2,3,4,5$,

    - the common vertex of $ L_{n-2} $ and $ M_{6} $ is the ending point of
    each of them,

    - if the starting and the ending point of $ L_{n-2} $
    are same, then they are the common vertex of $ L_{n-2} $ and
    $ M_{6} $,

    - if the starting and the ending point of
    $ L_{n-2} $ are different, then $ M_{1} $ and $ M_{6} $ are
    mutually opposite ones inside $ L_{n-2} $,

    - one side of $L_{n-2}$, which contains a side of $ M_1 $,
    is divided into two line segments, each of which is  one side of
    $ M_{i} $ for $i=1,2$, respectively.

    Notice all the above assumptions hold for n=2.

    Note that the line segment connecting the barycenter and the
    starting point of $ L_{n-2} $ is one side of $ M_{1} $ from the
    assumption that the common vertex of $ L_{n-2} $ and $ M_{1} $ is
    the starting point of each of them.

    Under the above assumptions, we have two choices such that the
    order of $ M_{1} $ and $ M_{2} $ is either the counter-clockwise
    order or the clockwise order with respect to the barycenter of
    $ L_{n-2} $ and the line segment connecting the barycenter and the
    starting point of $ L_{n-2} $.

    \bigskip

    Subcase 2-1 ) $ a_{n-1} = 1 $, that is,
                  $ L_{n-1} = M_{1} = T_{a_0 \cdots a_{n-2} 1 }$ :

      \begin{figure}[h]
        \centering{\includegraphics[width=2in]{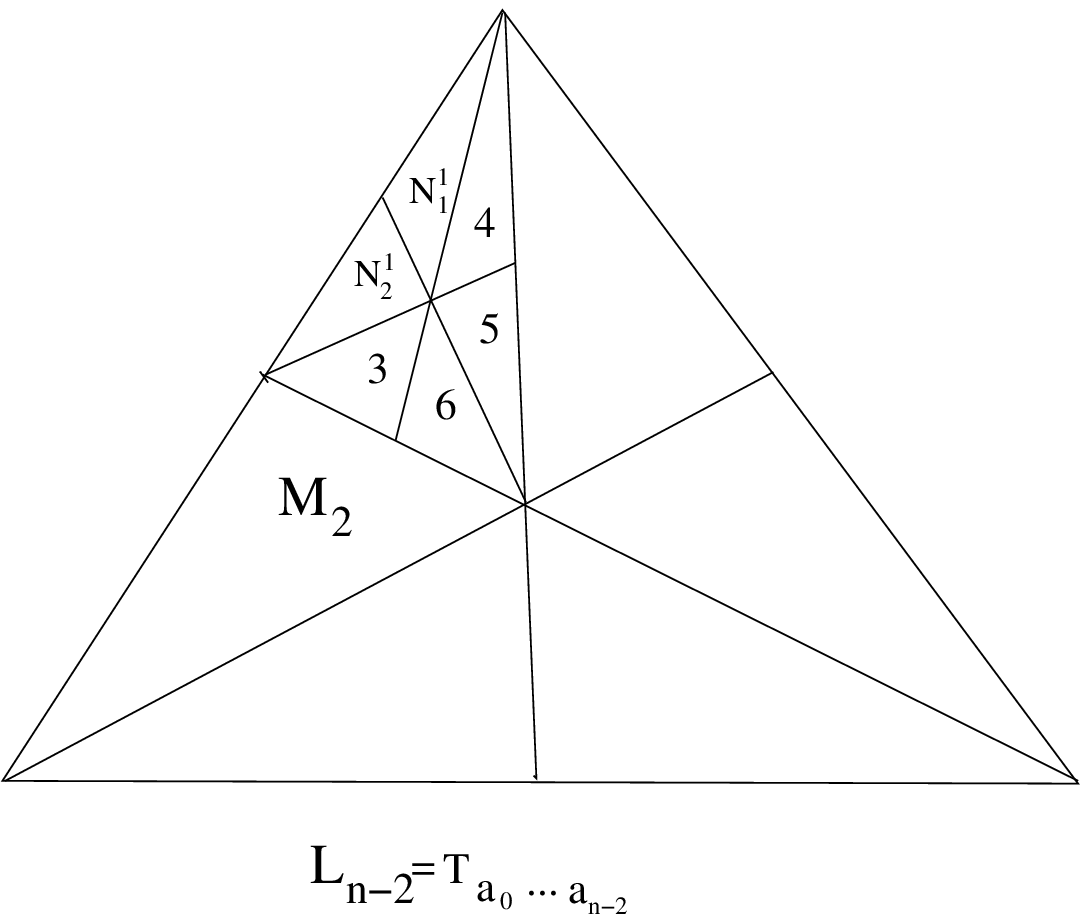}}
      \end{figure}

      Assume the order of $ L_{n-1} = M_{1} $ and $ M_{2} = T_{a_0
      \cdots a_{n-2} 2}$ is the counter-clockwise order with respect
      to the barycenter of $ L_{n-2} = T_{a_0 \cdots a_{n-2} }$ and
      the line segment connecting the barycenter and the starting
      point of $ L_{n-2} $. Out of six triangles obtained from the
      barycentric subdivision of $ L_{n-1} = M_{1} $, choose the
      triangle with a part of one side of $ L_{n-2} $ as its side and
      with the starting point and the barycenter of $ L_{n-1} = M_{1}
      $ as its vertices, and  call it  $ T_{a_0 \cdots a_{n-2} 1 1} $
      and let $N^{1}_{1} := T_{a_0 \cdots a_{n-2} 1 1}$. At the
      barycenter of $ L_{n-1} = M_1 $, consider the counter-clockwise
      order of the 6 triangles from the $ N^{1}_{1} $. The 5
      triangles from the next one of $ N^{1}_{1} $ will be called
      $$
       T_{a_0 \cdots a_{n-2} 1 2}, \: T_{a_0 \cdots a_{n-2} 1 3}, \:
       T_{a_0 \cdots a_{n-2} 1 6}, \: T_{a_0 \cdots a_{n-2} 1 5}, \:
       T_{a_0 \cdots a_{n-2} 1 4}
      $$
      in order. Let $N^{1}_{i} := T_{a_0 \cdots a_{n-2} 1 i}$ for
      $i=2,3,4,5,6.$

      \medskip

      If the order of $ L_{n-1} = M_{1} $ and $ M_{2} $ is the
      clockwise order, then the order will be given from the symmetry
      by the line connecting the barycenter and the starting point of
      $ L_{n-1} = M_{1} $:

      \begin{figure}[h]
        \centering{\includegraphics[width=2in]{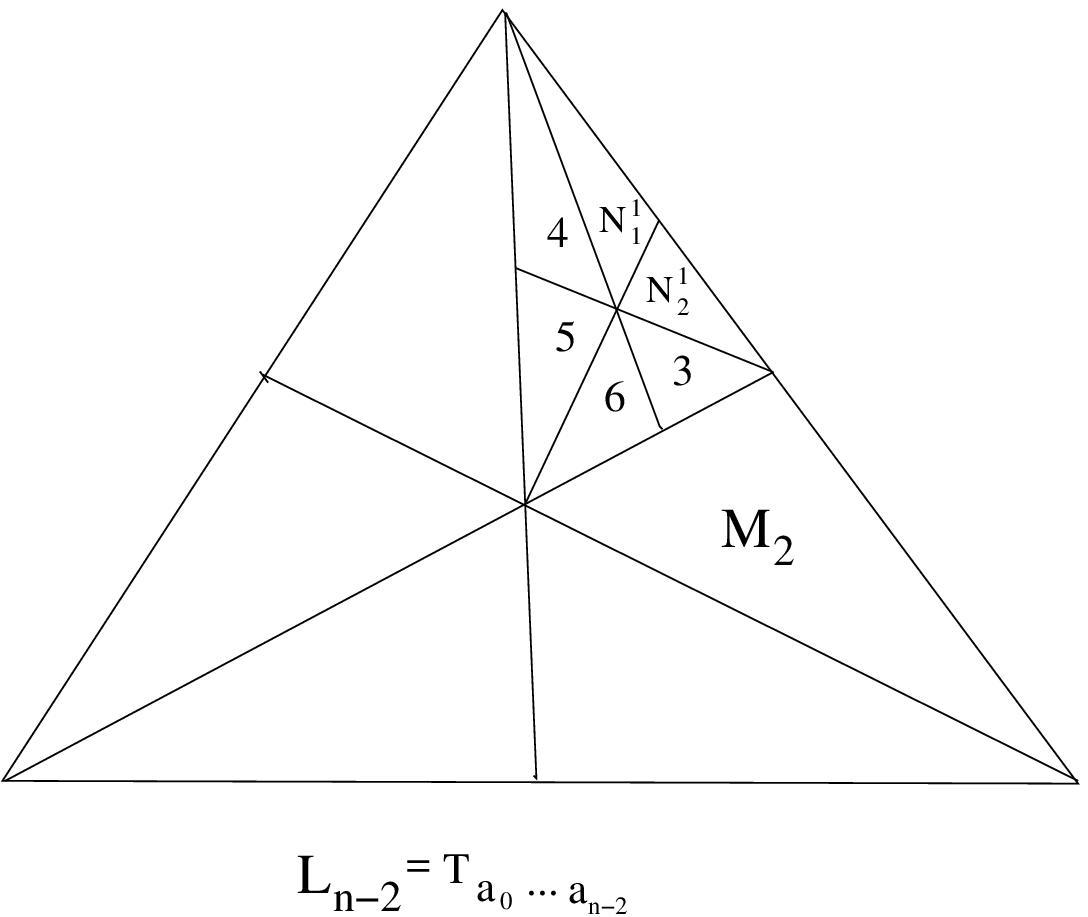}}
      \end{figure}

      Recall the assumptions for $ L_{n-1} = M_{1}, $ lying between
      the phrase `Case 2)   $ n \geq 2 $' and the one
      `Subcase 2-1 ) $ a_{n-1} = 1 $, $\cdots$,' and let
      $ L_n := N^1_j, \; j=1, \cdots , 6 .$

      Note the common vertex of $ L_{n-1} = M_{1} $ and $ N^{1}_{1} $
      is the starting point of $L_{n-1} = M_{1}$ from the definition
      of $ N^{1}_{1} \; .$ Now, call the vertex the starting point of
      $N^{1}_{1} \; .$ And call the barycenter of $L_{n-1} = M_{1}$
      the starting point of $N^{1}_{i}$ for $i=2,3,4,5,6.$ Also, call
      the barycenter the ending point of $N^{1}_{i}$ for
      $i=1,2,3,4,5.$

      Note the common vertex of $ L_{n-1} = M_{1} $ and $ N^{1}_{6}$
      is the barycenter of $ L_{n-2} $, so the ending point of
      $L_{n-1} = M_{1} $ from the assumption for $M_{1}.$ Call the vertex
      the ending point of $N^{1}_{6} $.

      Note that the starting and the ending point of
      $ L_{n-1} = M_{1} $ are different and the positions of
      $ N^{1}_{1} $ and $ N^{1}_{6} $ are mutually opposite inside
      $ L_{n-1} = M_{1} $.

      And the side of $ L_{n-1} = M_{1} $, which contains a side of
      $N^{1}_{1}$, is divided into two line segments, each of which is
      one side of $ N^{1}_{i} $ for $i=1,2$, respectively.

     \bigskip

    Subcase 2-2 ) $ a_{n-1} = 6 $, that is,
                  $ L_{n-1} = M_{6} = T_{a_0 \cdots a_{n-2} 6 }$ :

      \begin{figure}[h]
        \centering
          {
            \includegraphics[width=2in]{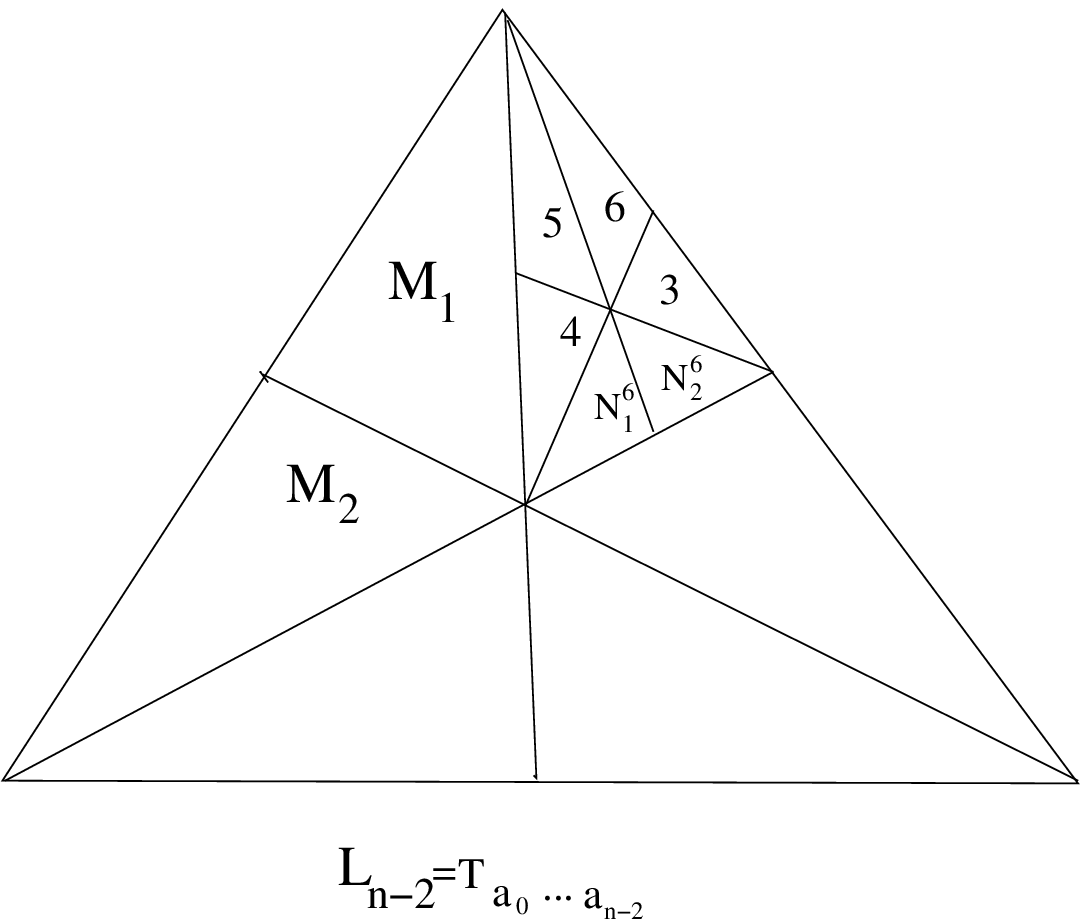}
            \hspace{0.5cm} \text{ or } \hspace{0.5cm}
            \includegraphics[width=2in]{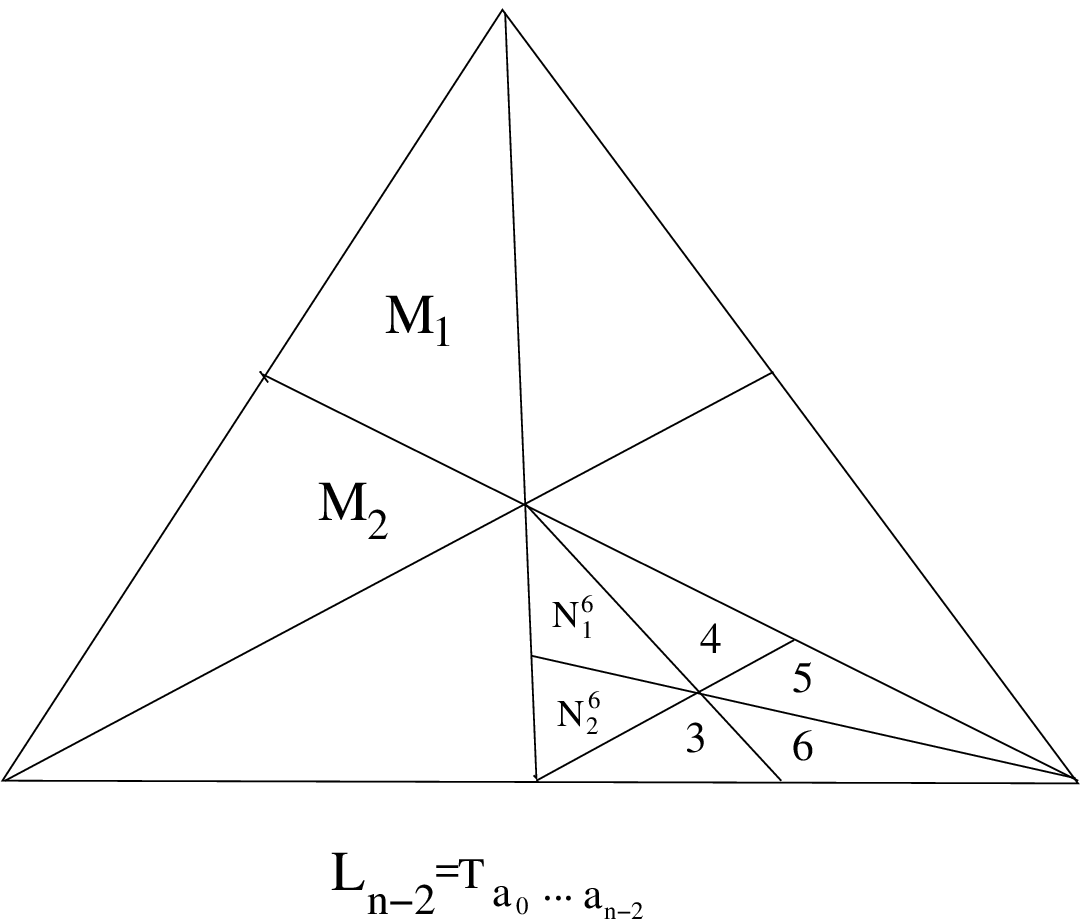}
          }
      \end{figure}

      Assume the order of $ M_{1} $ and $ M_{2} $ is the
      counter-clockwise order with respect to the barycenter of $
      L_{n-2} $ and the line segment connecting the barycenter and
      the starting point of $ L_{n-2} $. From the assumptions, lying
      between the phrase `Case 2)   $ n \geq 2 $' and the one
      `Subcase 2-1 ) $ a_{n-1} = 1 $, $\cdots$,' the vertex of
      $ L_{n-1} = M_{6} $, which is also the barycenter of
      $ L_{n-2} $, is the starting point of $ L_{n-1} = M_{6} $. The
      counter-clockwise angle of $ L_{n-1} = M_{6} $ at its starting
      point determines its initial side and the terminal side.  Out
      of six triangles obtained from the barycentric subdivision of
      $ L_{n-1} = M_{6} $, choose the triangle with a part of the
      initial side of $ L_{n-1} = M_{6} $ as its side and with the
      starting point and the barycenter of $ L_{n-1} = M_{6} $ as its
      vertices, and  call it $ T_{a_0 \cdots a_{n-2} 6 1} $ and let
      $N^{6}_{1} := T_{a_0 \cdots a_{n-2} 6 1}$. At the barycenter of
      $ L_{n-1} = M_{6} $, consider the counter-clockwise order of
      the 6 triangles from the $ N^{6}_{1} $. The 5 triangles from
      the next one of $ N^{6}_{1} = T_{a_0 \cdots a_{n-2} 6 1} $ will
      be called
      $$
         T_{a_0 \cdots a_{n-2} 6 2} , T_{a_0 \cdots a_{n-2} 6 3} ,
         T_{a_0 \cdots a_{n-2} 6 6} , T_{a_0 \cdots a_{n-2} 6 5} ,
         T_{a_0 \cdots a_{n-2} 6 4}
      $$
      in order. Let $N^{6}_{i} := T_{a_0 \cdots a_{n-2} 6 i}$ for
      $i=2,3,4,5,6.$

      If the order of $ M_{1} $ and $ M_{2} $ is the clockwise order,
      then the order will be given from the symmetry by the line
      connecting the barycenter and the starting point of
      $ L_{n-1} = M_{6} $:

      \begin{figure}[h]
       \centering
         {
           \includegraphics[width=1.9in]{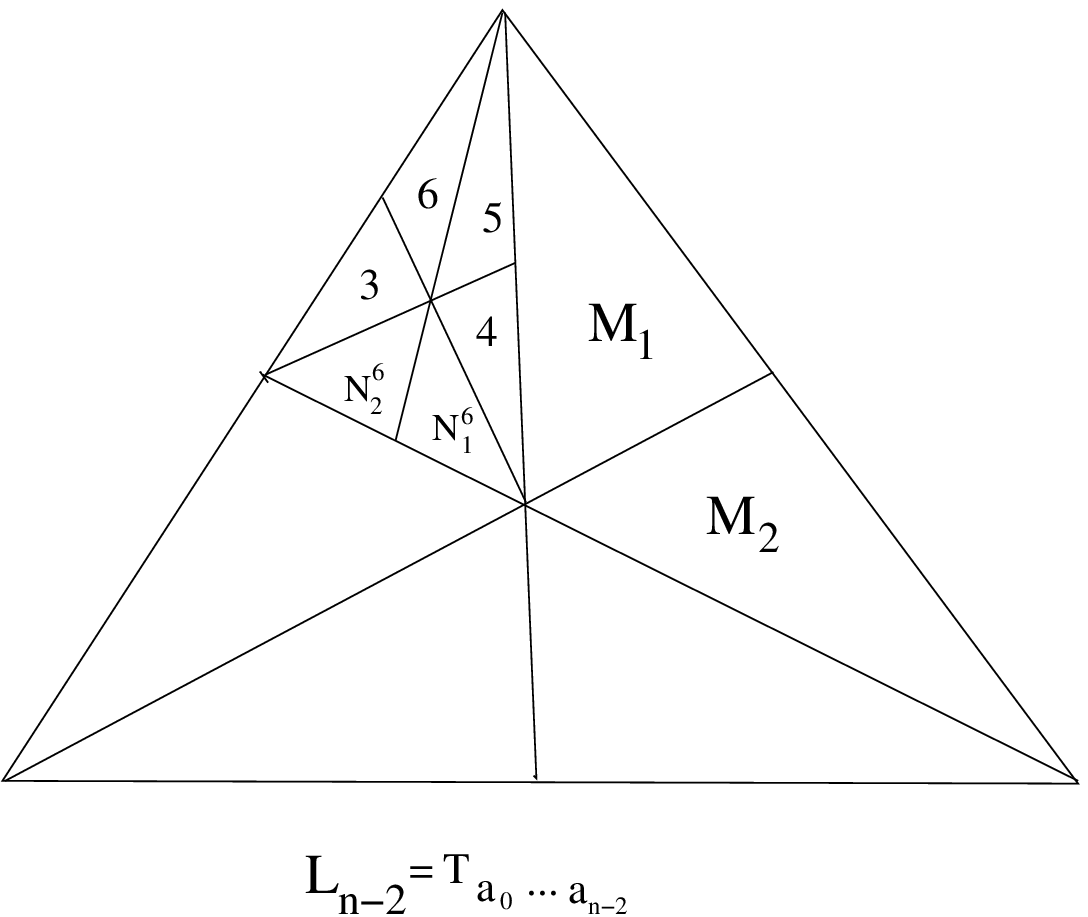}
           \hspace{0.5cm} \text{or} \hspace{0.5cm}
           \includegraphics[width=1.9in]{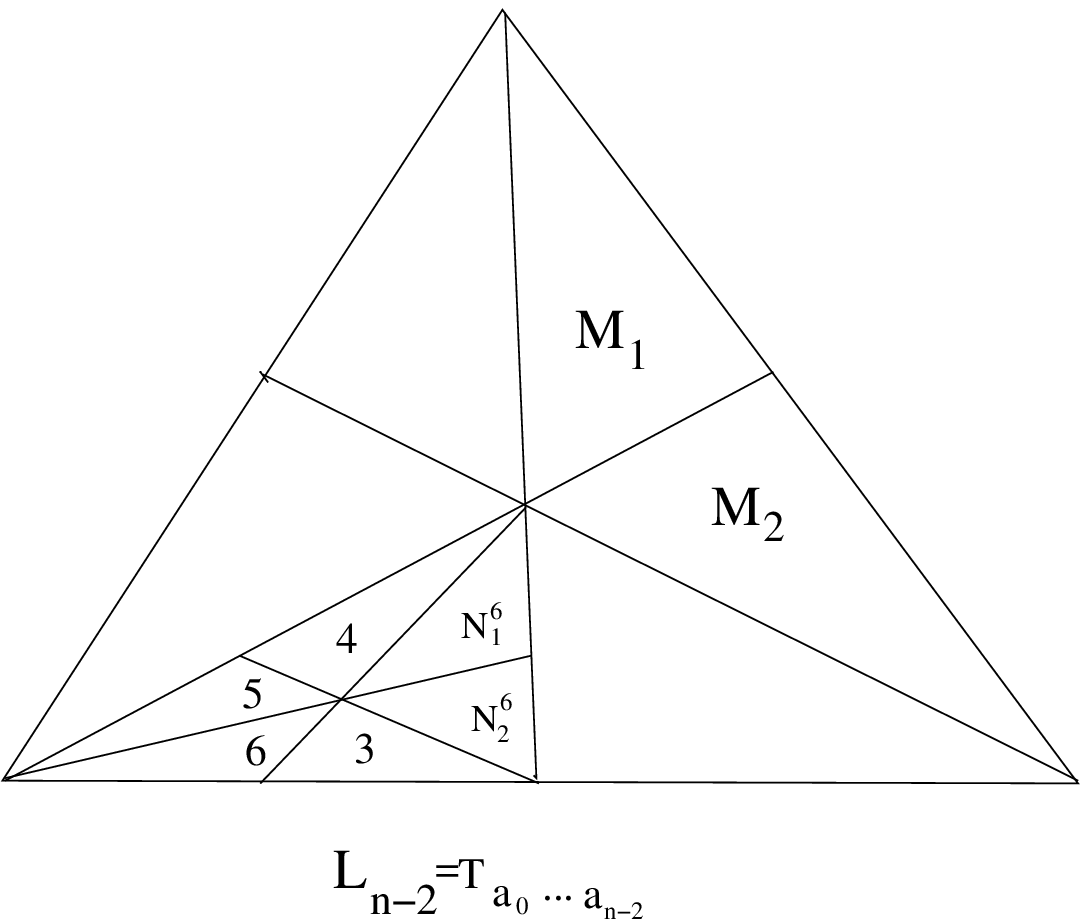}
         }
      \end{figure}

      Recall the assumptions for $ L_{n-1} = M_{6}, $ lying between
      the phrase `Case 2)   $ n \geq 2 $' and the one
      `Subcase 2-1 ) $ a_{n-1} = 1 $, $\cdots$,' and let
      $ L_n := N^6_j, \; j=1, \cdots , 6 .$

      Note the common vertex of $ L_{n-1} = M_{6} $ and $ N^{6}_{1} $
      is the starting point of $L_{n-1} = M_{6}$ from the definition
      of $ N^{6}_{1} $. Now, call the vertex the starting point of
      $ N^{6}_{1}$. And call the barycenter of $ L_{n-1} = M_{6} $
      the starting point of $N^{6}_{i} $ for $i=2,3,4,5,6.$ Also,
      call the barycenter the ending point of $N^{6}_{i} $ for
      $i=1,2,3,4,5.$

      To consider the common vertex of $ L_{n-1} = M_{6} $ and
      $ N^{6}_{6} $, we have the following two possibilities :

      The starting point and the ending point of $ L_{n-2} $ are same
      or different.

      But in any possibilities, the common vertex of
      $L_{n-1} = M_{6}$ and $ N^{6}_{6} $ is also the common vertex
      of $L_{n-2}$ and $ L_{n-1} = M_{6} $, so the ending point of
      $ L_{n-1} = M_{6} $ from the assumption for $ M_{6} $. Call the
      vertex the ending point of $N^{6}_{6} $.

      Note that the starting and the ending point of $L_{n-1} = M_{6}$
      are different and the positions of $ N^{6}_{1} $ and $N^{6}_{6}$
      are mutually opposite inside $ L_{n-1} = M_{6} $.

      Notice  the side of  $ L_{n-1} = M_{6} $, which contains a side
      of $ N^{6}_{1} $, is divided into two line segments, each of
      which is one side of $ N^{6}_{i} $ for $i=1,2$, respectively.

    \pagebreak

    Subcase 2-3) $a_{n-1} \in  \{ 2,3 \}$ or
    $(a_{n-1} \in  \{ 4,5 \} $ and $ a_{n-2} \in \{ 0,2,3,4,5 \} ),$
    \\
    that is,
    $$
      L_{n-1} = M_{i} = T_{a_0 \cdots a_{n-2} \, i } \text{ for } \:
      i=2,3
    $$
    or
    $$
      L_{n-1} = M_{i} = T_{a_0 \cdots a_{n-2} i } \text{ for } \:
      i=4,5 \text{ and }  a_{n-2} \in \{ 0,2,3,4,5 \} :
    $$

      \begin{figure}[h]
        \centering
          {
            \includegraphics[width=1.9in]{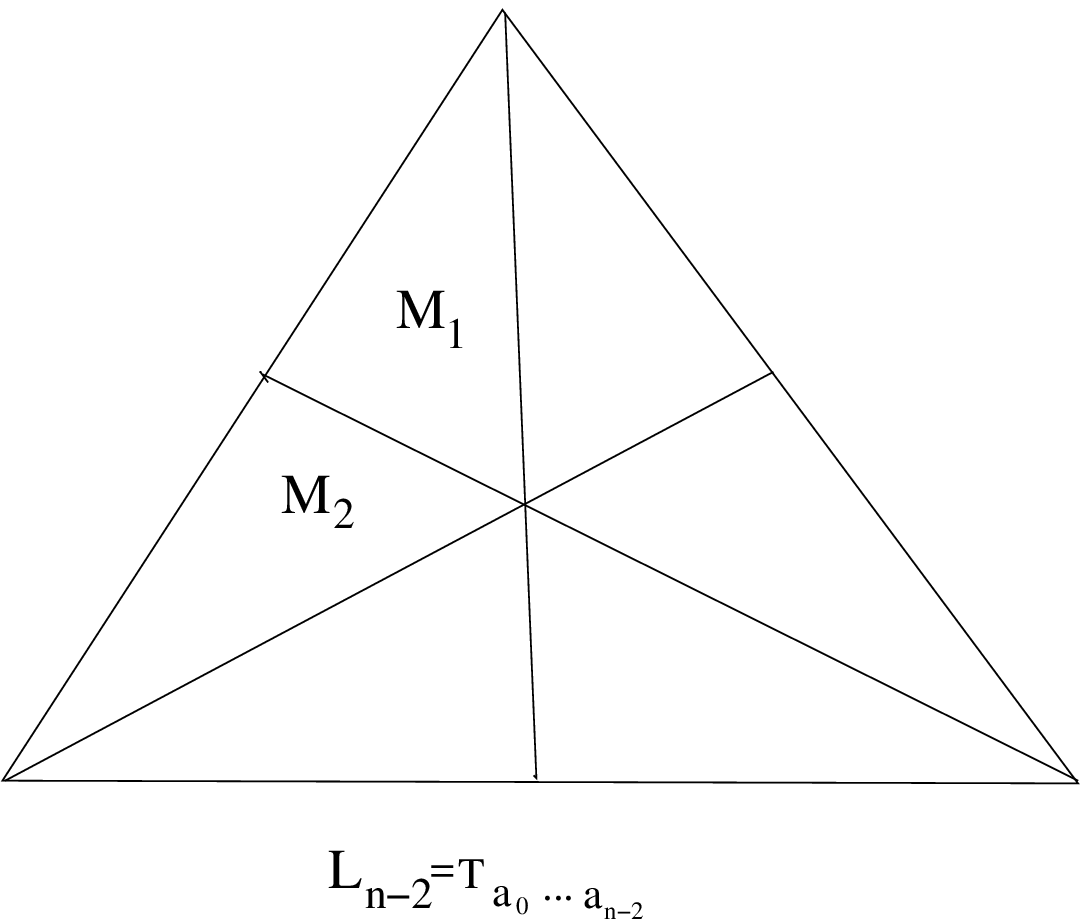}
            \hspace{1cm}
            \includegraphics[width=.5cm]{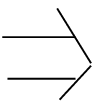}
            \hspace{1cm}
            \includegraphics[width=1.9in]{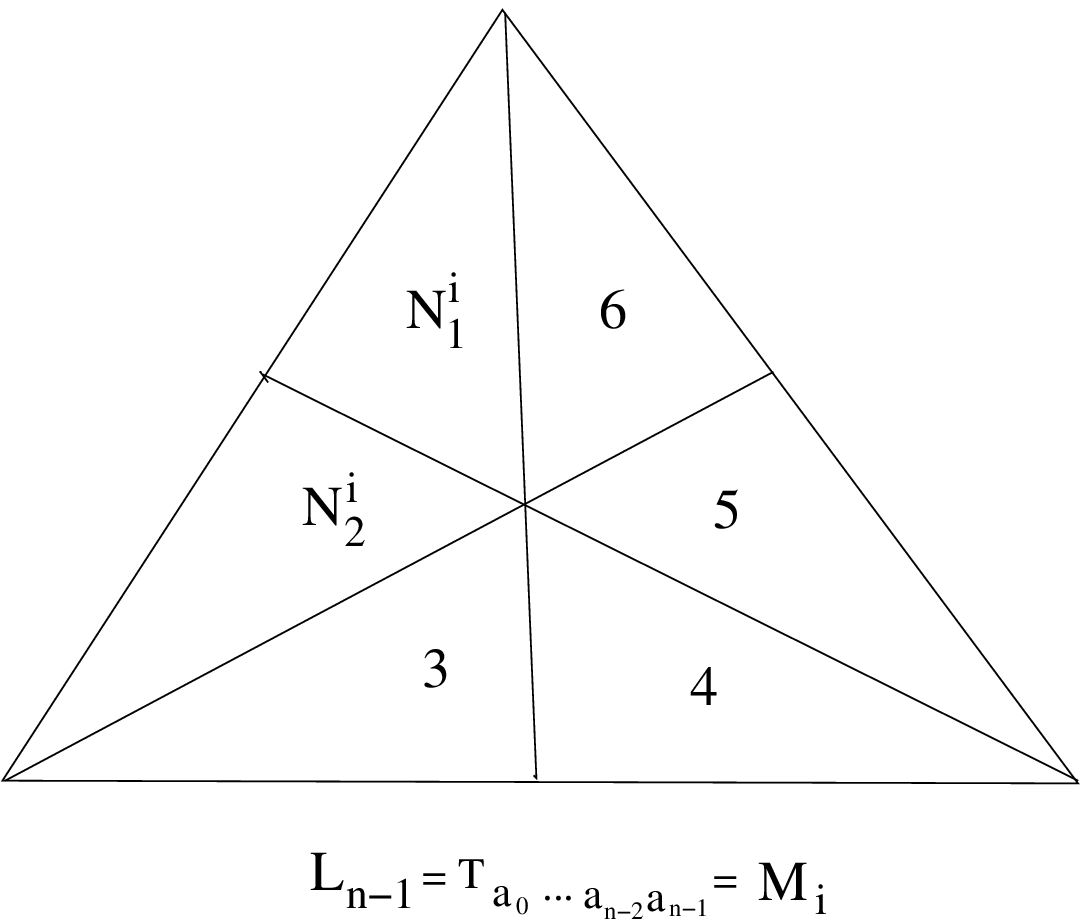}
          }
      \end{figure}

      Let $ i = a_{n-1}$.

      Assume the order of $ M_{1} $ and $ M_{2} $ is the
      counter-clockwise order with respect to the barycenter of
      $ L_{n-2} $ and the line segment connecting the barycenter and
      the starting point of $ L_{n-2} $. From the assumptions, lying
      between the phrase `Case 2)   $ n \geq 2 $' and the one
      `Subcase 2-1 ) $ a_{n-1} = 1 $, $\cdots$,' the
      vertex of $ L_{n-1} = M_{i} $, which is also the barycenter of
      $ L_{n-2} $, is the starting point of $ L_{n-1} = M_{i} $. The
      counter-clockwise angle of $ L_{n-1} = M_{i} $ at its starting
      point determines the initial side and the terminal side.  Out
      of six triangles obtained from the barycentric subdivision of
      $ L_{n-1} = M_{i} $, choose the triangle with a part of the
      initial side of $ L_{n-1} = M_{i} $ as its side and with the
      starting point and the barycenter of $ L_{n-1} = M_{i} $ as its
      vertices, and  call it  $ T_{a_0 \cdots a_{n-2} a_{n-1} 1} $,
      in other words, $ T_{a_0 \cdots a_{n-2} i 1} $, and let
      $N^{i}_{1} := T_{a_0 \cdots a_{n-2} i 1}$. At the barycenter of
      $ L_{n-1} = M_{i} $, consider the counter-clockwise order of
      the 6 triangles from the $ N^{i}_{1} $. The 5 triangles from
      the next one of $N^{i}_{1} = T_{a_0 \cdots a_{n-2} a_{n-1} 1}$
      will be called
      $$
        T_{a_0 \cdots a_{n-2} a_{n-1} 2} ,
        T_{a_0 \cdots a_{n-2} a_{n-1} 3} ,
        T_{a_0 \cdots a_{n-2} a_{n-1} 4} ,
        T_{a_0 \cdots a_{n-2} a_{n-1} 5} ,
        T_{a_0 \cdots a_{n-2} a_{n-1} 6}
      $$
      in order. Let
      $
        N^{i}_{j} := T_{a_0 \cdots a_{n-2} \, i j}
        = T_{a_0 \cdots a_{n-2} a_{n-1} \, j}
      $
      for $j=2,3,4,5,6.$

     \bigskip

     If the order of $ M_{1} $ and $ M_{2} $ is the clockwise order,
     then the order will be given from the symmetry by the line
     connecting the barycenter and the starting point of
     $ L_{n-1} = M_{i} = T_{a_0 \cdots a_{n-2} a_{n-1}} \; :$

     \pagebreak

     \begin{figure}[h]
       \centering
         {
           \includegraphics[width=1.9in]{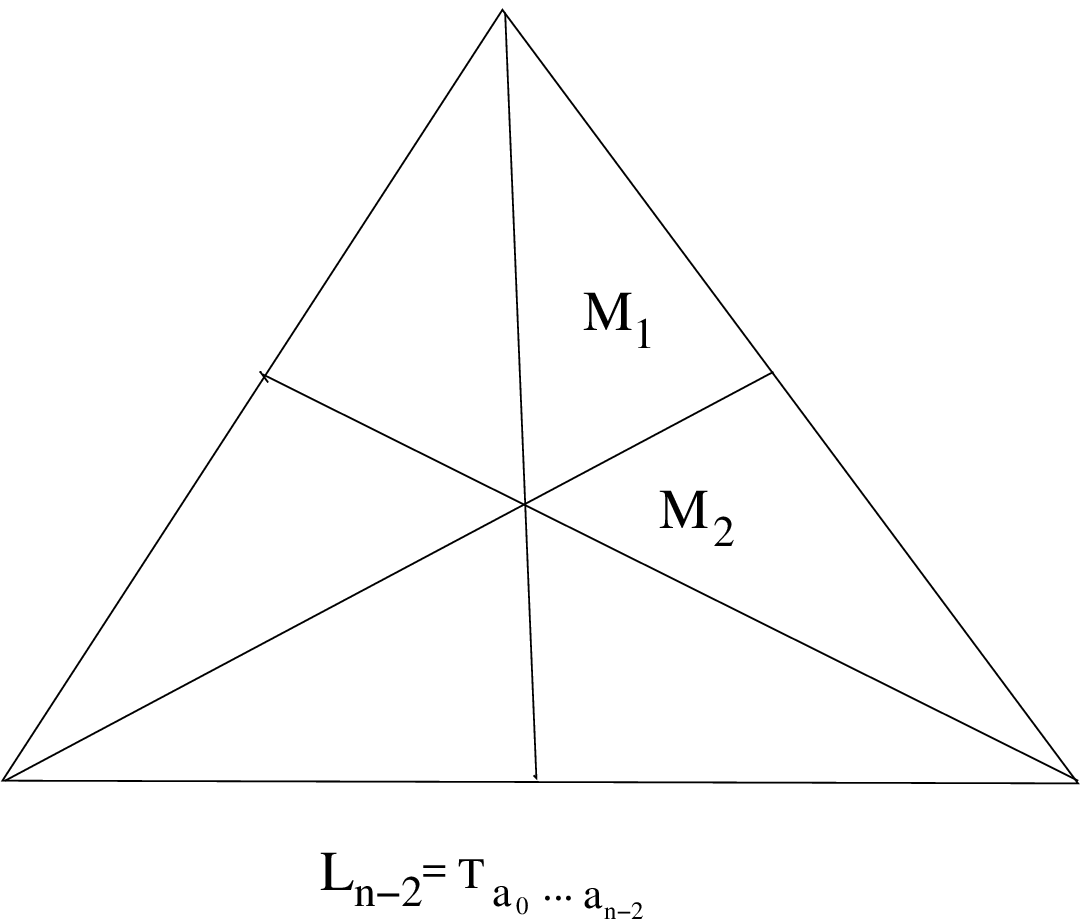}
           \hspace{1cm}
           \includegraphics[width=.5cm]{rightarrow}
           \hspace{1cm}
           \includegraphics[width=1.9in]{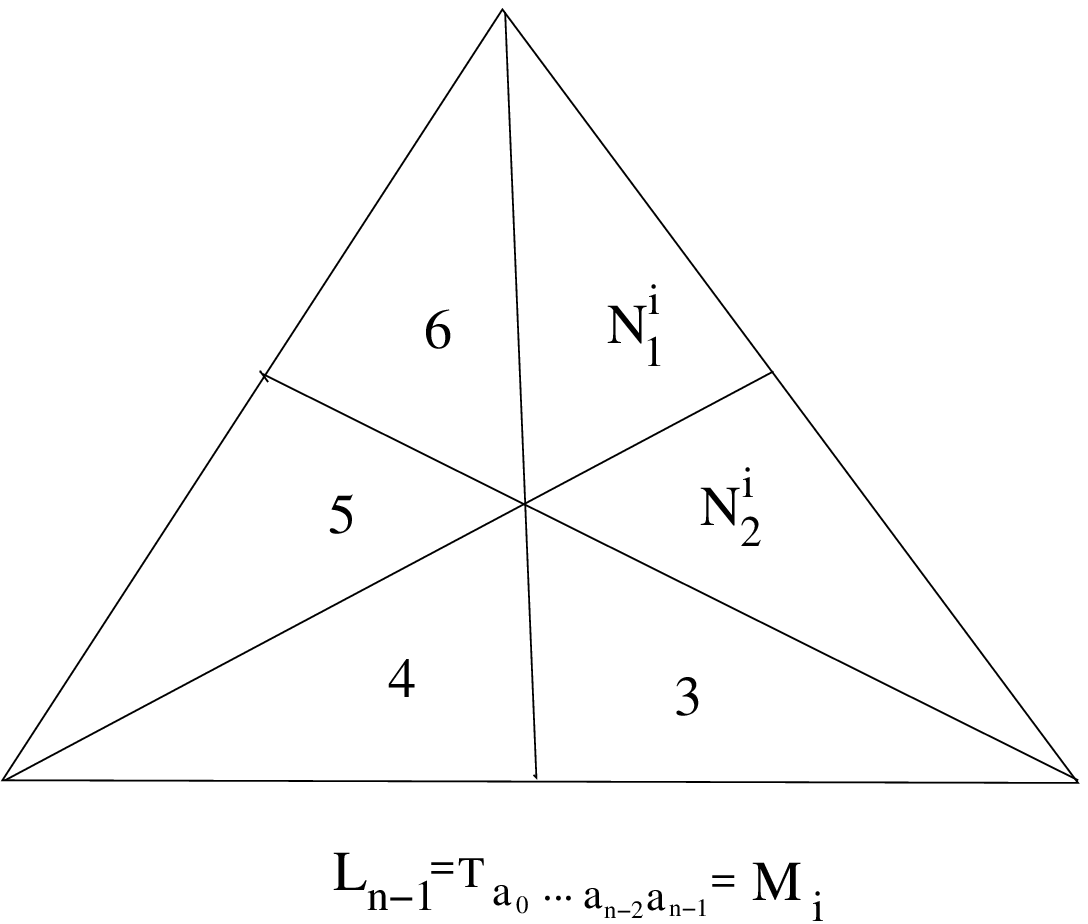}
         }
     \end{figure}

     Recall the assumptions for $ L_{n-1} = M_{i}, $ lying between
     the phrase `Case 2)   $ n \geq 2 $' and the one
     `Subcase 2-1 ) $ a_{n-1} = 1 $, $\cdots$,' and let
     $ L_n := N^i_j, \; j=1, \cdots , 6 .$

     Note the common vertex of $ L_{n-1} = M_{i} $ and $ N^{i}_{1} $
     is the starting point of $ L_{n-1} = M_{i} $ from the
     definition of $ N^{i}_{1} $. Now, call the vertex the starting
     point of  $ N^{i}_{1} \; .$ And call the barycenter of
     $ L_{n-1} = M_{i} $ the starting point of $N^{i}_{j} $ for
     $j=2,3,4,5,6.$ Also, call the barycenter the ending point of
     $N^{i}_{j} $ for $j=1,2,3,4,5.$

     Note the common vertex of $ L_{n-1} = M_{i} $ and
     $ N^{i}_{6} $ is the starting point of $ L_{n-1} = M_{i} \; , $
     so the barycenter of $ L_{n-2} $ and the ending point of
     $ L_{n-1} = M_{i} $ from the assumption for $M_{i}$. Call the
     vertex the ending point of $ N^{i}_{6} $.

     Note that the starting and the ending point of
     $ L_{n-1} = M_{i} $ are same and they are the common vertex of $
     L_{n-1} = M_{i} $ and $ N^{i}_{6}  $.

     And the side of $ L_{n-1} = M_{i} $, which contains a side of
     $ N^{i}_{1} $, is divided into two line segments, each of which
     is one side of $ N^{i}_{j}  $ for $j=1,2$, respectively.

     \bigskip

    Subcase 2-4 ) $ a_{n-1} \in  \{ 4,5 \} $ and
                   $ a_{n-2} \in  \{ 1,6 \} $, \\
    that is,
    $$
      L_{n-1} = M_{i} = T_{a_0 \cdots a_{n-2} i } \text{ for } \:
      i=4,5 \text{ and }  a_{n-2} \in \{ 1,6 \} :
    $$

      \begin{figure}[h]
        \centering
          {
            \includegraphics[width=2in]{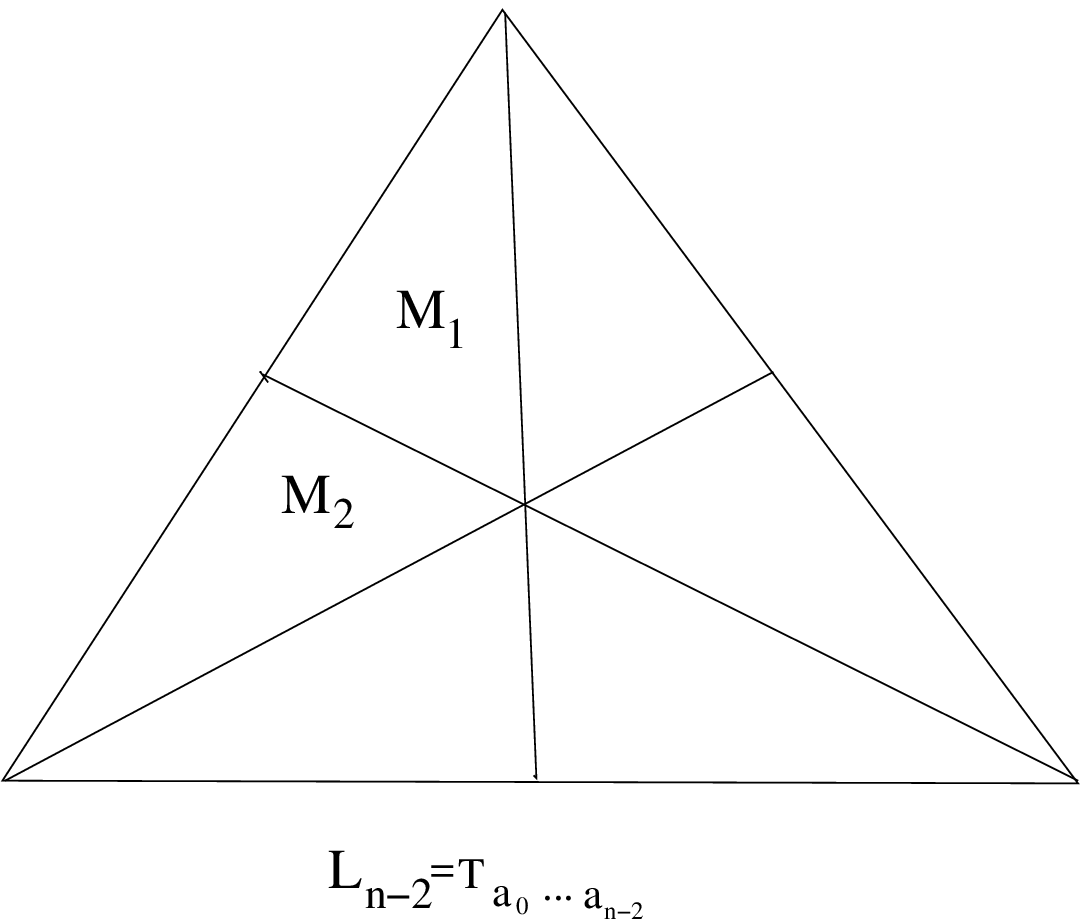}
            \hspace{.5cm}
            \includegraphics[width=.5cm]{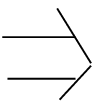}
            \hspace{.5cm}
            \includegraphics[width=2in]{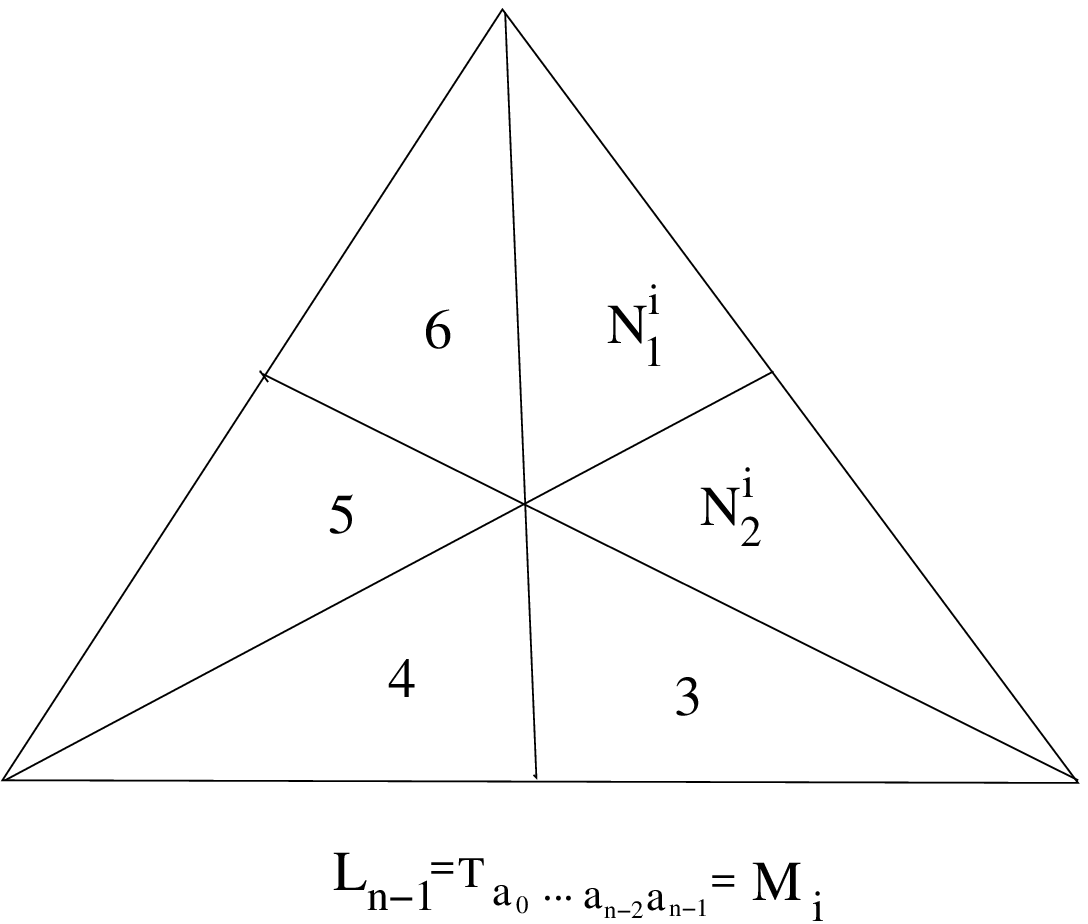}
          }
      \end{figure}

      Let $ i = a_{n-1}$.

      Assume the order of $ M_{1} $ and $ M_{2} $ is the
      counter-clockwise order with respect to the barycenter of
      $  L_{n-2} $ and the line segment connecting the barycenter and
      the starting point of $ L_{n-2} $. From the assumptions,lying
      between the phrase `Case 2)   $ n \geq 2 $' and the one
     `Subcase 2-1 ) $ a_{n-1} = 1 $, $\cdots$,' the
      vertex of $ L_{n-1} = M_{i} $, which is also the barycenter of
      $ L_{n-2} $, is the starting point of $ L_{n-1} = M_{i} $. The
      clockwise angle of $ L_{n-1} = M_{i} $ at its starting point
      determines the initial side and the terminal side.  Out of six
      triangles obtained from the barycentric subdivision of
      $ L_{n-1} = M_{i} $, choose the triangle with a part of the
      initial side of $ L_{n-1} = M_{i} $ as its side and with the
      starting point and the barycenter of $ L_{n-1} = M_{i} $ as its
      vertices, and  call it $T_{a_0 \cdots a_{n-2} a_{n-1} 1} \; ,$
      in other words, $ T_{a_0 \cdots a_{n-2} i 1} $, and let
      $N^{i}_{1} := T_{a_0 \cdots a_{n-2} i 1}.$ At the barycenter of
      $ L_{n-1} = M_{i} $, consider the clockwise order of the 6
      triangles from the $ N^{i}_{1} $. The 5 triangles from the next
      one of
      $
        N^{i}_{1} = T_{a_0 \cdots a_{n-2} i 1}
        = T_{a_0 \cdots a_{n-2} a_{n-1} 1}
      $
      will be called
      $$
        T_{a_0 \cdots a_{n-2} a_{n-1} 2} ,
        T_{a_0 \cdots a_{n-2} a_{n-1} 3} ,
        T_{a_0 \cdots a_{n-2} a_{n-1} 4} ,
        T_{a_0 \cdots a_{n-2} a_{n-1} 5} ,
        T_{a_0 \cdots a_{n-2} a_{n-1} 6}
      $$
      in order. Let
      $
        N^{i}_{j} := T_{a_0 \cdots a_{n-2} \, i j}
        = T_{a_0 \cdots a_{n-2} a_{n-1} \, j}
      $
      for $j=2,3,4,5,6.$

      \medskip

      If the order of $ M_{1} $ and $ M_{2} $ is the clockwise order,
      then the order will be given from the symmetry by the line
      connecting the barycenter and the starting point of
      $ L_{n-1} = M_{i} \; : $

      \begin{figure}[h]
        \centering
          {
            \includegraphics[width=1.7in]{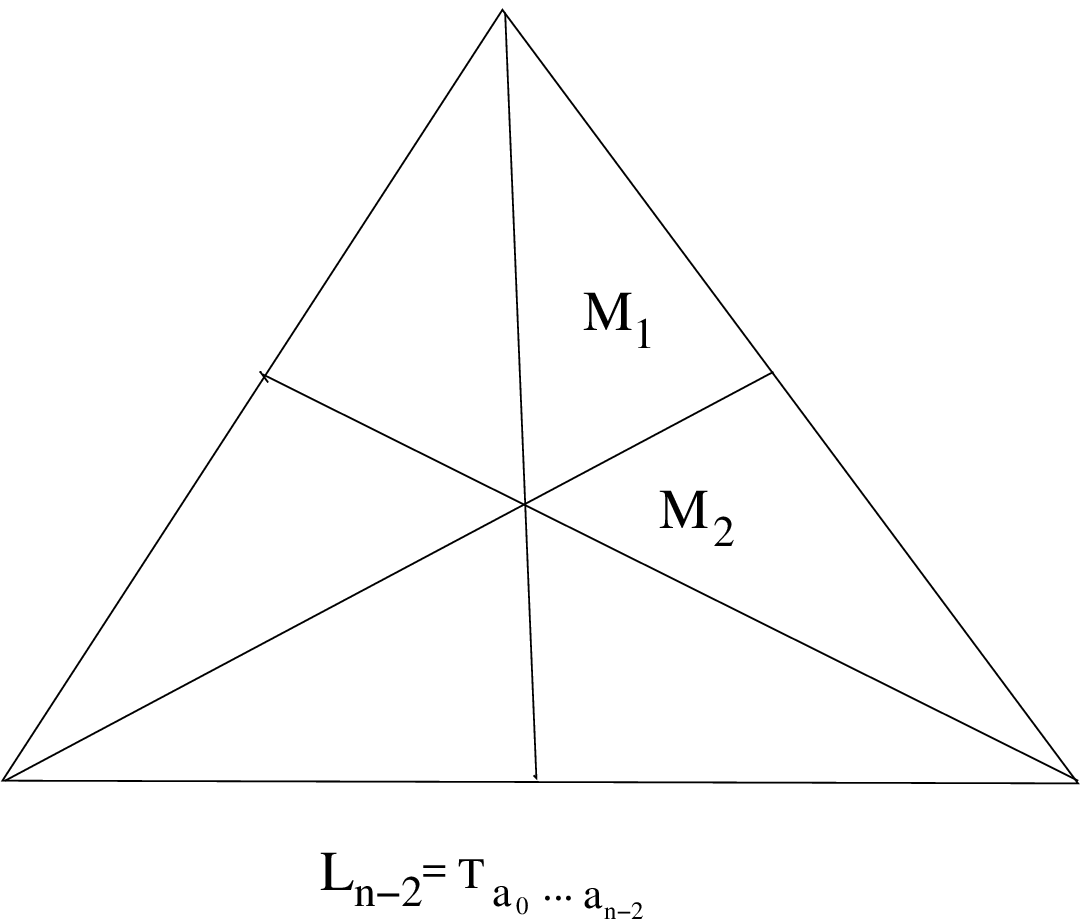}
            \hspace{.5cm}
            \includegraphics[width=.5cm]{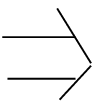}
            \hspace{.5cm}
            \includegraphics[width=1.7in]{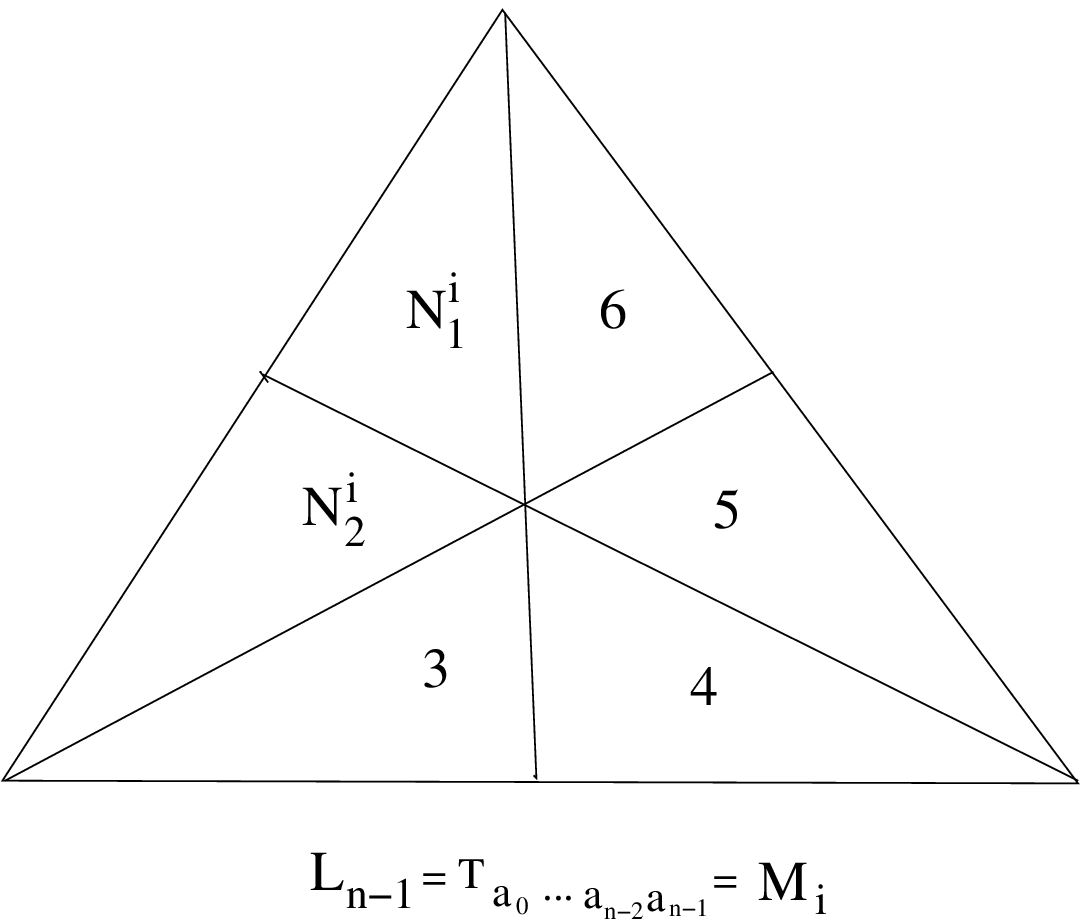}
          }
      \end{figure}

      Recall the assumptions for $ L_{n-1} = M_{i}, $ lying between
      the phrase `Case 2)   $ n \geq 2 $' and the one
      `Subcase 2-1 ) $ a_{n-1} = 1 $, $\cdots$,' and let
      $ L_n := N^i_j, \; j=1, \cdots , 6 .$

      Note the common vertex of $ L_{n-1} = M_{i} $ and $N^{i}_{1}$
      is the starting point of $L_{n-1} = M_{i}$ from the definition
      of $ N^{i}_{1} $. Now, call the vertex the starting point of
      $N^{i}_{1}$. And call the barycenter of $L_{n-1} = M_{i}$ the
      starting point of $N^{i}_{j}$ for $j=2,3,4,5,6.$ Also, call the
      barycenter the ending point of $N^{i}_{j} $ for $j=1,2,3,4,5.$

      Note the common vertex of $ L_{n-1} = M_{i} $ and $ N^{i}_{6} $
      is the starting point of $ L_{n-1} = M_{i} ,$ so the barycenter
      of $ L_{n-2} $ and the ending point of $ L_{n-1} = M_{i} $ from
      the assumption for $M_{i}.$ Call the vertex the ending point of
      $ N^{i}_{6} .$

      Note that the starting and the ending point of $L_{n-1}=M_{i}$
      are same and they are the common vertex of $ L_{n-1} = M_{i} $
      and $ N^{i}_{6} .$

      And the side of $ L_{n-1} = M_{i} $, which contains a side of
      $N^{i}_{1}$, is divided into two line segments, each of which
      is one side of $ N^{i}_{j} $ for $j=1,2,$ respectively.

    \bigskip

    Under the counterclockwise orientation, interior triangles for
    $n=2,3$ will be given as follows :

    \begin{figure}[h]
     \centering{\includegraphics[width=3.5in]{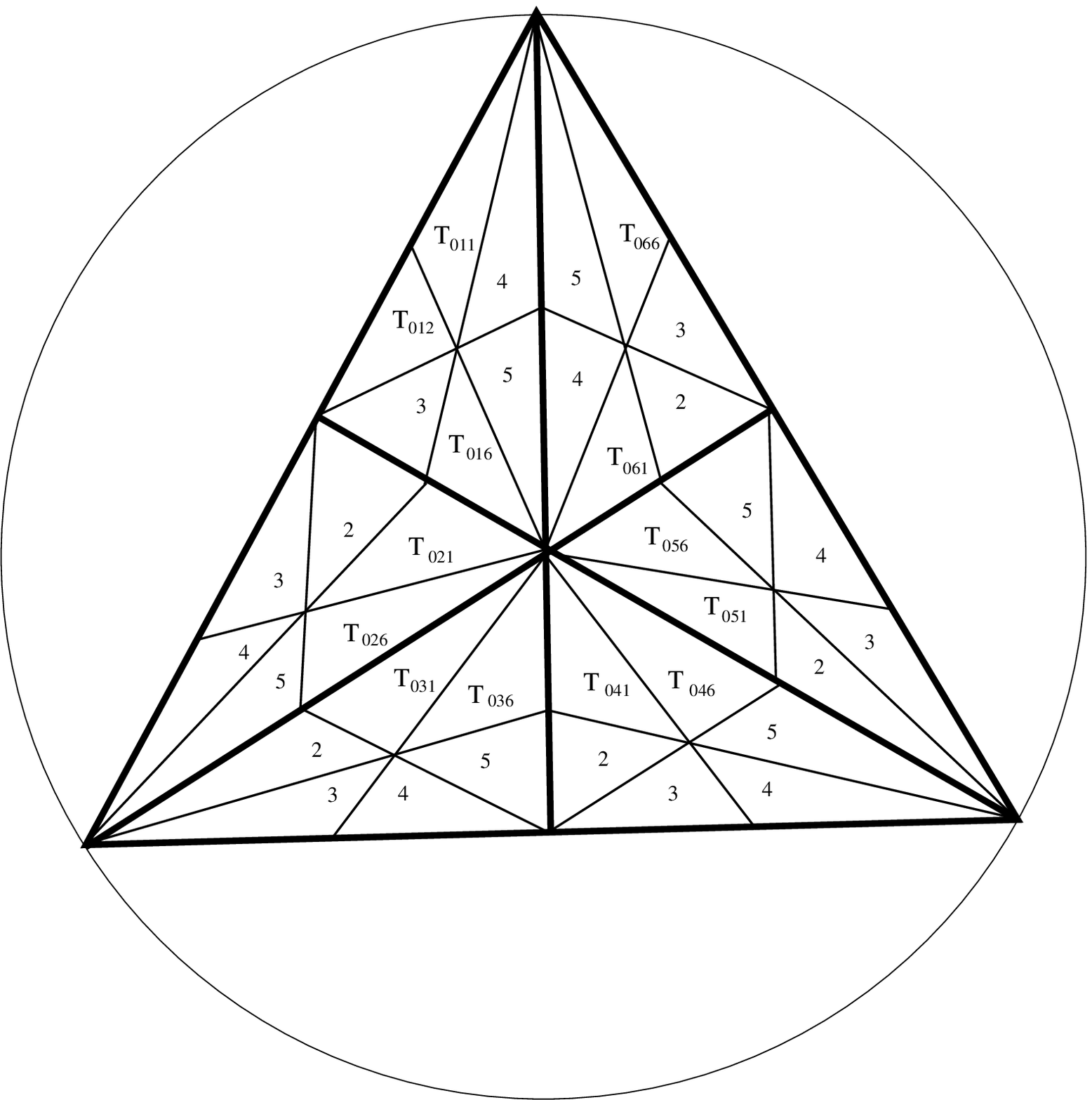}}
    \end{figure}

    \pagebreak

    \begin{figure}[h]
     \centering{\includegraphics[width=3.7in]{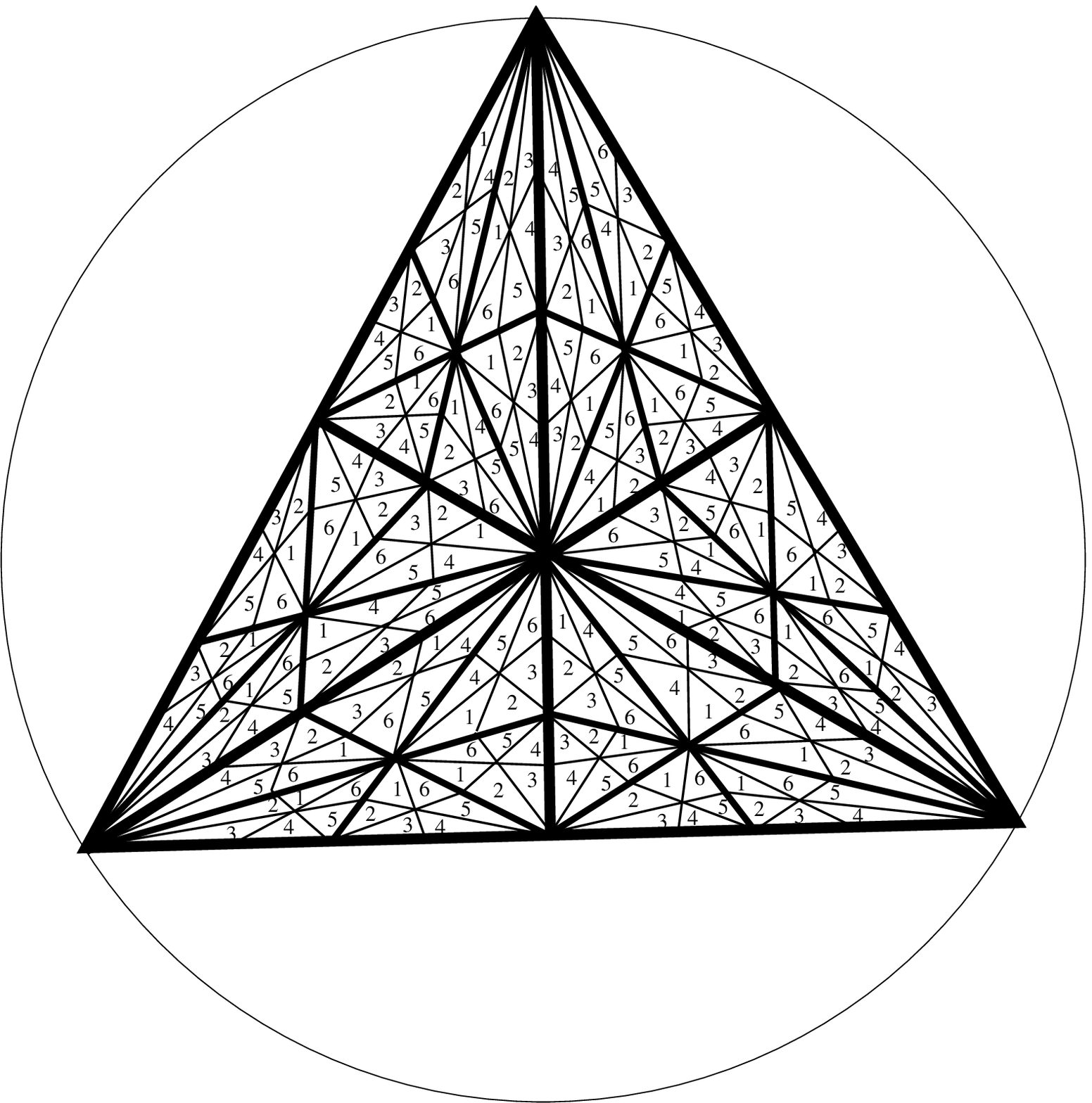}}
    \end{figure}

    Now let's review the definition of the starting points and
    ending points of triangles made right before as follows :

    To begin with, note that the common vertex of $\, T_{a_0 a_1
    \cdots a_{n-1}} $ and $\, T_{a_0 a_1 \cdots a_{n-1} 1} $ is the
    starting point of $ T_{a_0 a_1 \cdots a_{n-1} } $ and that the
    common vertex of $ T_{a_0 a_1 \cdots a_{n-1}} $ and
    $ T_{a_0 a_1 \cdots a_{n-1} 6} $ is the ending point of
    $ T_{a_0 a_1 \cdots a_{n-1} } .$

    For $ T_{a_0 a_1 \cdots a_{n-1} 1} $ , the starting point of
    $T_{a_0 a_1 \cdots a_{n-1} } $ is the common vertex with
    $ T_{a_0 a_1 \cdots a_{n-1} 1} $ and will be called the starting
    point of $ T_{a_0 a_1 \cdots a_{n-1} 1} $. And the barycenter of
    $ T_{a_0 a_1 \cdots a_{n-1} } $ will be called the ending point
    of $ T_{a_0 a_1 \cdots a_{n-1} 1} .$

    For $ T_{a_0 a_1 \cdots a_{n-1} i } $ , where $ i=2,3,4,5 $, the
    barycenter of $ T_{a_0 a_1 \cdots a_{n-1} } $ will be called the
    starting and ending point of $ T_{a_0 a_1 \cdots a_{n-1} i } .$

    For $ T_{a_0 a_1 \cdots a_{n-1} 6 } ,$   the barycenter of
    $T_{a_0 a_1 \cdots a_{n-1} } $ will be called the starting point
    of $ T_{a_0 a_1 \cdots a_{n-1} 6} .$ And the ending point of
    $T_{a_0 a_1 \cdots a_{n-1} } $ is the common vertex with
    $ T_{a_0 a_1 \cdots a_{n-1} 6} $  and will be called the ending
    point of $ T_{a_0 a_1 \cdots a_{n-1} 6 } .$

    \bigskip

    To check whether we can define triangles inductively:

    Recall the assumptions for $ T_{a_0 a_1 \cdots a_{n-1}}, $ lying
    between the phrase `Case 2)   $ n \geq 2 $' and the one
    `Subcase 2-1 ) $ a_{n-1} = 1 $, $\cdots$.'

    Note that the common vertex of $ T_{a_0 a_1 \cdots a_{n-1}} $ and
    $ T_{a_0 a_1 \cdots a_{n-1} 1} $ is the starting point of each of
    them. The barycenter of $ T_{a_0 a_1 \cdots a_{n-1}} $ is the
    starting point of  $ T_{a_0 a_1 \cdots a_{n-1} i} $ for
    $i=2,3,4,5,6,$ and the ending point of
    $ T_{a_0 a_1 \cdots a_{n-1} i} $ for $i=1,2,3,4,5.$ The common
    vertex of $ T_{a_0 a_1 \cdots a_{n-1}} $ and
    $ T_{a_0 a_1 \cdots a_{n-1} 6} $ is the ending point of each of
    them.

    Notice that if the starting and the ending point of
    $ T_{a_0 \cdots a_{n-2} a_{n-1}} $ are same then
    $a_{n-1} \neq 1,6$ and they are the common vertex of
    $ T_{a_0 \cdots a_{n-2} a_{n-1}} $ and
    $ T_{a_0 \cdots a_{n-2} a_{n-1} 6} .$ Also note that if the
    starting and the ending point of
    $ T_{a_0 \cdots a_{n-2} a_{n-1}} $ are different then
    $a_{n-1} \in \{ 1,6 \}$ and $ T_{a_0 \cdots a_{n-2} a_{n-1} 1} $
    and $ T_{a_0 \cdots a_{n-2} a_{n-1} 6} $ are mutually opposite
    ones inside  $ T_{a_0 \cdots a_{n-2} a_{n-1}} .$

    And one side of $ T_{a_0 a_1 \cdots a_{n-1}} $, which contains a
    side of $T_{a_0 a_1 a_2 \cdots a_{n-1} 1}$ is divided into two
    line segments, each of which is one side of
    $ T_{a_0 a_1 \cdots a_{n-1} i} , \text{ for } i=1,2,$
    respectively. Thus we can define triangles inductively.

  \bigskip

  \subsection{The definition of exterior triangles and the definition
  of their starting points and ending points}

  \subsubsection
      {
        \textbf{The definition of $\mathbf{S^{b_0 b_1 \cdots b_n}_0}$}
      }

  The given orientation at the center of $ D^2 $ and the base point,
  or equivalently the starting and ending point of $ T_0 $, will give
  the order $ b_0 $ of sides of $ T_0 $ , where $ b_0 = 1,2,3 $, as
  explained early in `Section 1.'

  \medskip

  Case 1) $ n = 1 $ :

    From the given orientation at the center of $ D^2 $, consider the
    direction of each side of $T_0$, which will give the starting
    point and the ending point of each side.

    For the side $ b_0 $ of $ T_0 $ and the (line) segment on the
    boundary of $ D^2 $ , which faces the side $ b_0 $ and has common
    terminal points with the side $ b_0 $ , consider the midpoint of
    the side $ b_0 $ and of the boundary segment, respectively. Then
    a given half of the side $ b_0 $, the straight line segment
    between the midpoint of the boundary segment and the common
    terminal point of the side $ b_0 $ and of the given half of the
    side $  b_0 $, and the straight line segment between the midpoint
    of the side $ b_0 $ and that of the boundary segment  will
    determine a triangle, so we can obtain two triangles from each
    half of the side $ b_0 $. Let's call them $ S^{b_0 1}_0 $ and
    $ S^{b_0 2}_0 $ , where for $ S^{b_0 i}_0 $ , $ i $ is determined
    by the order with respect to the orientation at the center of
    $ D^2 $  and the line segment connecting the center of $ D^2 $
    and the starting point of $T_0 .$

    Under the counterclockwise orientation, exterior triangles for
    $n=1$ will be given as follows :

    \begin{figure}[h]
     \centering{\includegraphics[width=2.5in]{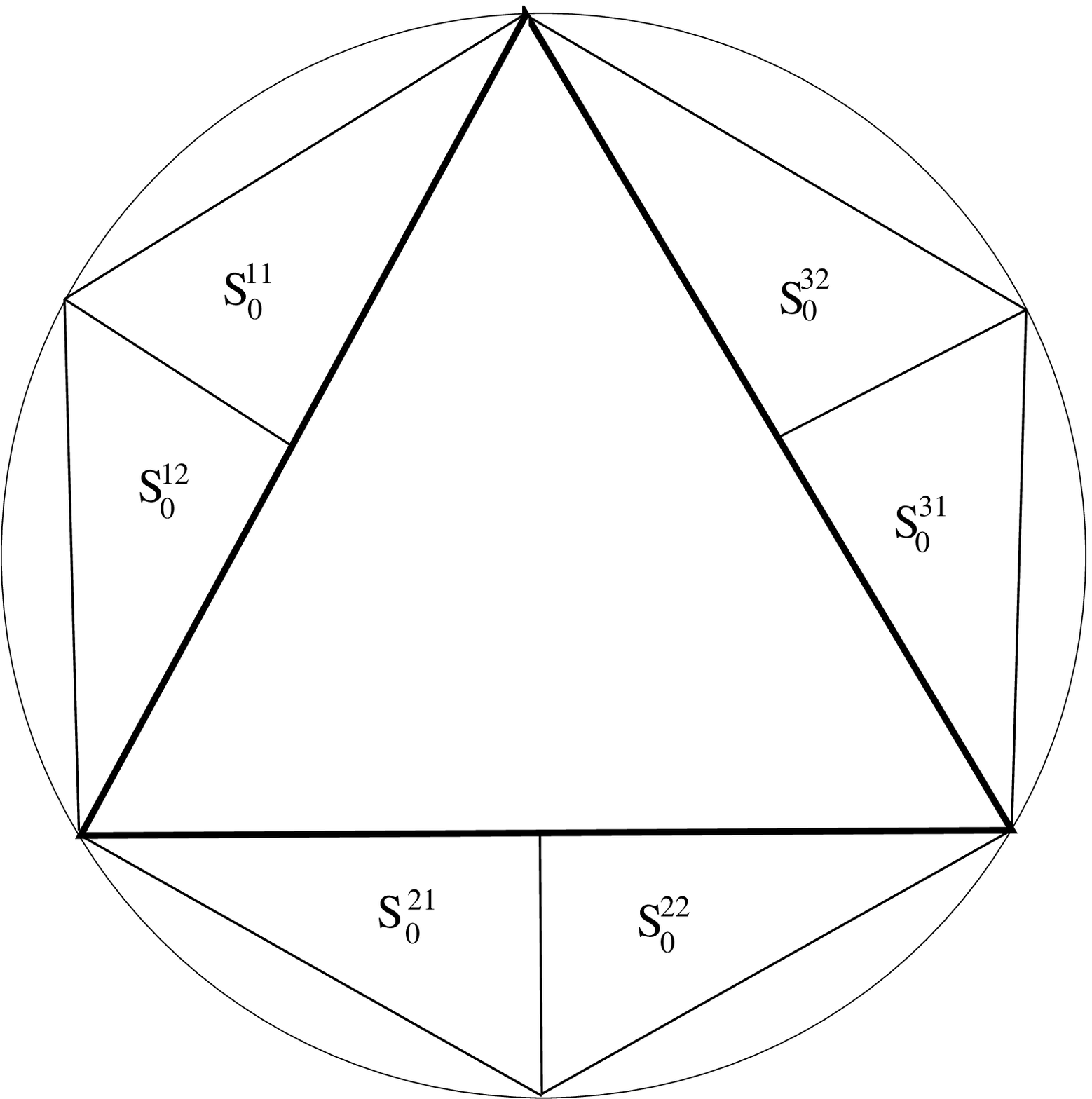}}
    \end{figure}

  \medskip

  Case 2) $ n  \geq 2 $ :

    Let  $ S^{b_0 b_1 \cdots b_{n-1}}_0 $ be given.

    The side of $ S^{b_0 b_1 \cdots b_{n-1}}_0 ,$ which faces the
    boundary of $ D^2 ,$ will give two triangles as follows :

    Consider the direction of the side of
    $ S^{b_0 b_1 \cdots b_{n-1}}_0 ,$ which faces the boundary of
    $ D^2 $,  and that of the line segment on the boundary of
    $ D^2 $, which is being faced by the side, respectively, from the
    orientation at the center of $ D^2 $ and the line segment
    connecting the center of $D^2$ and the starting point of $ T_0 .$
    Then we can think of the starting point, midpoint and ending
    point of the side of $ S^{b_0 b_1 \cdots b_{n-1}}_0 ,$ which
    faces the boundary of $ D^2 ,$ and those of the boundary segment,
    respectively. Now, refer to the construction of two triangles in
    `case 1.' Then the triangle with the midpoints and common
    starting point of the side and the boundary segment as vertices
    will be called $ S^{b_0 b_1 \cdots b_{n-1} 1}_0 $ and the
    triangle with the midpoints and common ending point of the side
    and the boundary segment as vertices will be called
    $ S^{b_0 b_1 \cdots b_{n-1} 2}_0 .$

    \bigskip

    Under the counterclockwise orientation, exterior triangles for
    $n=2,3$ will be given as follows :

    \pagebreak

    \begin{figure}[h]
      \centering{\includegraphics[width=2.9in]{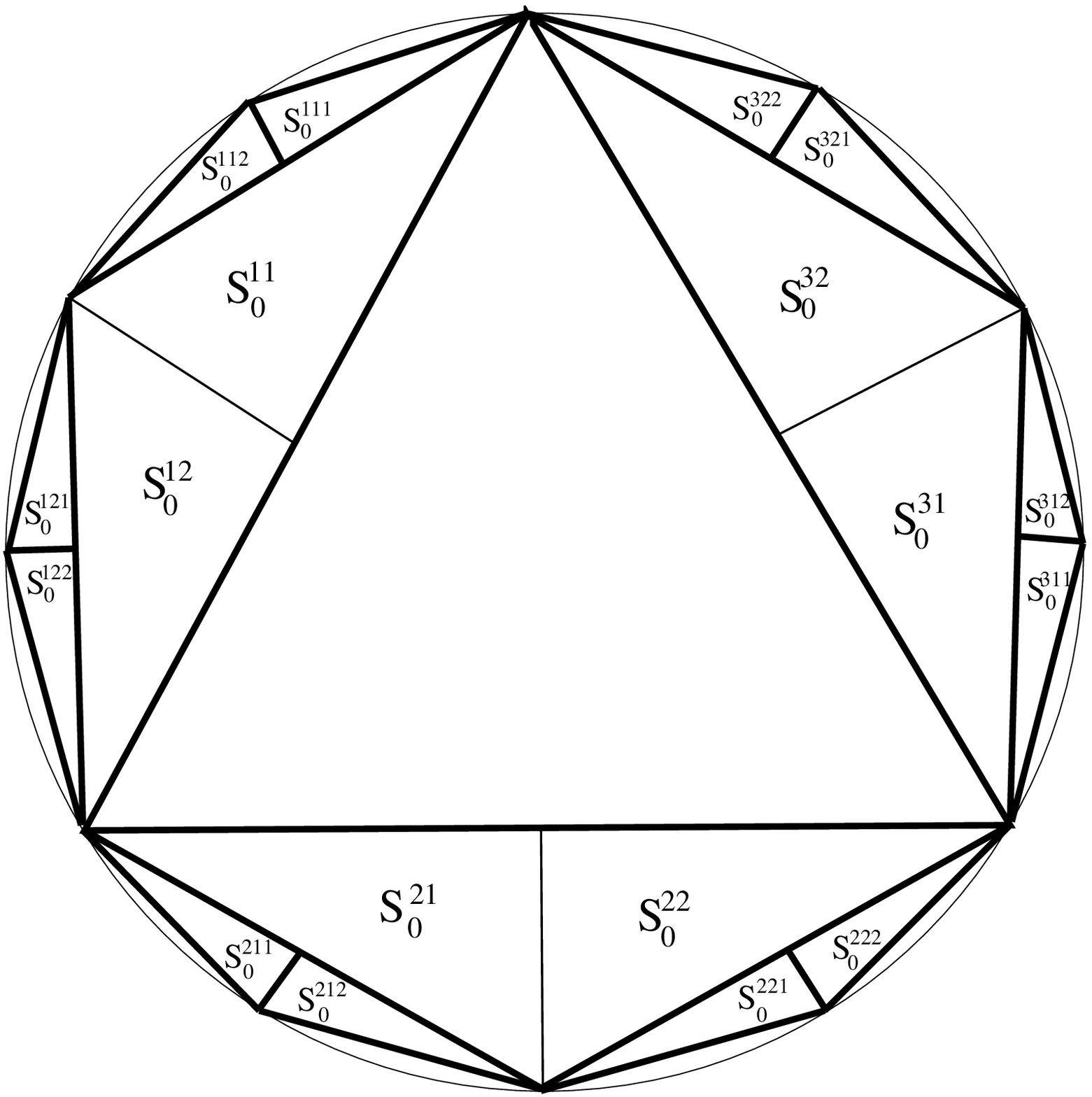}}
    \end{figure}

    \begin{figure}[h]
      \centering{\includegraphics[width=4.4in]{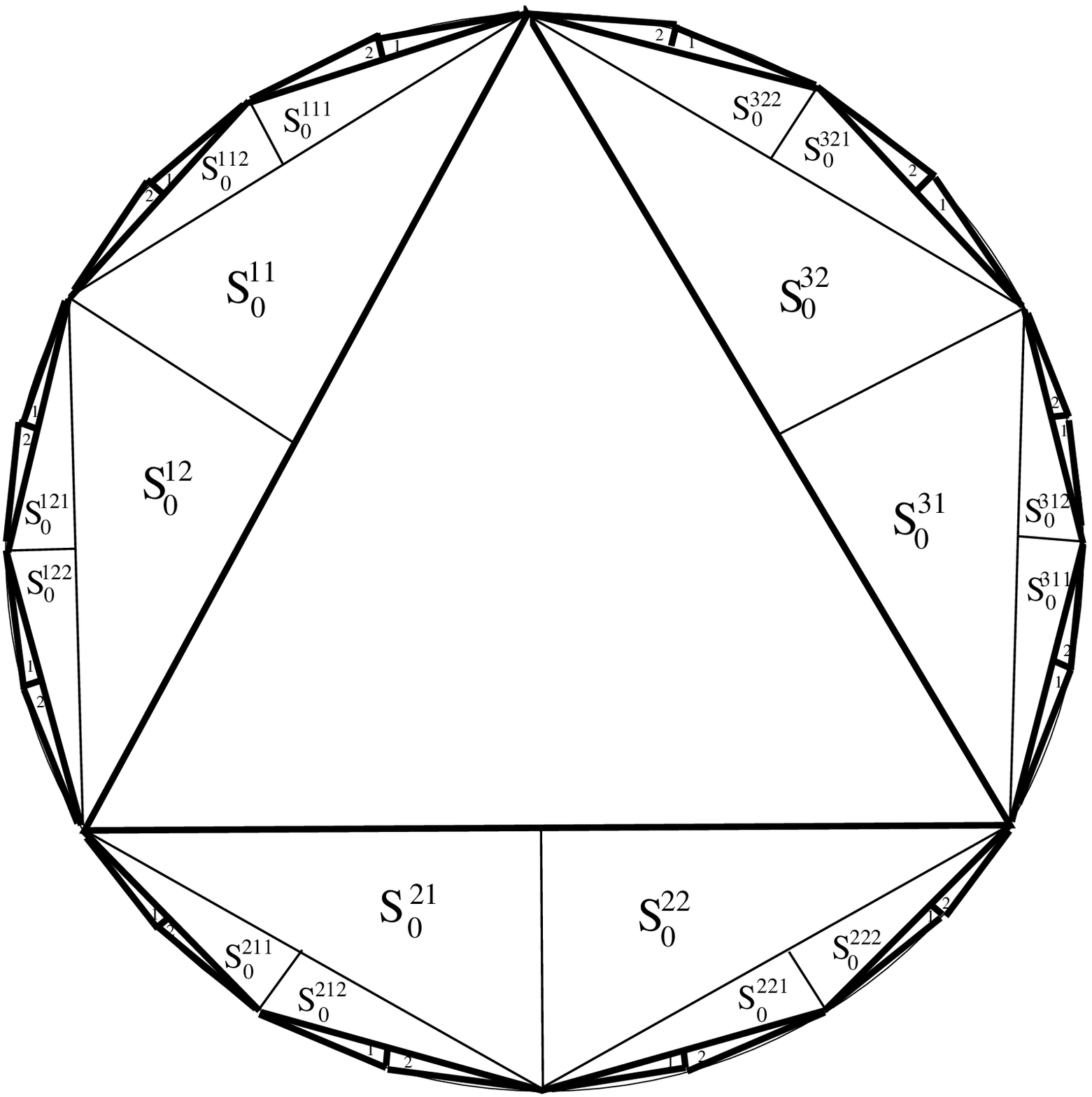}}
    \end{figure}

  \pagebreak

  Now, let's define the starting point and the ending point of the
  triangles made right before as follows:

  \bigskip

  Let $ n \geq 1.$

  For $ S^{b_0 b_1 \cdots b_{n-1} i}_0 , i=1,2, $  consider the
  direction of its side facing the boundary with respect to the
  orientation at the center of $ D^2 .$

  If n is odd,
  the ending point of the side, facing the boundary of $ D^2 $,
  will be called the starting point of $ S^{b_0 b_1 \cdots b_{n-1} i}_0 $
  and
  the starting point of the side, facing the boundary of $ D^2 $,
  will be called the ending point of $ S^{b_0 b_1 \cdots b_{n-1} i}_0 .$

  If n is even,
  the starting point of the side, facing the boundary of $ D^2 $,
  will be called the starting point of $ S^{b_0 b_1 \cdots b_{n-1} i}_0 $
  and
  the ending point of the side, facing the boundary of $ D^2 $,
  will be called the ending point of $ S^{b_0 b_1 \cdots b_{n-1} i}_0 .$

  \vspace{1cm}

  \subsubsection
      {
        \textbf
            {
              The definition of
              $\mathbf{S^{b_0 b_1 \cdots b_k}_{a_0 a_1 \cdots a_m}}$
            }
      }

  Let $ 1 \leq k < n $ be given. Let $ m = n-k $.
  To define
  $S^{b_0 b_1 \cdots b_k}_{a_0 a_1 \cdots a_m}$,
  consider a triangle $\tilde{T}_{0}$ whose orientation is the
  opposite one of $ T_{0} $ (, considering $ T_{014} $ might be
  helpful). Then the $m$-step barycentric subdivision makes us think of
  $ \tilde{T}_{0 a_1 \cdots a_m} $ , which is the mirror-symmetry
  of $ T_{0 a_1 \cdots a_m} $
  (, for example $T_{0141 a_1 \cdots a_m}$).
  Note the orientation of the triangle $\tilde{T}_{01}$ is the
  opposite one of $ T_{01} $ (,considering $ T_{0141} $ might be
  helpful), and its $m$-step barycentric subdivision
  $ \tilde{T}_{0 1 a_2 \cdots a_m} $ is also the
  mirror-symmetry of $ T_{01 a_2 \cdots a_m} $
  (, for example $T_{0141 a_2 \cdots a_m}$).

  We want to define
  $ S^{b_0 b_1 \cdots b_k}_{a_0 a_1 \cdots a_m} $ as follows :

  \medskip

  Case 1-1) $k$ is odd  and $b_k = 1$:

    Consider the barycentric subdivision of $S^{b_0 b_1 \cdots b_k}_0.$
    By comparing it with that of $T_{01},$ define

    - $S^{b_0 b_1 \cdots b_k}_{0j}$, which matches $T_{01j}$ \:
    for $j \in \{1,2,3,4\}$,

   - $S^{b_0 b_1 \cdots b_k}_{05}$, which matches $T_{016}$,

   - $S^{b_0 b_1 \cdots b_k}_{06}$, which matches $T_{015}$,

   \noindent
   and their starting and ending points.

    For $m \geq 2$, the respective identification of
    $$
      S^{b_0 b_1 \cdots b_k}_{01}, S^{b_0 b_1 \cdots b_k}_{0j},
      S^{b_0 b_1 \cdots b_k}_{04}, S^{b_0 b_1 \cdots b_k}_{06}
      \: \text{ with } \:
      T_{01}, T_{0}, \tilde{T}_0, \tilde{T}_{01},
    $$
    where $j \in \{2,3,5\},$ and their $m$-step
    barycentric subdivision can make us define
    $ S^{b_0 b_1 \cdots b_k}_{a_0 a_1 a_2 \cdots a_m}, $
    where $ a_0 = 0 .$

  \medskip
  Case 1-2) $k$ is odd  and $b_k = 2$:

    Identify $ S^{b_0 b_1 \cdots b_k}_0 $ with $ T_{01} ,$  where the
    starting point and ending point  of $ S^{b_0 b_1 \cdots b_k}_0 $
    is also identified to those of $ T_{01} .$

    Consider the m-step barycentric subdivision of
    $ S^{b_0 b_1 \cdots b_k}_0 $ and $ T_{01} $ respectively. The
    identification, then, can make us define
    $ S^{b_0 b_1 \cdots b_k}_{a_0 a_1 \cdots a_m} $ from
    $ T_{01 a_1 \cdots a_m} ,$ where $ a_0 = 0 .$

    \medskip
    Under the counterclockwise orientation, the triangles for $k=1$
    and $m=1$ will be given as follows :

    \begin{figure}[h]
      \centering{\includegraphics[width=3.2in]{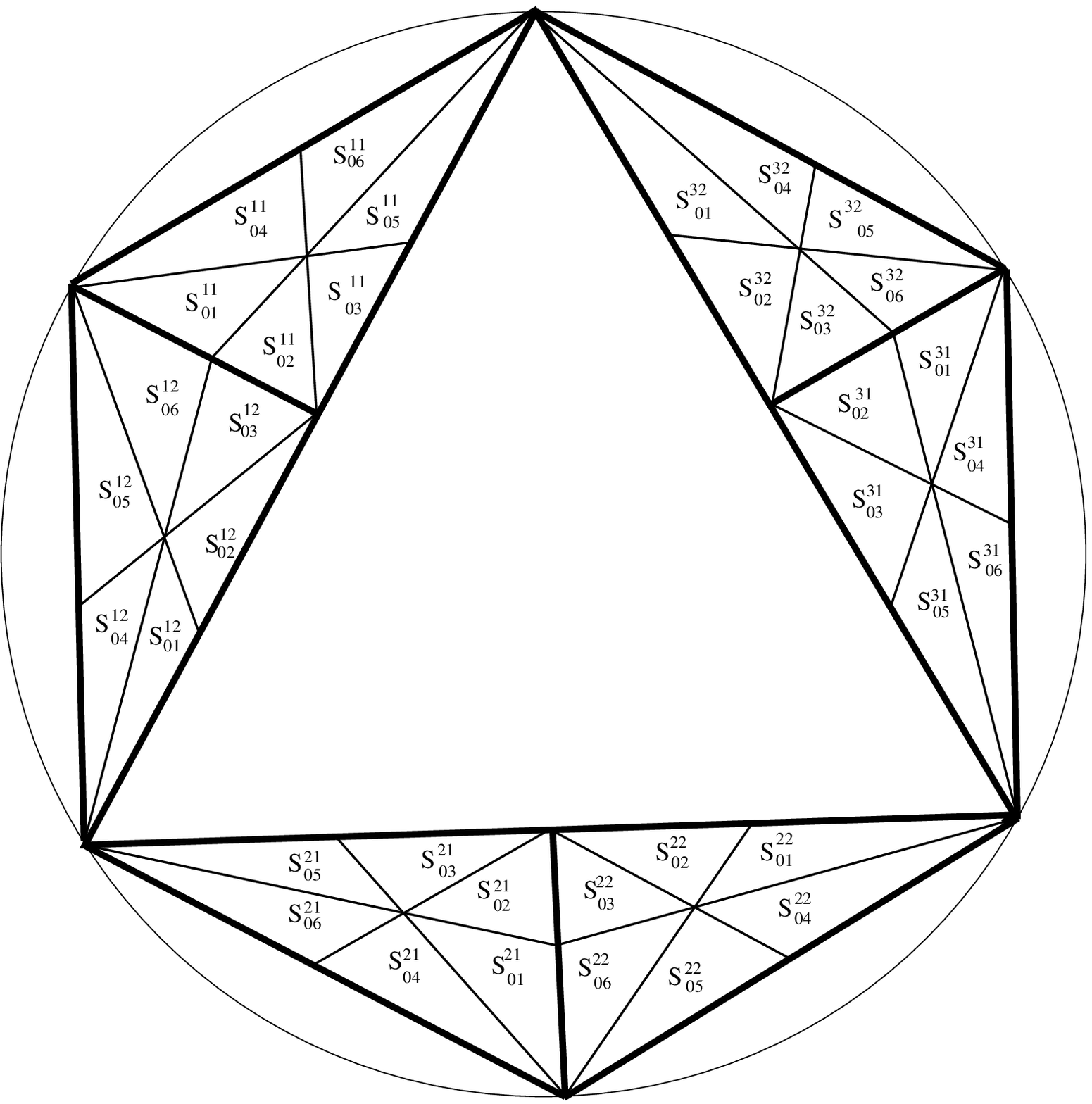}}
    \end{figure}

  \medskip
  Case 2-1) $k$ is even and $b_k=1$:

    By identifying
    $ S^{b_0 b_1 \cdots b_k}_{0} $ with $\tilde{T}_{01}$, we can
    define $ S^{b_0 b_1 \cdots b_k}_{a_0 a_1 \cdots a_m} $ from
    $ \tilde{T}_{01 a_1 \cdots a_m} $, where $ a_0 = 0 $ (, for
    example $ T_{0141 a_1 \cdots a_m} $) .

  \medskip

  Case 2-2) $k$ is even  and $b_k = 2$:

    Consider the barycentric subdivision of $S^{b_0 b_1 \cdots b_k}_0.$
    By comparing it with that of $\tilde{T}_{01},$ define

    - $S^{b_0 b_1 \cdots b_k}_{0j}$, which matches
    $\tilde{T}_{01j}$ \: for $j \in \{1,2,3,4\}$,

   - $S^{b_0 b_1 \cdots b_k}_{05}$, which matches $\tilde{T}_{016}$,

   - $S^{b_0 b_1 \cdots b_k}_{06}$, which matches $\tilde{T}_{015}$,

   \noindent
   and their starting and ending points.

    For $m \geq 2$, the respective identification of
    $$
      S^{b_0 b_1 \cdots b_k}_{01}, S^{b_0 b_1 \cdots b_k}_{0j},
      S^{b_0 b_1 \cdots b_k}_{04}, S^{b_0 b_1 \cdots b_k}_{06}
      \: \text{ with } \:
      \tilde{T}_{01}, \tilde{T}_{0}, T_0, T_{01},
    $$
    where $j \in \{2,3,5\},$ and their $m$-step
    barycentric subdivision can make us define
    $ S^{b_0 b_1 \cdots b_k}_{a_0 a_1 a_2 \cdots a_m}, $
    where $ a_0 = 0 .$

    \medskip
    Under the counterclockwise orientation, the triangles for $k=2$
    and $m=1$ will be given as follows :

    \begin{figure}[h]
      \centering{\includegraphics[width=3.5in]{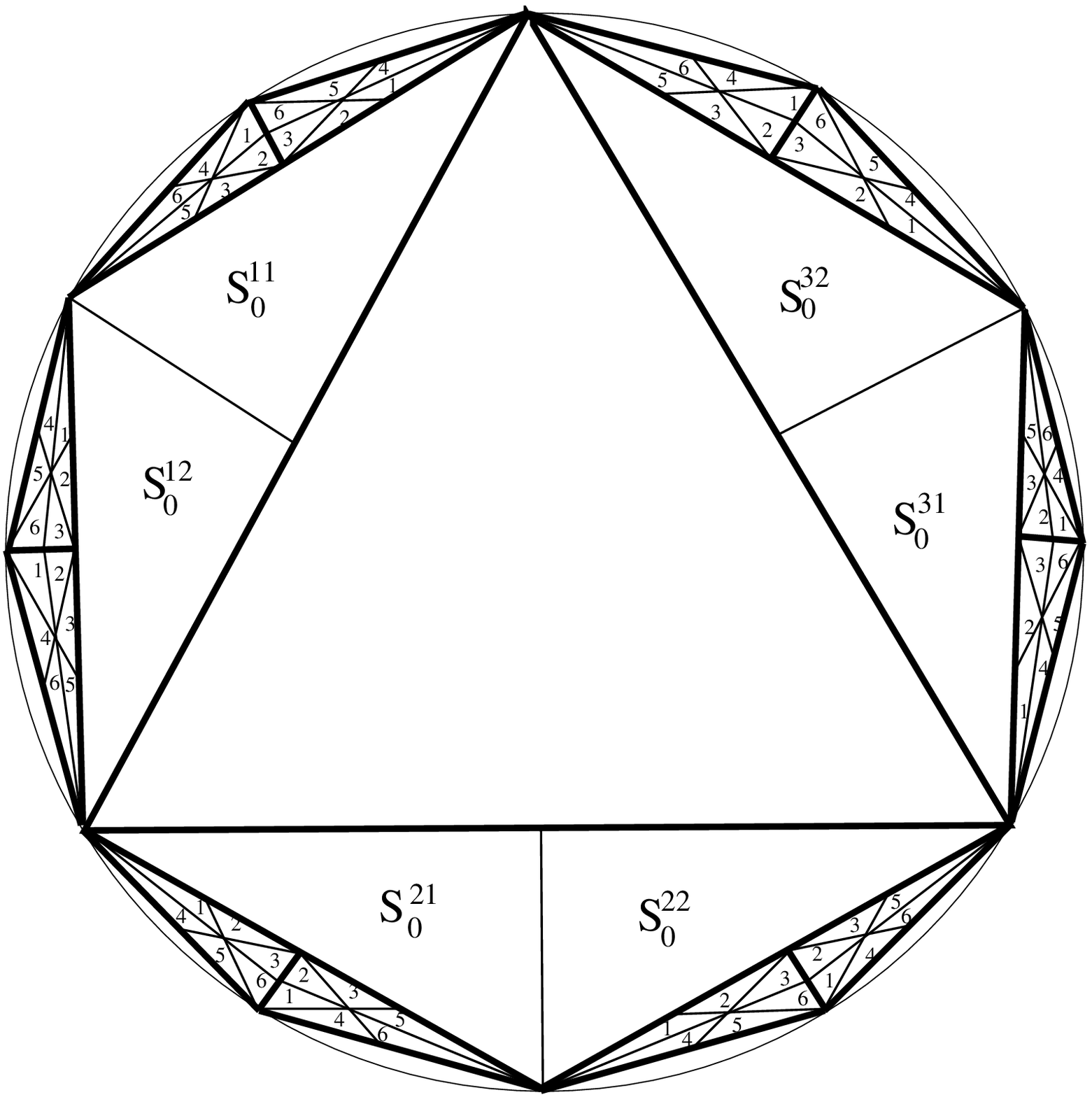}}
    \end{figure}

  \vspace{1cm}

  \subsection{The ordering of triangles in the $ n-th $ step}

  For $ n = 1,2, \cdots  $ , let

  \begin{align*}
    A_n =
    \{
      T_{a_0 a_1 \cdots a_n} \mid
      a_0=0, a_i \in \{ 1,2,3,4,5,6 \} \textrm{ for } i=1, \cdots , n
    \}
    \bigcup
    \\
    (
      \cup _{k+m=n, 1 \leq k \leq n , 0 \leq m \leq n-1 }
      \{
        S ^{ b_0 b_1 \cdots b_k }_{ c_0 \cdots c_m } \mid
        b_0 \in \{ 1,2,3 \} , b_i \in \{ 1,2 \} \textrm{ for }
        i=1, \cdots k ,
        \\
        c_0 =0 , c_j \in \{ 1, \cdots , 6 \} \textrm{ for }
        1 \leq j \leq m \textrm{ if } m \geq 1
      \}
    )
  \end{align*}

  \noindent
  , which is regarded as the set of all triangles in the $ n-th $
  step.

  \vspace{1cm}

  Now refer to the following pictures for 0th, 1st, 2nd and 3rd step
  under the counterclockwise orientation :

    \begin{figure}[h]
      \centering{
                 \includegraphics[width=2in]{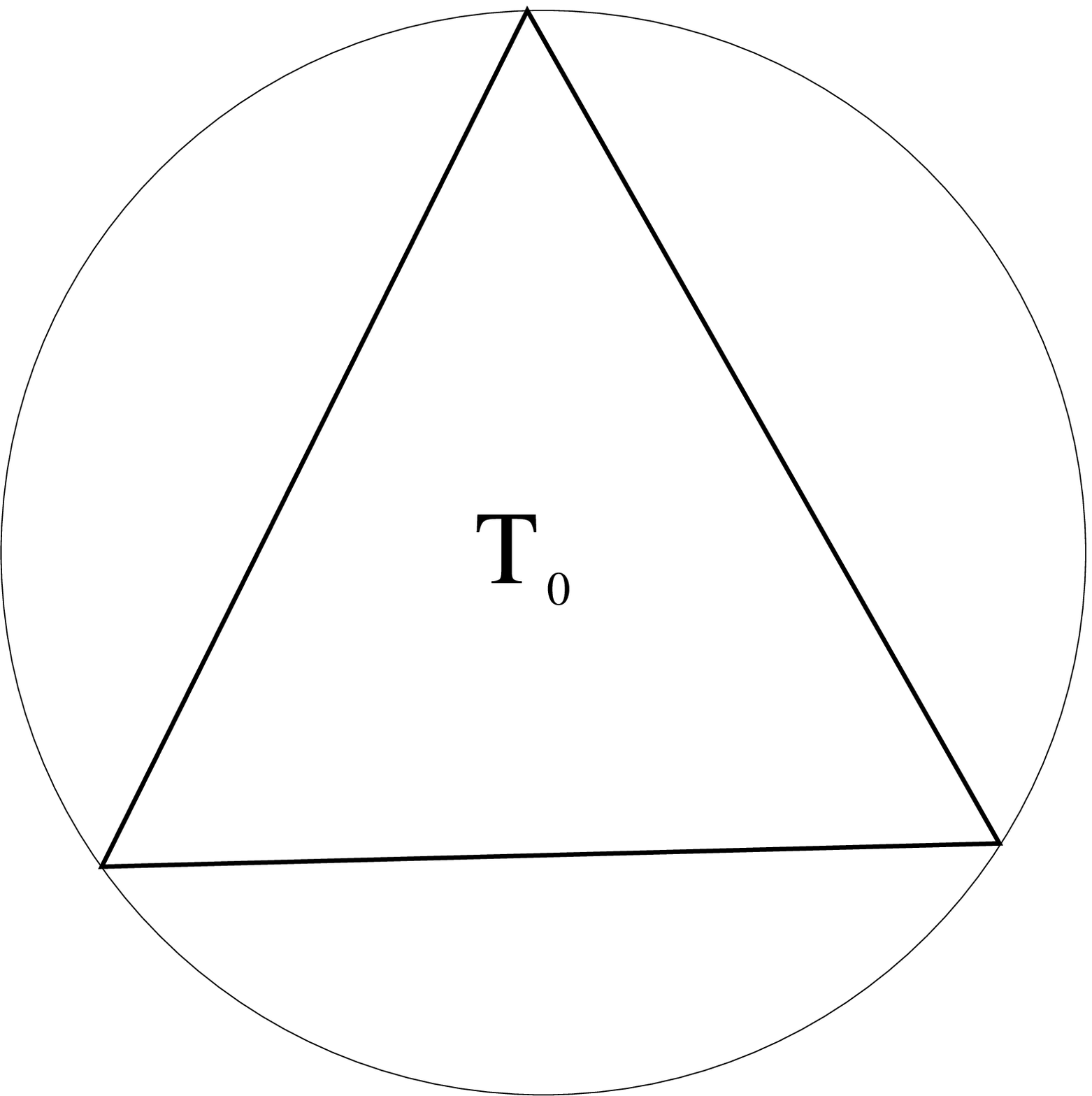}
                 \includegraphics[width=2in]{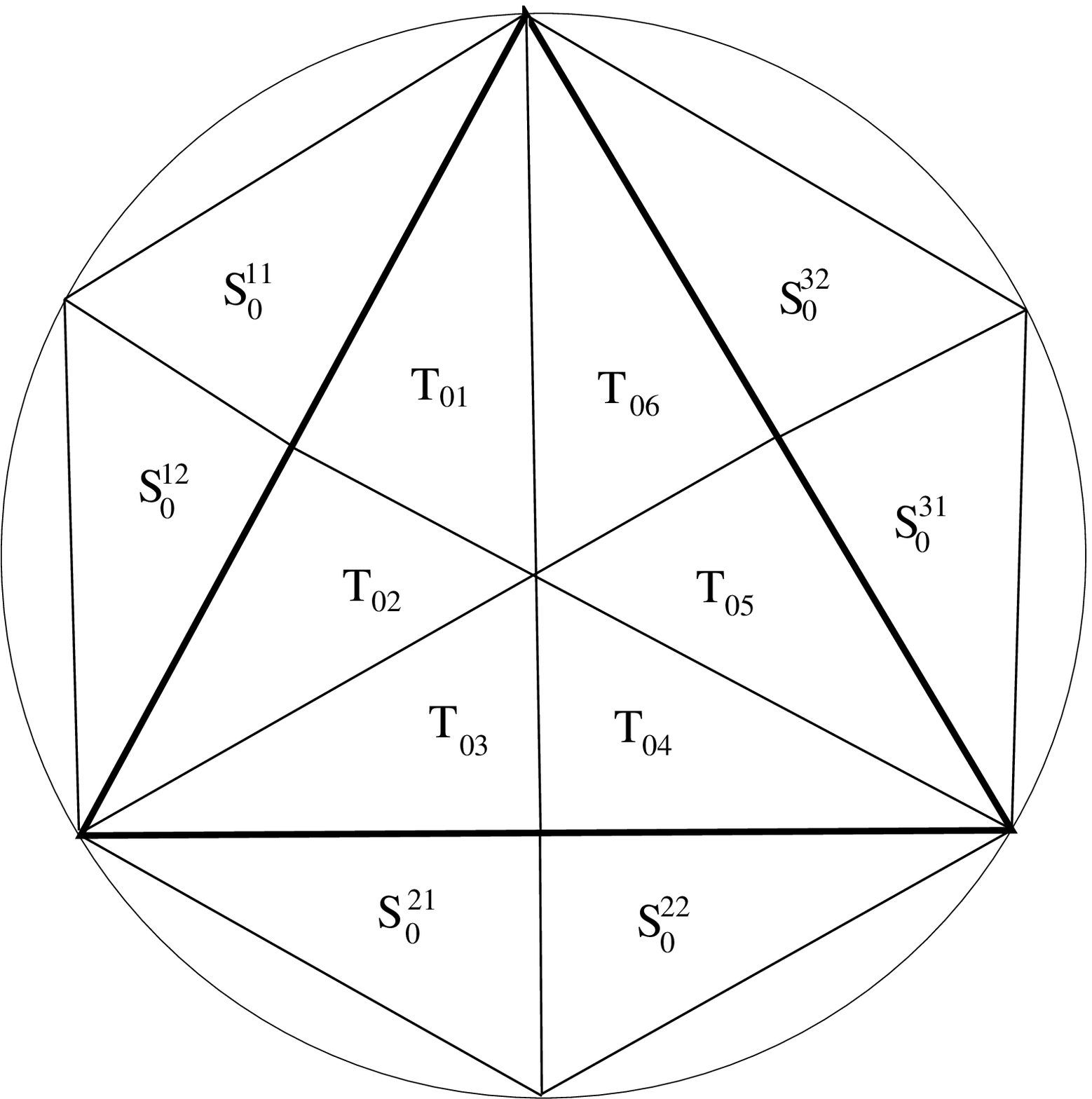}
                }
    \end{figure}

    \begin{figure}[h]
      \centering{\includegraphics[width=5in]{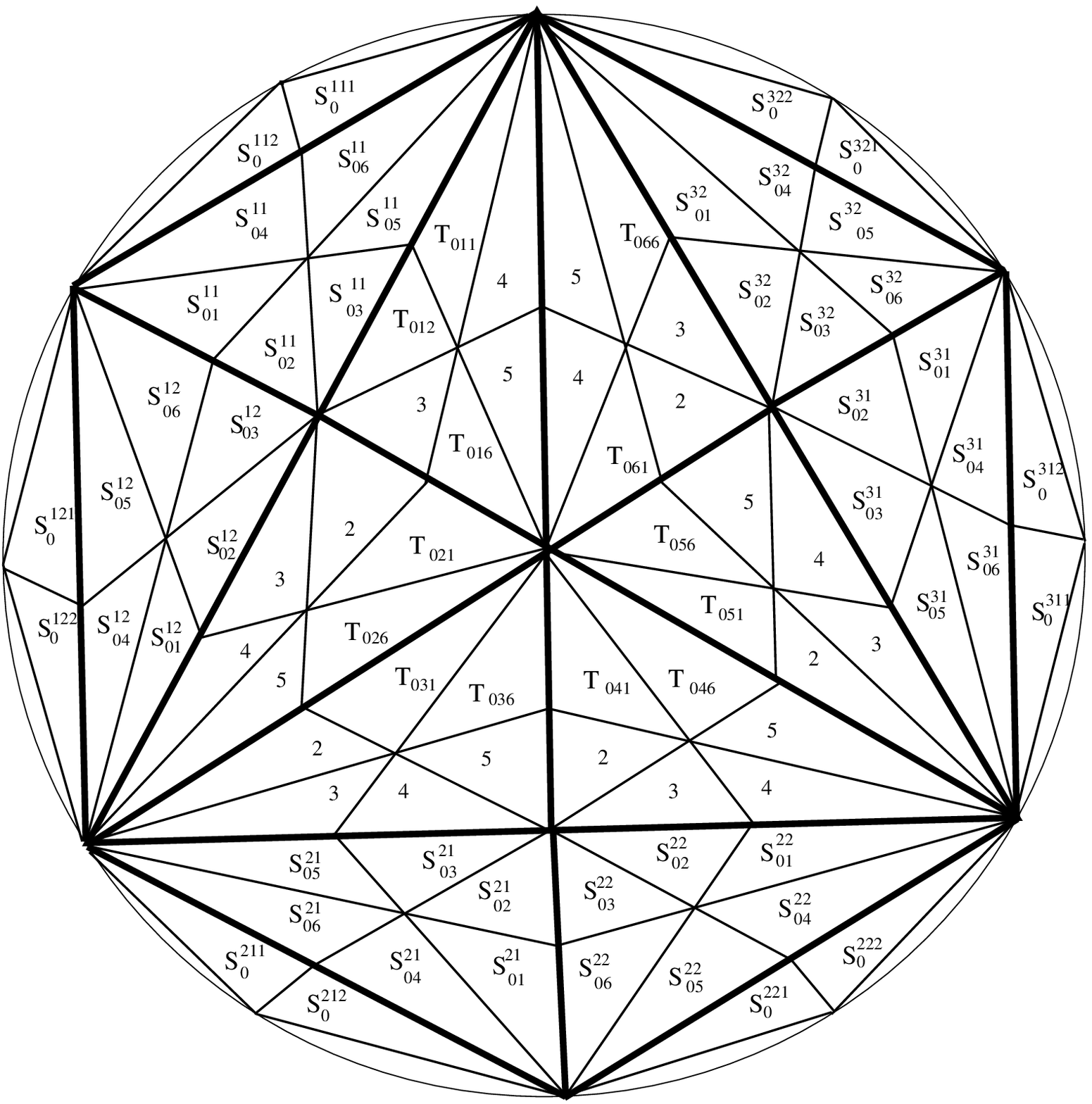}}
    \end{figure}

    \begin{figure}[h]
      \centering{\includegraphics[width=5in]{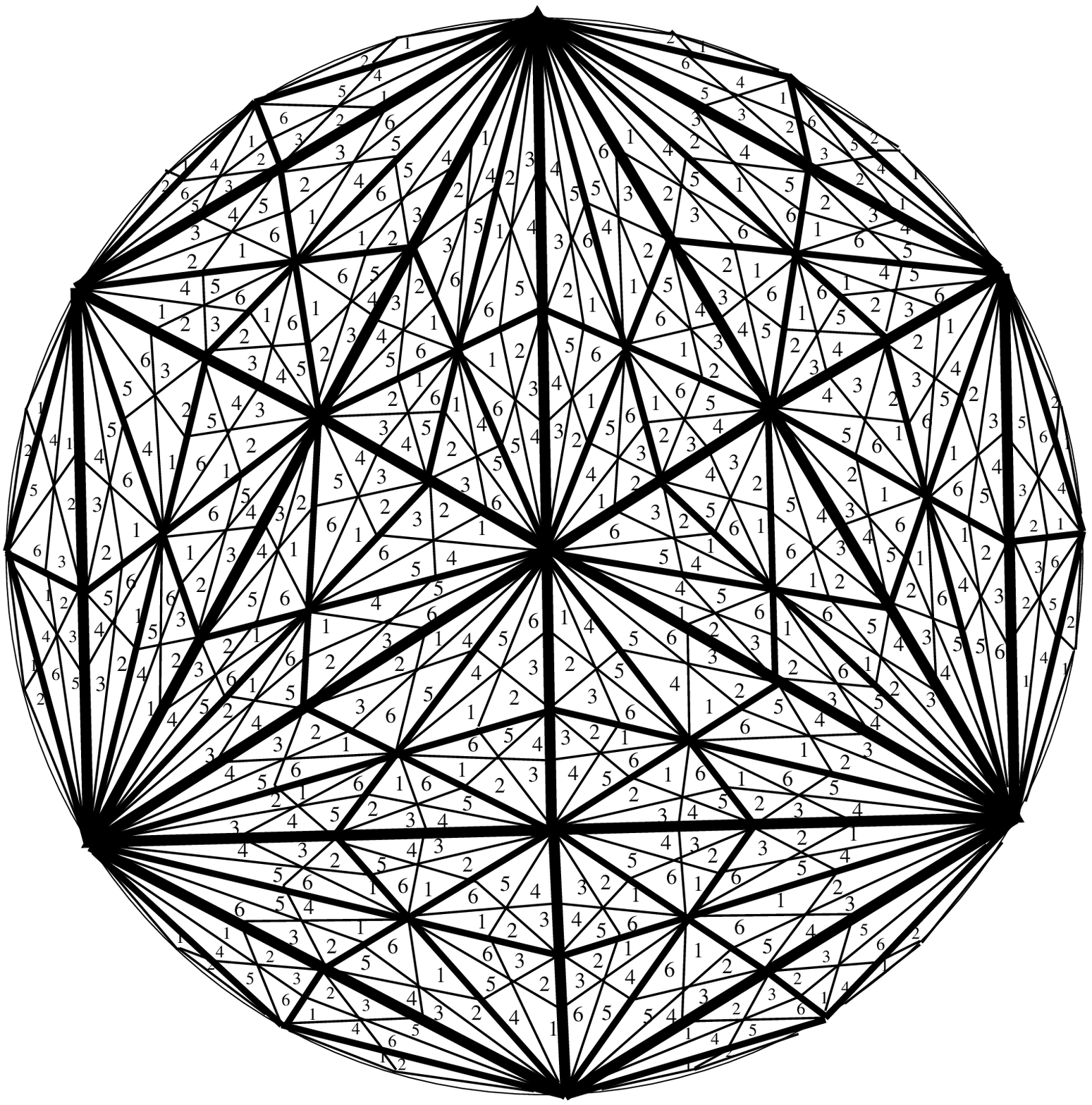}}
    \end{figure}

  Case 1)  $
            T_{a_0 \cdots a_n} < S^{b_0 b_1 \cdots b_k}_{c_0 \cdots
            c_m}
           $
           , where $ k+m = n $

  \bigskip

  Case 2)  $ T_{a_0 \cdots a_n} < T_{b_0 \cdots b_n} $  if
           $ (a_0, \cdots , a_n) < (b_0, \cdots , b_0) $
           with respect to the dictionary order

  \bigskip

  Case 3) The order of $ S^{b_0 b_1 \cdots b_k}_{a_0 \cdots a_m} $
          and $ S^{c_0 c_1 \cdots c_s}_{d_0 \cdots c_t} $
          , where $ k+m  = n = s+t $

    \bigskip

    Case 3-1) $ k<s $ :
            $$
              S^{b_0 b_1 \cdots b_k}_{a_0 \cdots a_m}
              < S^{c_0 c_1 \cdots c_s}_{d_0 \cdots d_t}
            $$

    \bigskip

     Case 3-2) $ k=s $ ( so, $ m=t $ ) and
               $(b_0, b_1, \cdots , b_k) < (c_0, c_1, \cdots , c_k)$
               with respect to the dictionary order :

            If $ k $ is odd, then
            $
              S^{b_0 b_1 \cdots b_k}_{a_0 \cdots a_m}
              > S^{c_0 c_1 \cdots c_k}_{d_0 \cdots d_m}
            $

            If $ k $ is even, then
            $
              S^{b_0 b_1 \cdots b_k}_{a_0 \cdots a_m}
              < S^{c_0 c_1 \cdots c_k}_{d_0 \cdots d_m}
            $

    \bigskip

    Case 3-3) $ k=s $,
              $(b_0, b_1, \cdots , b_k) = (c_0, c_1, \cdots , c_k)$
              and $(a_0, \cdots , a_m) < (d_0, \cdots , d_m)$
              with respect to the dictionary order :

            $$
              S^{b_0 b_1 \cdots b_k}_{a_0 \cdots a_m}
              < S^{c_0 c_1 \cdots c_k}_{d_0 \cdots d_m}
            $$

  \vspace*{1cm}

  \subsection{The properties of triangles in $ A_n $}
  \label{sec:prop-of-triangles}

  We can easily check the following three properties from the
  definition of triangles.

   \bigskip

  Property 1.) Given a non-first element $L$ in $A_n$, the boundary
  of $ \bigcup \{ M \in A_n | M < L \}  $ contains a side of L, which
  will be divided into two line segments in its barycentric
  subdivision, where one of two line segments will become a side of
  the first triangle and the other one will become a side of the
  second triangle in the barycentric subdivision of L.

  \bigskip

  Property 2.) Given $L \in A_n$ ,
  $ \bigcup \{ M \in A_n | M \leq L \}  $ is diffeomorphic to the
  disk $D^2$.
  \bigskip

  Property 3.)
  Assume $L \in A_n$ and six triangles
  $M_1, M_2, \cdots, M_6 \in A_{n+1}$, obtained from the
  barycentric subdivision of $L$, follows the order
  of $i=1,2, \cdots 6$ in $A_{n+1}$. Then the starting points of
  $M_1$ and $L$ are same. Also are the ending points of $M_6$ and
  $L$.

  \bigskip

  And we also have the next property :

  \bigskip

  Property 4.)
  Assume $ L, M \in A_n $ and that $M$ is the next element of $L$ in
  $A_n$ for $n \geq 1$.

  Then,  The ending point of $L$ and the starting point of $M$ are
  same.

  \bigskip

  Proof )

  Case 1) $ L = T_{a_0 \cdots a_{n-1} a_n} $ for some
          $(a_0, \cdots , a_{n-1} , a_n)$

    Subcase 1-1 ) $a_n \neq 6$

      Note $ M = T_{b_0 \cdots b_{n-1} b_n} $ , where

      $$
        b_i = a_i \text{ for } 0 \leq i < n \text{ and } b_n  = a_n +1 .
      $$

      Then inside the triangle $ T{a_0 \cdots a_{n-1}} $, the
      barycenter of $ T_{a_0 \cdots a_{n-1} } $ is the ending point
      of $ L = T_{a_0 \cdots a_{n-1} a_n} $ and the starting point of
      $ M = T_{b_0 \cdots b_{n-1} b_n} $ at the same time.

    \bigskip

    Subcase 1-2 ) $ a_n = 6 $

      If $ a_0 = 0  $ and $ a_1 = \cdots = a_n = 6 $, then the ending
      point of $ L = T_{a_0 \cdots a_{n-1} a_n} $ is the ending point
      of $T_{06}$ by induction and also the ending point of $T_0,$
      that is, the basepoint, which is the starting point of
      $S^{32}_0$ and so the starting point of $ M = S^{32}_0 $ if
      $n=1$ and the starting point of
      $ M= S^{32}_{b_0 \cdots b_{n-1}} $ with $b_0 = 0$ and
      $b_1 = \cdots = b_{n-1} = 1$ if $n \geq 2 .$

      Now assume $n \geq 2$ and $a_i \neq 6$ for some $i$ with
      $1 \leq i < n.$

      We can find $i_0$ satisfying
      $ 1 \leq i_0 < n, \; a_{i_0} \neq 6 $ and $a_i = 6$ for all
      $i_0 < i \leq n.$ Then $ M = T_{b_0 \cdots b_{n-1} b_n} $
      satisfies

      \begin{align*}
        b_i \text{ } &= a_i \text{ for all } 0 \leq i < i_0    \\
        b_{i_0} &= a_{i_0} + 1                                 \\
        b_i \text{ } &= 1   \text{ for all } i_0 < i \leq n
      \end{align*}

      Note the ending point of $ L = T_{a_0 \cdots a_{n-1} a_n} $ is
      the ending point of $ T_{a_0 \cdots a_{i_0} } $ by induction.

      Notice the starting point of $ M = T_{b_0 \cdots b_{n-1} b_n} $
      is the starting point of $ T_{b_0 \cdots b_{i_0} } $ by
      induction.

      Since $a_i = b_i $ for $0 \leq i < i_0$ and
      $b_{i_0} = a_{i_0} + 1$, the ending point of
      $ L = T_{a_0 \cdots a_{i_0 -1} a_{i_0}} $ is the barycenter of
      $ T_{a_0 \cdots a_{i_0 -1} } $, which is the starting point of
      $
        T_{a_0 \cdots a_{i_0 -1} b_{i_0}}
        = T_{b_0 \cdots b_{i_0 -1} b_{i_0}} = M .
      $
      Thus, we get

      $$
        \text{the ending point of } L
        \text{ is the starting point of } M.
      $$

  \bigskip

  Case 2) $L = S^{b_0 b_1 \cdots b_k}_{c_0 \cdots c_m}$ where $k+m=n.$

    \bigskip

    Subcase 2-1 ) $m=0$ and $k=n$ is odd.

      \bigskip

      Note $ (b_0, b_1, \cdots , b_n) \neq (1, 1, \cdots , 1) $,
      because if $ (b_0, b_1, \cdots , b_n) = (1, 1, \cdots , 1) $
      then $ L = S^{b_0 b_1 \cdots b_n}_0 = S^{1 1 \cdots 1}_0 $ is
      the last element in $A_n$.

      If $n=1$, then we can trivially obtain that the ending point of
      $L$ is the starting point of $M$ from the definition of
      triangles.

      Assume $ n \geq 2 $.

      If $ M = S^{d_0 d_1 \cdots d_n}_0 $, then we get
      $(b_0, b_1, \cdots, b_n) > (d_0, d_1, \cdots, d_n)$ and so

      $$
        \text{ either }
        (
          d_0 = b_0 -1, d_1 = \cdots = d_n = 2
          \text{ and }
          b_1 = \cdots = b_n = 1
        )
      $$

      or

      $$
        \exists i_0  \text{ with } 1 \leq i_0 \leq n
        \text{ such that }
      $$

      \begin{align*}
        d_i  \: &= b_i \text{ for all } 0 \leq i < i_0      \\
        d_{i_0} &= 1 , b_{i_0} = 2                          \\
        d_i  \: &= 2, b_i = 1    \text{ for all } i_0 < i \leq n
        \text{ if } 1 \leq i_0 < n .
      \end{align*}

      In the first possibility, the side of $L$ and the side of $M,$
      both of which faces the boundary, meet at one point where
      $ S^{b_0 b_1}_0 = S^{b_0 1}_0 $ and
      $ S^{d_0 d_1}_0 = S^{d_0 2}_0 $ meet.

      In the second possibility, the side of $L$ and the side of $M,$
      both of which faces the boundary, meet at one point which
      is contained in such line segment as the intersection of
      $
        S^{b_0 b_1 \cdots b_{i_0 -1} b_{i_0}}_0
        = S^{b_0 b_1 \cdots b_{i_0 -1} 2}_0
      $
      and
      $
         S^{b_0 b_1 \cdots b_{i_0 -1} 1}_0
         = S^{d_0 d_1 \cdots d_{i_0 -1} d_{i_0}}_0 .
      $

      In any possibilities, the side of $L$ and the side of $M$,
      both of which faces the boundary, meet at one point.
      Since $n$ is odd, the point
      is the starting point of the side of $L$, facing the boundary,
      and the ending point of the side of $M$, facing the boundary, so

      $$
        \text{the ending point of } L
        \text{ is the starting point of } M.
      $$

    \bigskip

    Subcase 2-2 ) $ m = 0 $ and $ k=n $ is even

      \bigskip

      Note $ (b_0, b_1, \cdots , b_n) \neq (3, 2, \cdots , 2) $,
      because if $ (b_0, b_1, \cdots , b_n) = (3, 2, \cdots , 2) $
      then $ L = S^{b_0 b_1 \cdots b_n}_0 = S^{3 2 \cdots 2}_0 $ is
      the last element in $A_n .$

      If $ M = S^{d_0 d_1 \cdots d_n}_0 $, then we get
      $(b_0, b_1, \cdots, b_n) < (d_0, d_1, \cdots, d_n)$ and so

      $$
        \text{ either }
        (
          d_0 = b_0 +1, b_1 = \cdots = b_n = 2
          \text{ and }  d_1 = \cdots = d_n = 1
        )
      $$

      or

      $$
        \exists i_0  \text{ with }  1 \leq i_0 \leq n
        \text{ such that }
      $$

      \begin{align*}
        d_i \:  &= b_i \text{ for all } 0 \leq i < i_0  \\
        d_{i_0} &= 2 , b_{i_0} = 1                      \\
        d_i \:  &= 1, b_i = 2    \text{ for all } i_0 < i \leq n
        \text{ if } 1 \leq i_0 < n   \qquad .
      \end{align*}

      In the first possibility, the side of $L$ and the side of $M,$
      both of which faces the boundary, meet at one point where
      $ S^{b_0 b_1}_0 = S^{b_0 2}_0 $ and
      $ S^{d_0 d_1}_0 = S^{d_0 1}_0 $ meet.

      In the second possibility, the side of $L$ and the side of $M$,
      both of which faces the boundary, meet at one point which is
      contained in such line segment as the intersection of
      $
        S^{b_0 b_1 \cdots b_{i_0 -1} b_{i_0}}_0
        = S^{b_0 b_1 \cdots b_{i_0 -1} 1}_0
      $
      and
      $
        S^{b_0 b_1 \cdots b_{i_0 -1} 2}_0
        = S^{d_0 d_1 \cdots d_{i_0 -1} d_{i_0}}_0  .
      $

      In any possibilities, the side of $L$ and the side of $M$, both
      of which faces the boundary, meet at one point. Since $n$ is
      even, the above condition implies that the point is the ending
      point of the side of $L$, facing the boundary and the starting
      point of the side of $M$, facing the boundary, so

      $$
        \text{the ending point of } L
        \text{ is the starting point of } M.
      $$

    \bigskip

    Subcase 2-3 ) $ m \geq 1 $ and $ c_m \neq 6 $

      \bigskip
      If $m=1$, then the barycenter is both the ending point of
      $L$ and the starting point of $M$ from the definition.

      Assume $m \geq 2.$
      Note
      $L = S^{b_0 b_1 \cdots b_k}_{c_0 c_1 c_2 \cdots c_{m-1} c_m}$
      and its next element $M$ are inside the triangle
      $ S^{b_0 b_1 \cdots b_k}_{c_0 c_1} ,$
      which is one of the triangles obtained by the barycentric
      subdivision
      $ S^{b_0 b_1 \cdots b_k}_{c_0} = S^{b_0 b_1 \cdots b_k}_0 .$

      Compare it with the proper one of $T_0, T_{01}, \tilde{T}_0 $
      and $\tilde{T}_{01}$.  By referring to subcase 1-1,
      - by restricting it to the first triangle if needed-,
      we get

      $$
        \text{the ending point of } L
        \text{ is the starting point of } M.
      $$

    \bigskip

    Subcase 2-4 ) $ m \geq 1 $ and $ c_m = 6 $

      \bigskip

      Note $L = S^{b_0 b_1 \cdots b_k}_{c_0 c_1 \cdots c_m}$ is inside
      the triangle $ S^{b_0 b_1 \cdots b_k}_0 $, where $c_0=0.$

      \medskip
      Assume $ c_1 = c_2= \cdots = c_{m-1}= c_m = 6 $.

      If $m=1,$ then the ending point of
      $
        L = S^{b_0 b_1 \cdots b_k}_{c_0 c_m}
          = S^{b_0 b_1 \cdots b_k}_{0 6}
      $
      will be the ending point of $S^{b_0 \cdots b_k}_{06}$
      tautologically.
      If $m \geq 2,$ then compare $ S^{b_0 b_1 \cdots b_k}_{c_0 c_1} $
      with the proper one of
      $T_0, T_{01}, \tilde{T}_0 $ and $\tilde{T}_{01}$.
      Then from the comparison,
      the ending point of L will be the ending point of
      $S^{b_0 \cdots b_k}_{0 6 6},$ which is also the ending point
      of $S^{b_0 \cdots b_k}_{0 6}.$

      If $M = S^{d_0 d_1 \cdots d_s}_{a_0 \cdots a_t}$ with $ a_0=0$,
      then we get

      \noindent
      either
      $$
        ( k=s, a_1 = \cdots = a_t = 1 \text { and }
        S^{d_0 d_1 \cdots d_s}_{01} \text{ is the next element of }
        S^{b_0 \cdots b_k}_{06}  \text{ in } A_{k+1} )
      $$
      or
      $$
        (  s=k+1, \: t=m-1, \: a_i = 1 \text{ for } 1 \leq i \leq t
        \text{ in case of } m \geq 2 \text{ and }
      $$
      $$
        S^{d_0 d_1 \cdots d_s}_{0} \text{ is the next element of }
        S^{b_0 \cdots b_k}_{06}  \text{ in } A_{k+1} ) .
      $$

      In the first possibility, $ S^{d_0 d_1 \cdots d_s}_0 $ will be
      also the next element of $ S^{b_0  b_1 \cdots b_k}_0 $ in $A_k$
      and so the ending point of $ S^{b_0 b_1 \cdots b_k}_0 $ will be
      the starting point of $ S^{d_0 d_1 \cdots d_s}_0 $ from subcase
      2-1 and 2-2 , which implies

      $$
        \text{the ending point of } S^{b_0 b_1 \cdots b_k}_{06}
        \text{ will be the starting point of }
        S^{d_0 d_1 \cdots d_s}_{01}.
      $$

      \medskip

      In the second possibility, note one of $k$ and $s$ is odd and
      the other one is even, which implies that
      $S^{b_0 b_1 \cdots b_k}_{0}$ is the last element in $A_k$ and
      that $S^{d_0 d_1 \cdots d_s}_{0}$ is the first element in the
      subset

      \begin{align*}
        \{
           S ^{ x_0 x_1 \cdots x_s }_{0} \mid
           x_0 \in \{ 1,2,3 \} \; , \; x_i \in \{ 1,2 \} \textrm{ for }
           i=1, \cdots s
        \}
      \end{align*}

      \noindent
      of $A_s = A_{k+1}$. Also, notice that the ending point of
      $S^{b_0 b_1 \cdots b_k}_{06}$ is also the ending point of
      $S^{b_0 b_1 \cdots b_k}_{0}$ from the definition of triangles.
      By thinking of the side of $S^{b_0 b_1 \cdots b_k}_{0} ,$
      which faces the boundary, and the side of
      $S^{d_0 d_1 \cdots d_s}_{0},$ which faces the boundary, we get

      $$
        \text{the ending point of } S^{b_0 b_1 \cdots b_k}_{0}
        \text{ will be the starting point of }
        S^{d_0 d_1 \cdots d_s}_{0},
      $$

      so

      $$
        \text{the ending point of } S^{b_0 b_1 \cdots b_k}_{06}
        \text{ will be the starting point of }
        S^{d_0 d_1 \cdots d_s}_{0}.
      $$

      In any possibilities, the ending point of $ S^{b_0 \cdots
      b_k}_{06} $, which is also the ending point of $L$, is the
      starting point of its next element in $ A_{k+1} $, which will
      be the starting point of $M$ from $ a_1 = \cdots = a_t = 1 $ if
      $ t \geq 1 $. Thus, we get

      $$
        \text{the ending point of } L
        \text{ is the starting point of } M.
      $$

      \medskip
      Now , assume $ m \geq 2 $ and $ c_i \neq 6$ for some
      $ 1 \leq i < m $.
      From the comparison of $ S^{b_0 b_1 \cdots b_k}_{c_0 c_1} $
      with the proper one of
      $T_0, T_{01}, \tilde{T}_0 $ and $\tilde{T}_{01}$,
      we get $L$ and $M$ are
      inside the triangle $S^{b_0 \cdots b_k}_{c_0 c_1}$
      and from subcase 1-2, we get
      $$
        \text{the ending point of } L
        \text{ is the starting point of } M.
      $$

  \bigskip

\section{About Curves}\label{sec:curves}

  \subsection{Notations}

  $f*g : [0,1] \rightarrow \mathbb{H}^n$ is an ordinary
  juxtaposition of curves $f,g : [0,1] \rightarrow \mathbb{H}^n$.
  And, for a given curve $c:[0,1] \rightarrow \mathbb{H}^n$, $\bar{c}$
  represents a curve whose direction is opposite to that of c, that is,
  $\bar{c}:[0,1] \rightarrow \mathbb{H}^n$ is given by
  $\bar{c}(t) = c(1-t).$

  \subsection{Simplification $ \mathbf{\gamma} $ of a curve
  $ \mathbf{g : [a,b] \rightarrow \mathbb{H}^n}$}

  Given a curve $ g : [a,b] \rightarrow S $ , we can think of a
  curve $ \gamma : [a,b] \rightarrow S $ whose direction is
  one-sided as follows :

  If we can find $ c, d, e \in (a,b) $ such that $ a<c<d<e<b $ and
  $ Im (g| _{[c,d]}) = Im (g| _{[d,e]}) $ and that the directions of
  $ g| _{[c,d]} $ and $ g| _{[d,e]} $ are one-sided but opposite
  from each other, then we can think of the new curve
  $ \tilde{g} : [a,b] \rightarrow D^2 $ from the remaining part
  $ g| _{[a,c]} $ and $ g| _{[e,b]} $ by translating in the domain
  and rescaling as follows :

  Note $g(c) = g(e).$

  Consider two curves $g_1 : [a, d] \rightarrow \mathbb{H}^n  $ and
  $g_2 : [d, b] \rightarrow \mathbb{H}^n $ given by
  $$
    \displaystyle
    g \left( \frac{c-a}{d-a} (t-a) + a \right)
    = g_1 (t)  \text{ for } t \in [a, d]
  $$
  and
  $$
    \displaystyle
    g \left( \frac{b-e}{b-d} (t-b) + b \right)
    = g_2 (t)  \text{ for } t \in [d, b],
  $$
  and then let $\tilde{g} = g_1 * g_2.$

  From a curve obtained by doing this work again and again, we can
  think of a constant speed curve
  $ \gamma : [a,b] \rightarrow S $ which we want.

  \subsection{The definition of $ D_n, j_n, t^n_1, t^n_2 $}

  \begin{align*}
        D_n =
        \displaystyle
        \left\{
           \dfrac{1}{2} \cdot \frac{j}{6^n} \mid j
           = 0, 1, 2, \cdots , 6^n \right
        \}
        \bigcup
        \\
        \left(
          \cup ^n_{k=1}
          \left\{
            \sum ^k_{i=1} \frac{1}{2^i}
            + \frac{1}{2^{k+1}} \cdot
            \frac{j}{2^{k-1} \cdot 6^{n-k+1}}
            \mid
            j=0,1,2, \cdots , 2^{k-1} \cdot 6^{n-k+1}
          \right\}
        \right)
  \end{align*}

  Think of the usual order $D_n$ and regard
  $$
    0, \: \frac{1}{2} \cdot \frac{1}{6^n}, \:
    \frac{1}{2} \cdot \frac{2}{6^n}, \:
    \cdots , \:
    \frac{1}{2} = \frac{1}{2} \cdot \frac{6^n}{6^n}, \:
    \frac{1}{2} + \frac{1}{2^2} \cdot \frac{1}{2^0 \cdot 6^n}, \:
    \cdots \quad
    \in D_n
  $$

  as 0th, 1st, 2nd, $ \cdots , \: 6^n $th, $6^{n+1}$th, $\cdots$
  element, respectively.

  Now, define functions
  $$
    j_n : D_n \rightarrow \{ 0, 1,2,3, \cdots \}
  $$
  $$
    t^n_1 : D_n - \{ 0 \} \rightarrow D_n
  $$
  $$
    t^n_2 :
    D_n - \{ \text{ the last element of } D_n  \} \rightarrow D_n
  $$

  as follows : \\

  $j_n(s) = j$ \quad for the $j$-th element $s \in D_n.$ \\

  $ t^n_1 (s) $ is the $(j-1)$-th element in $ D_n  $ for
  a given $j$-th element $ s \in D_n -\{ 0 \}.$\\

  $ t^n_2 (s) $ is the $(j+1)$-th element in $ D_n $ for a given $j$-th
  element $ s \in D_n -  \{ $ the last element of $ D_n  \} .$ \\

  \bigskip

  \subsection{Definition of $ \mathbf{\gamma ^n_{t_0} , c^n_{t_0} ,
  \bar{c}^n_{t_0}, {_1}c^n_{t_0}, {_1}\bar{c}^n_{t_0},
  \varphi ^n_{t_0}} $ and $\mathbf{\psi ^n_{t_0}}$ on the disk $ \mathbf{D^2} $}

  Let $ n \in \{1,2,3, \cdots \} $ and $ t_0 \in D_n $ be given.
  With respect to the ordering of $D_n$, we'll define
  $ \gamma ^n_{t_0} , c^n_{t_0} , \bar{c}^n_{t_0} $ and
  $ \varphi ^n_{t_0} $ inductively for each fixed n:

  Case 1) $ t_0 $ is the first element in $ D_n $ , in fact,
          $ t_0 = \frac{1}{2} \cdot \frac{1}{{6^n}} $

    \bigskip

    The orientation at the barycenter of $ T_0 \in A_0 $ will give the
    direction of the boundary curve of the first triangle in $A_n.$

    Then
    $$ c^n_{t_0} : [0, 1] \rightarrow \{ basepoint \} \subset D^2 $$

    $$
      \bar{c} ^n_{t_0} : [ 0, 1] \rightarrow
      \{ basepoint \} \subset D^2
    $$

    $$ \varphi ^n_{t_0} : [ 0, 1] \rightarrow D^2 $$

    and

    $$ \gamma ^n_{t_0} : [0, 1 ] \rightarrow D^2 $$

    \noindent
    can be thought, where $ \varphi ^n_{t_0} $ and
    $ \gamma ^n_{t_0} $ are the piecewise smooth boundary curve of
    the first triangle in $A_n$ with constant speed
    and the direction of the boundary curve is induced from the given
    orientation.

    Note $ \gamma ^n_{t_0} $ can be regarded as the simplification of  $ c^n_{t_0} * \varphi ^n_{t_0} *
    \bar{c} ^n_{t_0} \: .$

    We will call $ \gamma ^n_{t_0} $ the \emph{holonomy curve at time
    $ t = t_0 $}.

    Now,  consider the path from the basepoint to the ending point of
    the first triangle in $n$-step along the opposite
    direction of the holonomy curve $ \gamma ^n_{t_0} $ at $t={t_0}$
    , which is a piecewise smooth curve with constant speed. Then from
    the path, we can define a piecewise smooth curve

    $$ _1 c ^n_{t_0} : [0,1] \rightarrow D^2 $$

    \noindent
    with constant speed.
    And its opposite direction can make us define

    $$ _1 \bar{c} ^n_{t_0} : [0,1] \rightarrow D^2 .$$

    Define a piecewise smooth curve

    $$ \psi ^n_{t_0} : [0,1] \rightarrow D^2 $$

    \noindent
    with constant speed as the boundary curve of the 1st triangle in
    the $n$-th step, where the curve is a loop at the ending point of
    the first triangle and the direction of the boundary curve is
    induced from the given orientation.

  \bigskip

  Case 2) $ t_0 $ is the $j$-th element in $ D_n $ ,
          where $ j \geq 2 $

    \bigskip

    Let $ t_1 $ be the $(j-1)$-th element in $ D_n $ ,
    where $ j-1 \geq 1 .$

    Consider the path from the basepoint  to the starting point of
    the $j$-th triangle in the $n$-th step along the opposite
    direction of the holonomy curve $ \gamma ^n_{t_1} $ at
    $ t = t_1 $ , which is a piecewise smooth one with constant speed
    .  Then from the path, we can define a piecewise smooth curve

    $$ c^n_{t_0} : [0, 1] \rightarrow  \partial U_{j-1} \subset D^2 $$

    \noindent  with constant speed, where $U_{j-1}$ is the union
    of triangle in $A_n$ from 1st one to $(j-1)$-th one.

    And its opposite direction can make us define

    $$ \bar{c} ^n_{t_0} : [0 , 1] \rightarrow  \partial U_{j-1} \subset D^2 .$$

    Define a piecewise smooth curve

    $$ \varphi ^n_{t_0} : [0,1] \rightarrow D^2 $$

    \noindent
    with constant speed as the boundary curve of the $j$-th
    triangle in the $n$-th step, where the curve is a loop at the
    starting point of the triangle and  the direction of the boundary
    curve is induced from the given orientation.

    Now define a piecewise smooth curve

    $$ \gamma ^n_{t_0} : [0, 1 ] \rightarrow  \partial U_{j} \subset D^2 $$

    \noindent
    with constant speed from the simplification of
    $
      \gamma ^n_{t_1} * c^n_{t_0} * \varphi ^n_{t_0} * \bar{c}
      ^n_{t_0},
    $
    where $U_j$ is the union
    of triangle in $A_n$ from 1st one to $j$-th one.
    The new curve will be also
    called the \emph{holonomy curve at time $ t=t_0 $ }.


    Now,  consider the path from the basepoint to the ending point of
    the $j$-th triangle in the $n$-th step along the opposite
    direction of the holonomy curve $ \gamma ^n_{t_0} $ at $t={t_0}$
    , which is a piecewise smooth one with constant speed. Then from
    the path, we can define a piecewise smooth curve

    $$ _1 c ^n_{t_0} : [0,1] \rightarrow \partial U_j \subset D^2 $$

    \noindent
    with constant speed.
    And its opposite direction can make us define

    $$ _1 \bar{c} ^n_{t_0} : [0,1] \rightarrow \partial U_j \subset D^2 .$$

    Define a piecewise smooth curve

    $$ \psi ^n_{t_0} : [0,1] \rightarrow D^2 $$

    \noindent
    with constant speed as the boundary curve of the $j$-th triangle in
    the $n$-th step, where the curve is a loop at the ending point of the
    $j$-th triangle and the direction of the boundary curve is induced
    from the given orientation.

  \subsection{
               the simplification of
               $\mathbf{\bar{c}^n_{t_0} * {_1 c^n_{t_0}}}$
             }

  \bigskip
%
%

  For each $n \geq 1$ and $t_0 \neq 0$, where  $t_0$ is the
  $j_n (t_0)$-th element in $D_n$, the simplification of
  $\bar{c}^n_{t_0} * {_1 c^n_{t_0}}$ is a curve along the boundary
  curve  of $j_n (t_0)$-th triangle in $A_n$ with opposite
  direction to the given orientation such that it
  starts from the starting point of the triangle and that its image
  consists of the following sets :

  one point, one side, two sides or the boundary of the triangle.

  \medskip

  Proof )

  If n=1, then it can be easily checked.

  Assume $ n \geq 2 $. If $t_0 $ is greater than the maximum of
  $D_{n-1} - \{ 1 \}$, then the above property can be easily checked.

  Now assume $ n \geq 2 $ and $t_0 $ is less than or equal to the
  maximum of $D_{n-1} - \{ 1 \}$. Now find $\delta _n (t_0) \in D_{n-1}$
  such that $t^{n-1}_1 ( \delta _n (t_0) ) < t_0 \leq \delta _n (t_0)$,
  where $t^{n-1}_1 ( \delta _n (t_0) )$ is the previous element of
  $\delta _n (t_0)$ in $D_{n-1}.$ Then, the $j_n(t_0)$-th triangle
  in $A_n$ is one of the barycentric subdivision of the
  $j_{n-1}(\delta _n (t_0))$-th triangle in $A_{n-1}$.

  And find a value $\epsilon ( j_n (t_0) ) $ such that , for the
  given $ L = j_n (t_0)$-th triangle in $A_n$,

  $$
    \text{ if } L = T_{a_0 a_1 \cdots a_n}
    \text{, then }  \epsilon ( j_n (t_0) ) = a_n
  $$
  and
  $$
    \text{ if } L = S^{b_0 b_1 \cdots b_k}_{a_0 a_1 \cdots a_s}
    \text{, where } n = k+s,
    \text{ then }  \epsilon ( j_n (t_0) ) = a_s .
  $$

  Assume $n \geq 2$ and that the property, mentioned early in this
  subsection, holds for $n-1$.

  Then we obtained the following result.

  Case 1) Assume the image of the simplification of
         $
           \bar{c}^{n-1}_{ \delta _n (t_0) } *
           {_1 c^{n-1}_{ \delta _n (t_0)}}
         $
  consists of one point.

    Now refer to the following picture under the counterclockwise
    orientation. The thick line is a part of the image of
    $\gamma ^{n-1}_{t^{n-1}_1 (\delta_n (t_0))}$
    and the outer triangle is the $j_{n-1} (\delta _n (t_0))$-th
    triangle in $A_{n-1}.$


    \begin{figure}[h]
     \centering
       {
         \includegraphics[width=1.7cm]{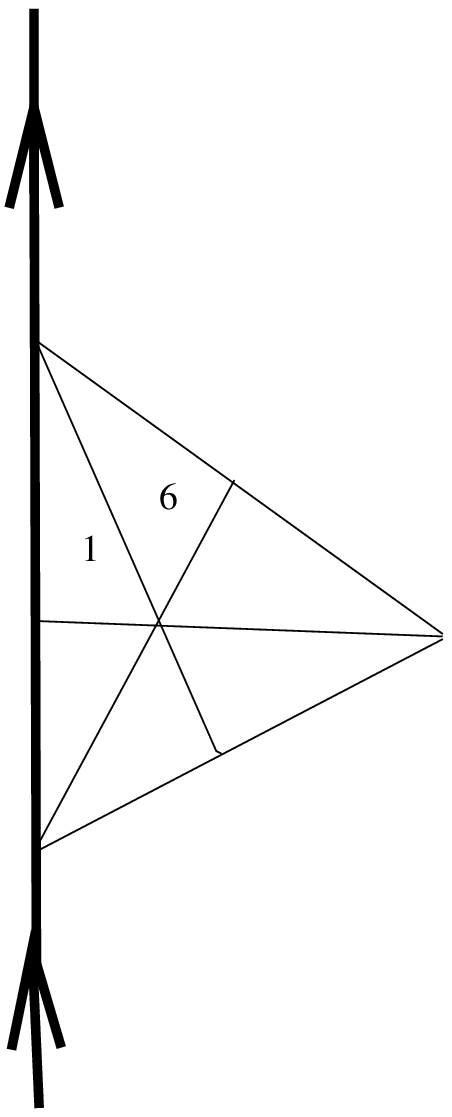}
         \hspace{1cm}
         \includegraphics[width=1.7cm]{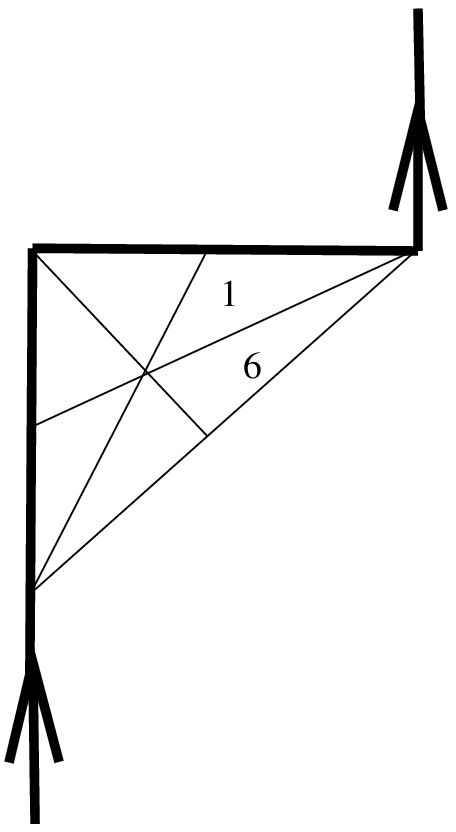}
       }
    \end{figure}

    Note the direction of the line segment of the
    $j_{n-1} (\delta_n (t_0)) $-th triangle along
    $\gamma ^{n-1}_{t^{n-1}_1 (\delta_n (t_0))}$, mentioned in the
    Property 1 in the subsection ~\ref{sec:prop-of-triangles} of the
    section ~\ref{sec:triangles},
    lying on the boundary curve
    $\gamma ^{n-1}_{t^{n-1} _1(\delta_n (t_0))}$, is
    from the common vertex of
    the $j_{n-1} (\delta _n (t_0))$-th triangle in $A_{n-1}$
    with the second triangle of its barycentric subdivision
    to its common vertex with the first triangle of its
    barycentric subdivision , and

    \begin{align*}
      \epsilon ( j_n (t_0) ) &= 1,6 \Rightarrow
      \bar{c}^n_{t_0} * {_1 c^n_{t_0}} \text{ consists of one side}
      \\
      \epsilon ( j_n (t_0) ) &= 2,3 \Rightarrow
      \bar{c}^n_{t_0} * {_1 c^n_{t_0}} \text{ consists of one point}
      \\
      \epsilon ( j_n (t_0) ) &= 4,5 \Rightarrow
      \bar{c}^n_{t_0} * {_1 c^n_{t_0}} \text{ consists of one point}
    \end{align*}

    \medskip





  Case 2) Assume the image of the simplification of
          $
            \bar{c}^{n-1}_{ \delta _n (t_0) } *
            {_1 c^{n-1}_{ \delta _n (t_0)}}
          $
  consists of one side.

    Now refer to the following picture under the counterclockwise
    orientation. The thick line is a part of the image of
    $\gamma ^{n-1}_{t^{n-1}_1(\delta_n (t_0))}$
    and the outer triangle is the $j_{n-1} (\delta _n (t_0))$-th
    triangle in $A_{n-1}.$ Don't forget that
    the ending point of $j_{n-1}(\delta_{n} (t_0))$ in $A_{n-1}$ will
    lie on the image of $\gamma ^{n-1}_{\delta_n (t_0)},$
    even though it might not lie on the image of
    $\gamma ^{n-1}_{t^{n-1} _1 (\delta_n (t_0))}$.

    \begin{figure}[h]
      \centering
        {
          \includegraphics[width=2cm]{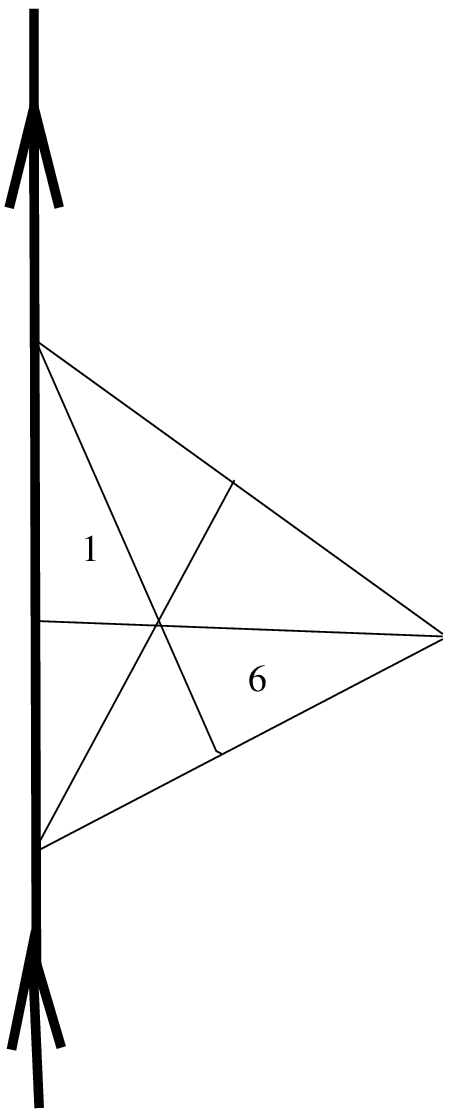} \:
          \includegraphics[width=2cm]{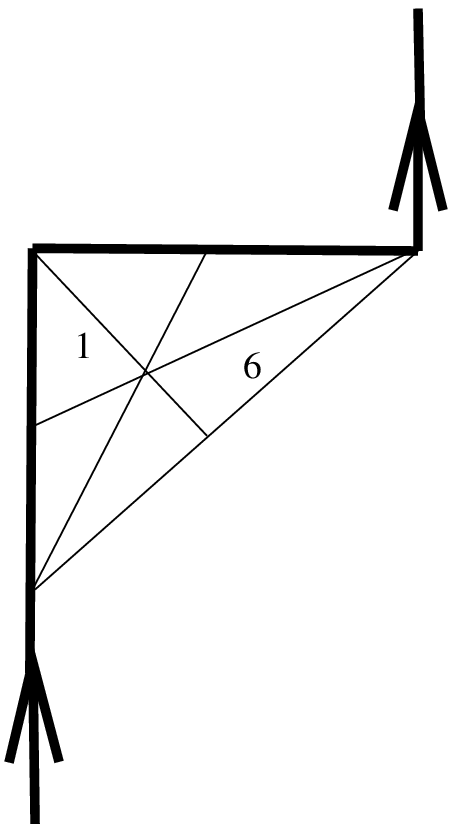} \:
          \includegraphics[width=2cm]{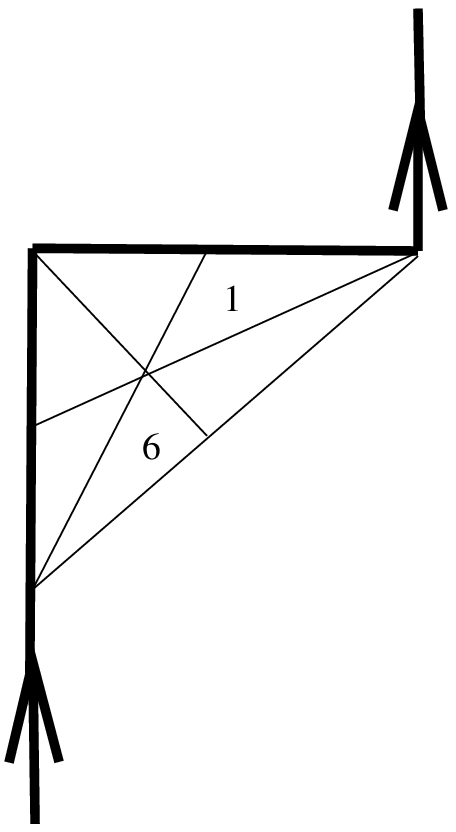}
        }
      \\
      \centering
        {
          \includegraphics[width=2cm]{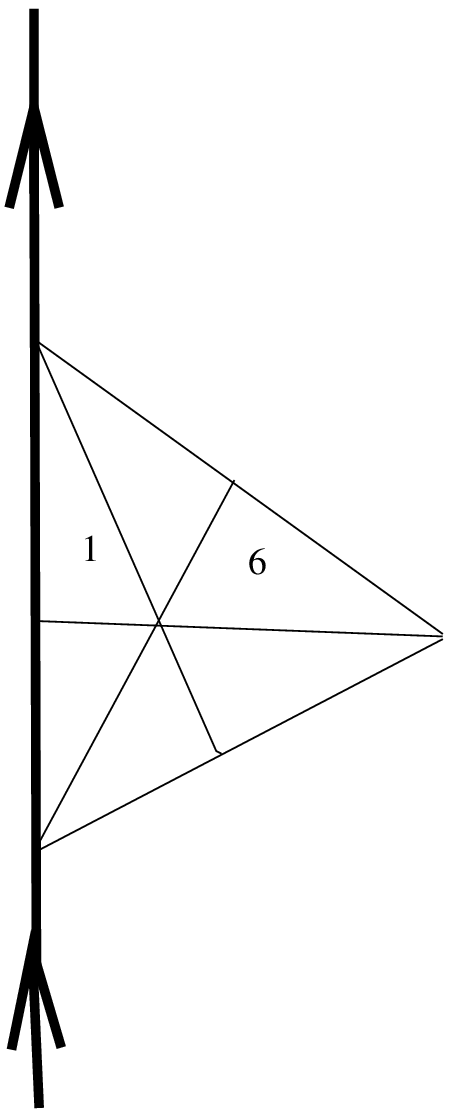} \:
          \includegraphics[width=2cm]{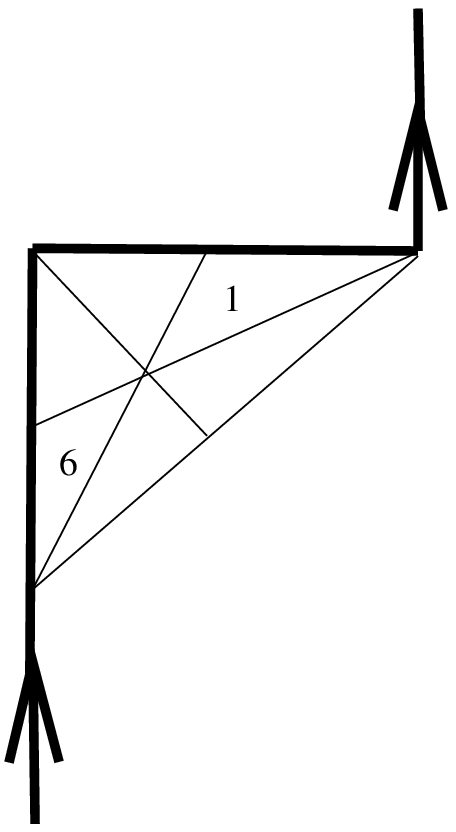}
        }
    \end{figure}

    Note the direction of the line segment of the
    $j_{n-1} (\delta_n (t_0)) $-th triangle along
    $\gamma ^{n-1}_{t^{n-1}_1 (\delta_n (t_0))}$, mentioned in the
    Property 1 in the subsection ~\ref{sec:prop-of-triangles} of the
    section ~\ref{sec:triangles},
    lying on the boundary curve
    $\gamma ^{n-1}_{t^{n-1} _1(\delta_n (t_0))}$, is
    from the common vertex of
    the $j_{n-1} (\delta _n (t_0))$-th triangle in $A_{n-1}$
    with the second triangle of its barycentric subdivision
    to its common vertex with the first triangle of its
    barycentric subdivision , and

    \begin{align*}
      \epsilon ( j_n (t_0) ) &= 1 \Rightarrow
      \bar{c}^n_{t_0} * {_1 c^n_{t_0}}
      \text{ consists of one side}
      \\
      \epsilon ( j_n (t_0) ) &= 2,3 \Rightarrow
      \bar{c}^n_{t_0} * {_1 c^n_{t_0}}
      \text{ consists of one point}
      \\
      \epsilon ( j_n (t_0) ) &= 4 \Rightarrow
      \bar{c}^n_{t_0} * {_1 c^n_{t_0}}
      \text{ consists of the boundary}
      \\
      \epsilon ( j_n (t_0) ) &= 5 \Rightarrow
      \bar{c}^n_{t_0} * {_1 c^n_{t_0}}
      \text{ consists of either the boundary or one side}
      \\
      \epsilon ( j_n (t_0) ) &= 6 \Rightarrow
      \bar{c}^n_{t_0} * {_1 c^n_{t_0}}
      \text{ consists of either one side or two sides}
    \end{align*}

    \begin{remark}
      The last 2 pictures in the bottom seem to be possible under the
      induction hypothesis. But it might not happen in fact.
    \end{remark}

    \begin{remark}
      The following picture in the bottom can't happen
      from Property 2 in the subsection ~\ref{sec:prop-of-triangles} of
      the section ~\ref{sec:triangles}.
    \end{remark}
    \begin{figure}[h]
      \centering
        {
          \includegraphics[width=2.5cm]{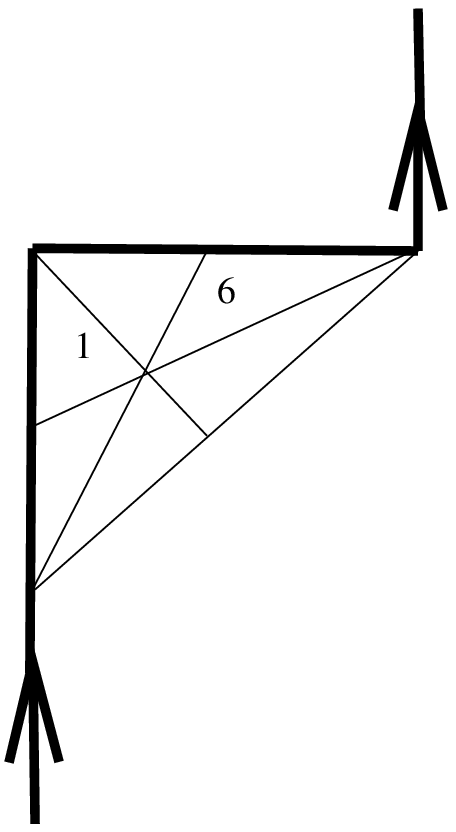}
        }
    \end{figure}

  Case 3) Assume the image of the simplification of
          $
            \bar{c}^{n-1}_{ \delta _n (t_0) } *
            {_1 c^{n-1}_{ \delta _n (t_0)}}
          $
  consists of two sides.

    Now refer to the following picture under the counterclockwise
    orientation. The thick line is a part of the image of
    $\gamma ^{n-1}_{\delta_n (t_0)}$
    and the outer triangle is the $j_{n-1} (\delta _n (t_0))$-th
    triangle in $A_{n-1}.$ Don't forget that
    the ending point of $j_{n-1}(\delta_{n} (t_0))$ in $A_{n-1}$ will
    lie on the image of $\gamma ^{n-1}_{\delta_n (t_0)},$
    even though it might not lie on the image of
    $\gamma ^{n-1}_ {t^{n-1} _1 (\delta _n (t_0))}$.
   \pagebreak
    \begin{figure}[h]
      \centering
        {
          \includegraphics[width=2cm]{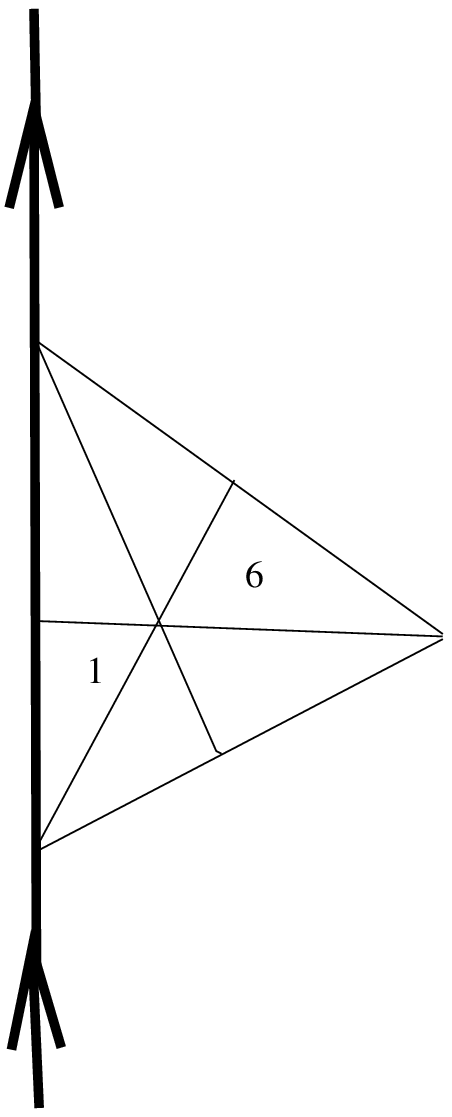} \:
          \includegraphics[width=2cm]{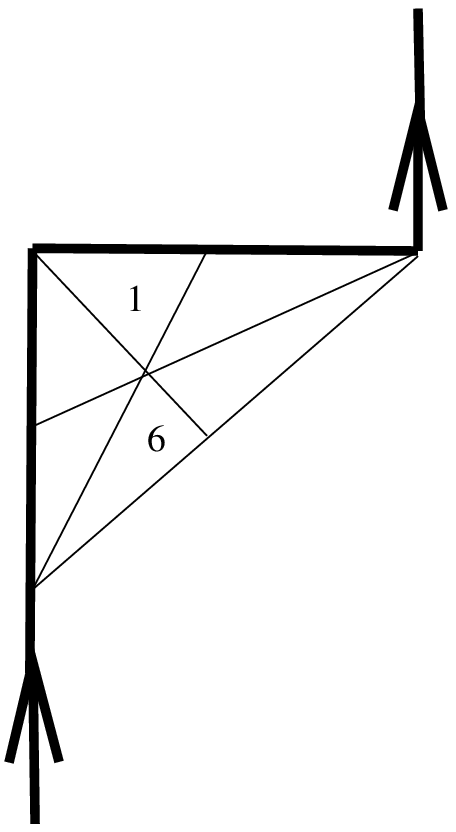} \:
          \includegraphics[width=2cm]{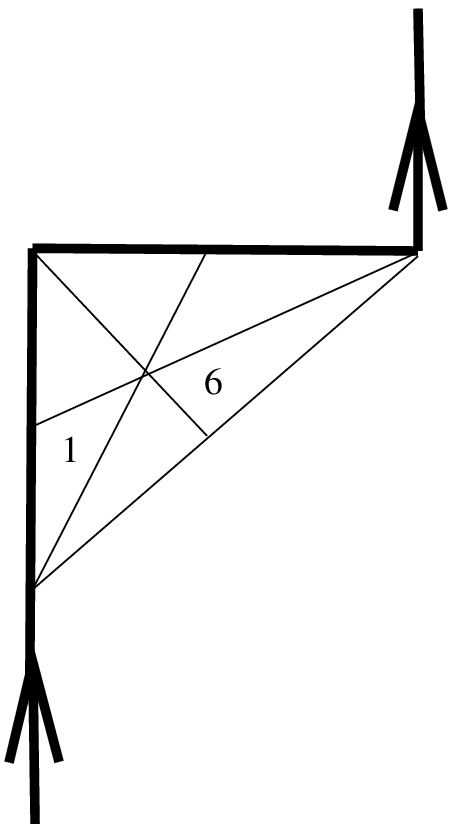}
        }
     \\
      \centering
        {
          \includegraphics[width=2cm]{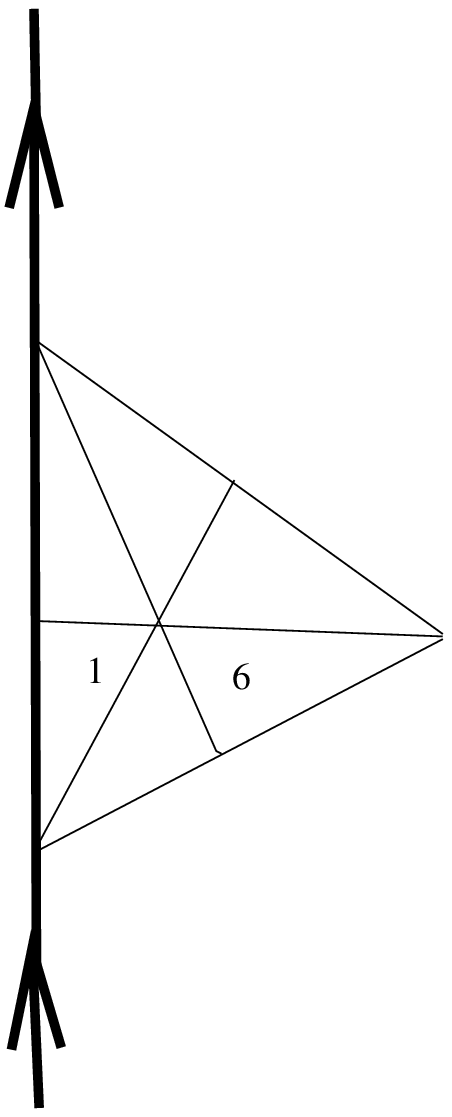} \:
          \includegraphics[width=2cm]{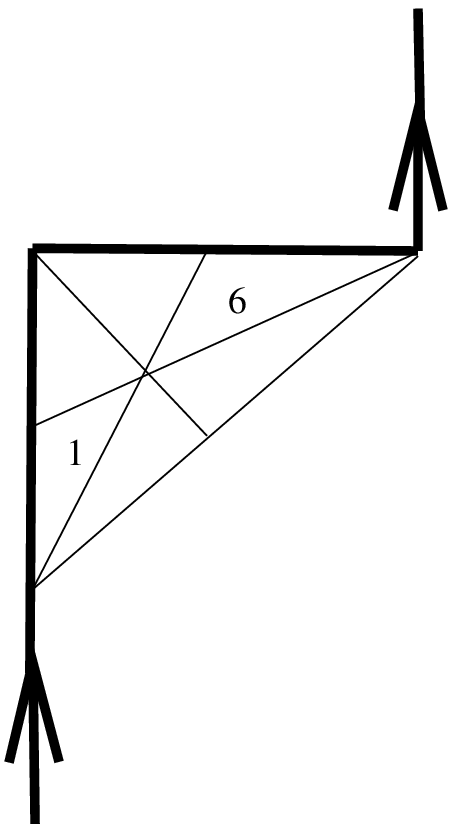}
        }
    \end{figure}

    Note the direction of the line segment of the
    $j_{n-1} (\delta_n (t_0)) $-th triangle along
    $\gamma ^{n-1}_{t^{n-1}_1 (\delta_n (t_0))}$, mentioned in the
    Property 1 in the subsection ~\ref{sec:prop-of-triangles} of the
    section ~\ref{sec:triangles},
    lying on the boundary curve
    $\gamma ^{n-1}_{t^{n-1} _1(\delta_n (t_0))}$, is
    from the common vertex of
    the $j_{n-1} (\delta _n (t_0))$-th triangle in $A_{n-1}$
    with the first triangle of its barycentric subdivision
    to its common vertex with the second triangle of its
    barycentric subdivision , and

    \begin{align*}
      \epsilon ( j_n (t_0) ) &= 1 \Rightarrow
      \bar{c}^n_{t_0} * {_1 c^n_{t_0}}
      \text{ consists of two sides}
      \\
      \epsilon ( j_n (t_0) ) &= 2,3 \Rightarrow
      \bar{c}^n_{t_0} * {_1 c^n_{t_0}}
      \text{ consists of the boundary}
      \\
      \epsilon ( j_n (t_0) ) &= 4 \Rightarrow
      \bar{c}^n_{t_0} * {_1 c^n_{t_0}}
      \text{ consists of one point}
      \\
      \epsilon ( j_n (t_0) ) &= 5 \Rightarrow
      \bar{c}^n_{t_0} * {_1 c^n_{t_0}}
      \text{ consists of either one point or the boundary}
      \\
      \epsilon ( j_n (t_0) ) &= 6 \Rightarrow
      \bar{c}^n_{t_0} * {_1 c^n_{t_0}}
      \text{ consists of either two sides or one side}
    \end{align*}

    \begin{remark}
      The last 3 pictures in the bottom seem to be possible under the
      induction hypothesis. But it might not happen in fact.
    \end{remark}

    \begin{remark}
      The following picture in the bottom can't happen
      from Property 2 in the subsection ~\ref{sec:prop-of-triangles} of
      the section ~\ref{sec:triangles}.
    \end{remark}
    \begin{figure}[h]
      \centering
        {
          \includegraphics[width=2cm]{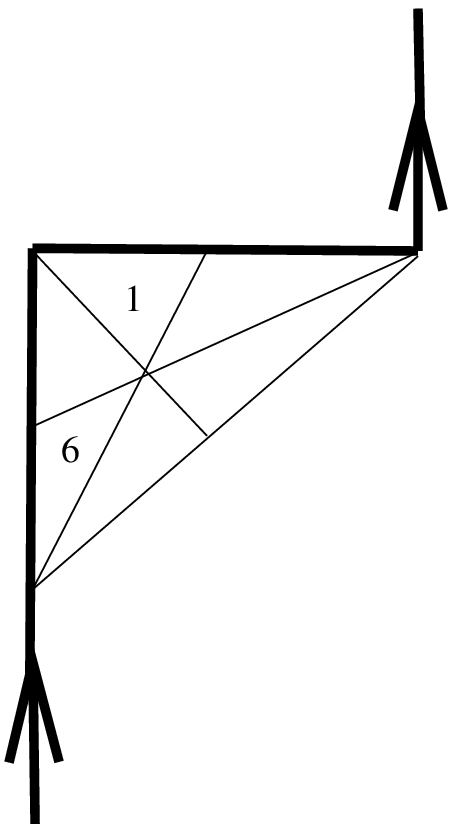}
        }
    \end{figure}

  \pagebreak

  Case 4) Assume the image of the simplification of
          $
            \bar{c}^{n-1}_{ \delta _n (t_0) } *
            {_1 c^{n-1}_{ \delta _n (t_0)}}
          $
  consists of the boudary.

    Now refer to the following picture under the counterclockwise
    orientation. The thick line is a part of the image of
    $\gamma ^{n-1}_{t^{n-1}_1(\delta_n (t_0))} $
    and the outer triangle is the $j_{n-1} (\delta _n (t_0))$-th
    triangle in $A_{n-1}.$

    \begin{figure}[h]
      \centering
        {
          \includegraphics[width=2cm]{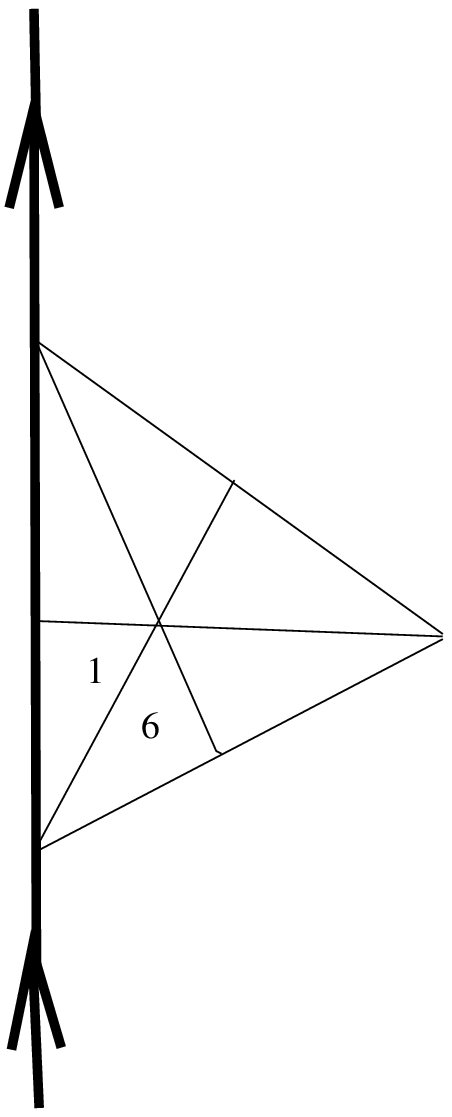}
          \includegraphics[width=2cm]{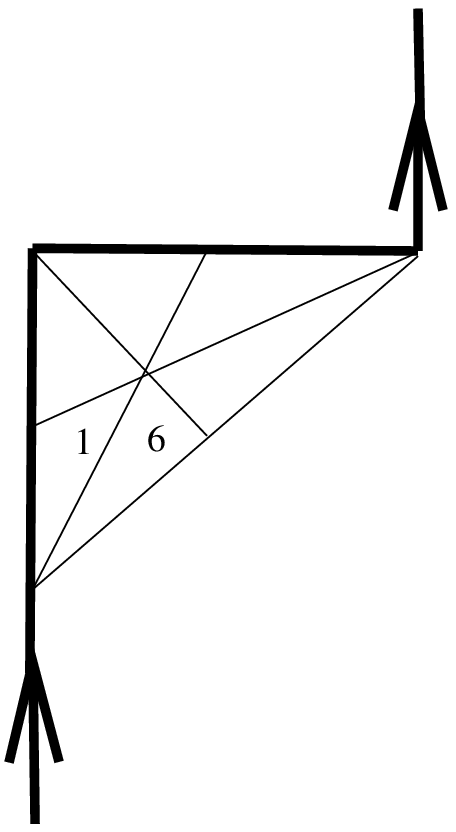}
        }
    \end{figure}

    Note the direction of the line segment of the
    $j_{n-1} (\delta_n (t_0)) $-th triangle along
    $\gamma ^{n-1}_{t^{n-1}_1 (\delta_n (t_0))}$, mentioned in the
    Property 1 in the subsection ~\ref{sec:prop-of-triangles} of the
    section ~\ref{sec:triangles},
    lying on the boundary curve
    $\gamma ^{n-1}_{t^{n-1} _1(\delta_n (t_0))}$, is
    from the common vertex of
    the $j_{n-1} (\delta _n (t_0))$-th triangle in $A_{n-1}$
    with the first triangle of its barycentric subdivision
    to its common vertex with the second triangle of its
    barycentric subdivision , and

    \begin{align*}
      \epsilon ( j_n (t_0) ) &= 1,6 \Rightarrow
      \bar{c}^n_{t_0} * {_1 c^n_{t_0}}
      \text{ consists of two sides}
      \\
      \epsilon ( j_n (t_0) ) &= 2,3 \Rightarrow
      \bar{c}^n_{t_0} * {_1 c^n_{t_0}}
      \text{ consists of the boundary}
      \\
      \epsilon ( j_n (t_0) ) &= 4,5 \Rightarrow
      \bar{c}^n_{t_0} * {_1 c^n_{t_0}}
      \text{ consists of the boundary}
    \end{align*}




  \bigskip
\end{appendices}


\bigskip

\begin{thebibliography}{9}
\bibitem{CL} Y.Choi and K.Lee, \emph{Holonomy of homogeneous
    submersions of} ${\rm PSL}_2 \mathbb{R}$, ???

\bibitem{GW} D.Gromoll and G.Walschap, \emph{Metric foliations and
    curvature}, Progress in Mathematics, 268. Birkh$\ddot{\rm a}$user,
    Verlag, Basel,2009

\bibitem{KN} S.Kobayashi and K.Nomizu, \emph{Foundations of
    differential geometry, volume 1, 2}, Jonh Wiley and Sons, Inc.,
    New York, 1969

\bibitem{W} G.Walschap, \emph{Metric structures in differential
    geometry}, Graduate Texts in Mathematics, 224, Springer-Verlag,
    New York, 2004
\end{thebibliography}
\end{document}